\documentclass{article}
\usepackage{amsfonts}
\usepackage{graphicx}
\usepackage[dvips]{color}
\usepackage{amssymb}
 
\begin{document}
\title{On the Finiteness Problem for Relative Minima of Area Bounded by a Real-Analytic Jordan Curve}
\date{\today}
\author{Michael Beeson\\
Department of Mathematics\\
San Jose State University\\
San Jose, California 95192\\
USA\\
email:  beesonpublic@gmail.com\\
}
\maketitle

\medskip
\newtheorem{theorem}{Theorem}
\newtheorem{lemma}{Lemma}
\newtheorem{corollary}{Corollary}
\newtheorem{definition}{Definition}
\newtheorem{conjecture}{Conjecture}
\def\arclength{{\bigg(\int \AA^2 \SS\, d\w\bigg)}}
\def\Arclength#1{{\bigg(\int_0^{#1} \AA^2 \SS\, d\w\bigg)}}
\def\p{\mathcal P}
\def\q{\mathcal Q}
\def\r{\mathcal R}
\def\ss{\mathcal S}
\def\mm{\mathcal M}
\def\St{{\bf St}}
\def\Ntwo{{}^2N}
\def\Nthree{{}^3N}
\def\None{{}^1N}
\def\uone{{}^1u}
\def\utwo{{}^2u}
\def\uthree{{}^3u}
\def\eq{eq}
\def\AA{{\mathbb A}}
\def\BB{{\mathbb B}}
\def\SS{{\mathbb S}}
\def\HH{{\mathbb H}}
\def\NN{{\mathbb N}}
\def\tAA{{\tilde \AA}}
\def\tBB{{\tilde \BB}}
\def\tSS{{\tilde \SS}}
\def\tHH{{\tilde \HH}}
\def\alphat{{\hat \alpha}}
\def\betat{{\hat \beta}}
\def\zetat{{\hat \zeta}}

\def\grad{{\rm grad\ }}
\def\fprime{f^{\prime}}
\def\gprime{g^{\prime}}
\def\Zprime{Z^{\prime}}
\def\cprime{c^{\prime}}
\def\tA{{\tilde\A}}
\def\tB{{\tilde\B}}
\def\tS{{\tilde\S}}
\def\vector#1#2#3{\left[
\begin{array}{l} 
#1\\ 
#2\\ 
#3 
\end{array} \right]}
\def\threecases#1#2#3{\left\{
\begin{array}{ll} 
#1\\ 
{}\\
#2\\
{}\\ 
#3 
\end{array} \right.}  

\def\vectortwo#1#2{\left[
\begin{array}{l} 
#1\\ 
#2  
\end{array} \right]} 

\def\matrixtwo#1#2#3#4{{
 \left( \begin{array}{cc}
#1 & #2   \\
#3 & #4   
\end{array} \right)}}

\def\twocases#1#2{\left\{
\begin{array}{ll} 
#1\\ 
{}\\
#2
\end{array} \right.}  

\def\fourcases#1#2#3#4{\left\{
\begin{array}{ll} 
#1\\ 
{}\\
#2\\
{}\\ 
#3\\
{}\\
#4
\end{array} \right.}  

\def\onehalf{\frac{1}{2}}
\def\iover2{\frac{i}{2}}
\def\Re{{\rm Re}}
\def\Im{{\rm Im}}
\def\A{{\bf A}}
\def\B{{\bf B}}
\def\S{{\bf S}}
\def\Q{{\bf Q}}
\def\a{{\bf a}}
\def\b{{\bf b}}
\def\s{{\bf s}}
\def\O{{\bf O}}
\def\R{{\bf R}}
\def\C{{\bf C}}
\def\P{{\mathcal P}}
\def\K{{\mathcal K}}
\def\Rot{{\rm Rot}}
\def\w{{w}}   
\def\degp{{{\rm deg}\ \p}}
\def\degr{{{\rm deg}\ \r}}
\def\degs{{{\rm deg}\ \ss}}

\begin{abstract}
Earlier versions of this paper claimed that a real-analytic Jordan curve cannot bound infinitely many relative minima of area. That claim is not proved. The first error is in the proof of Lemma 24 (Section 10); the material before that lemma is not affected.
\end{abstract}

\section{Introduction}

Minimal surfaces are mathematical objects that are intimately related to the physical surfaces
formed by thin soap films.  Aside from the idealization to zero thickness, there are some 
other differences between physical soap films and minimal surfaces.  First, physical soap films
are {\em stable}, in the sense that if they are disturbed slightly, they regain their shape.
Secondly,  physical soap films sometimes have internal edges where several films meet, or 
other kinds of more complicated topology.  Here we are concerned with the classical 
{\em problem of Plateau},  which requires,  given a Jordan curve $\Gamma$, to find 
(or at least prove the existence of) a surface of the topological type of the disk bounded
by $\Gamma$ and minimizing area among such surfaces.   Surfaces forming a {\em relative minimum}
of area correspond to soap films stable in the physical sense just described.  

Plateau's problem was solved independently by Douglas and Rado in the early 1930s.  Their
solution methods produced minimal surfaces that might have certain 
singularities known as {\em branch points}.
 It was not until the seventies that the regularity (lack of 
branch points) of solutions providing an absolute minimum of area with given real-analytic boundary was proved.
Readers not already intimately familiar with branch points may want to view the animated pictures posted to the 
Web at \cite{pictures} for a visual introduction to branch points.  One picture of a boundary branch point
is given here in Figure 1. 

The ``finiteness problem'' addressed in this paper is to prove that for a given Jordan curve $\Gamma$,
only finitely many disk-type soap films are bounded by $\Gamma$.  The reason that minimal surfaces
are called ``minimal'' is that a soap film furnishes a relative minimum of area; that is, 
perturbing the surfaces slightly will not decrease the area.  Thus the mathematical form of 
the finiteness problem is to prove that a given Jordan curve $\Gamma$ cannot bound
infinitely many relative minima of area.  We prove that theorem in this paper (for real-analytic
boundaries and disk-type surfaces).   

Under those same assumptions,
Tomi proved \cite{tomi2} that there cannot exist infinitely many {\em absolute minima} of area. 
But generally there will exist many relative minima of area that are not absolute minima,
so Tomi's theorem does not rule out infinitely many physical soap films.
Tromba proved \cite{tromba} {\em generic finiteness}; that is, the set of boundary 
curves for which finiteness holds is open and dense in a suitable topology.  Here we prove
a theorem about {\em every} real-analytic Jordan boundary.  
The finiteness problem for relative minima that we solve here
is a long-standing open problem.%
\footnote{See e.g. \cite{nitsche}, pages~252 and 379, for discussion and many references.
Although the date of that book is 1989,  there has been no progress on finiteness since then.
See also \cite{bohme2} for a somewhat more detailed discussion.
}

Branch points of minimal surfaces not furnishing a relative minimum of area are also of considerable 
interest, because they are intimately connected with the finiteness problem.   It follows (as shown 
in \cite{part1}, p. 121)  from the 
work of B\"ohme \cite{bohme} that if a real-analytic Jordan curve $\Gamma$ bounds infinitely many
disk-type relative minima of area, then it bounds an analytic one-parameter family of such
minimal surfaces, all with the same Dirichlet's integral,  which can be continued until it either loops 
or terminates in a branched minimal surface.  
In \cite{tomi2}, it is proved that a loop of disk-type minimal surfaces contains one which is not a relative
minimum of area.
Therefore, to prove that $\Gamma$ 
does not bound infinitely many relative minima of area,  we only have to rule out the possibility that $\Gamma$
bounds a 
one-parameter family of minimal surfaces $u(t,z)$, terminating in a
branched minimal surface when the parameter $t$ is zero, and furnishing a relative minimum of area 
for each positive $t$.  Tomi solved the finiteness problem 
for absolute minima of area and real-analytic boundaries \cite{tomi2} in this way, since the limit of absolute
minima of area is again an absolute minimum, and hence has no branch points.  But the problem for relative minima
has remained open.  (A discussion of other open problems in this area is given at the end of this paper.)  

\begin{figure}[htp]
\centering
\includegraphics[width=0.9\textwidth]{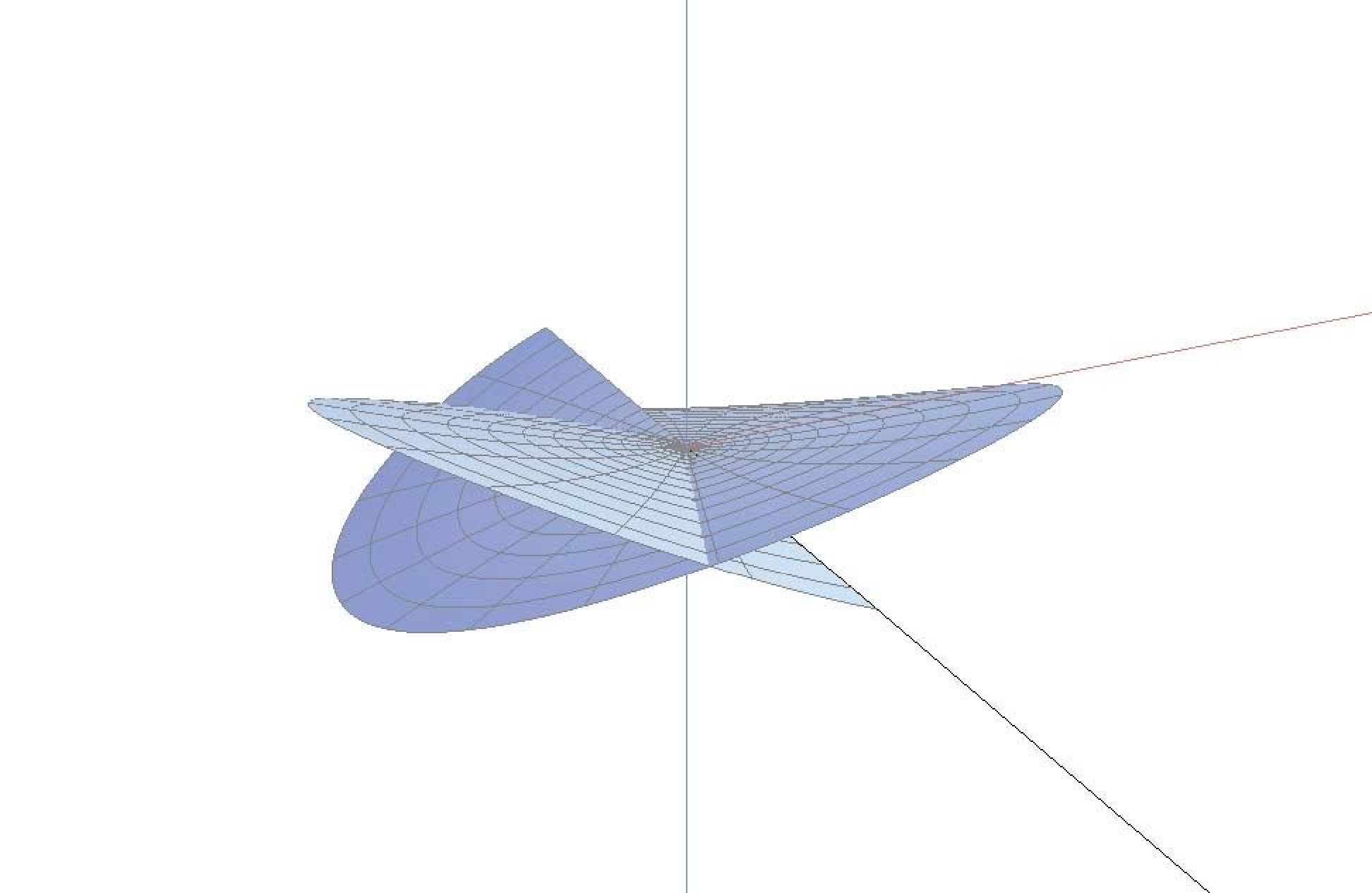}
\caption{\small This is a boundary branch point of order 2 and index 1.  A boundary branch point of order $2m$
goes ``around'' $(2m+1)/2$ times while going ``up and down'' $(2m+k+1)/2$ times.  Could a one-parameter family 
of surfaces, each of which is a relative minima of area, and all with the same Jordan boundary, end in a surface like this?
}
\end{figure}

Suppose we have a boundary branch point;  orient the boundary so it is tangent to the $X$-axis with the branch point
at origin.  Consider the unit normal $N$ restricted to the boundary.  As the parameter $t$ varies, we do not know 
{\em a priori} how $N$ behaves, but at least it must remain perpendicular to the boundary. Close to the branch
point it must almost lie in the $YZ$ plane.  Some of our paper contains 
calculations concerning the possible behavior of $N$ as the parameter $t$ goes to zero.

The finiteness problem for relative minima of area was first attacked in \cite{part1}, 
where interior branch points arising as the limit of one-parameter families of relative minima were 
shown not to exist, and in \cite{part2}, such boundary branch points were eliminated under special
conditions on $\Gamma$, but not in general; so the finiteness problem for relative minima remained open
until now, except for special boundaries.  

It is well known that there is a natural eigenvalue problem associated to the second variation
of the area functional.  That problem, and its first eigenvalue and the associated eigenfunction,
play an important role in our work.
The main line of attack in \cite{part1} is an eigenvalue argument, that in the case of an interior branch point,
the Gaussian image of $u$ for small positive $t$ must include more than a hemisphere, and hence $u$ cannot be 
a relative minimum of area.  The reason why it must include more than a hemisphere is the 
 Gauss-Bonnet-Sasaki-Nitsche formula,  which expresses the total boundary curvature as a sum of contributions
from the Gaussian area and the branch points.  This formula implies that as soon as $t > 0$, the 
contribution of the branch point term to the total curvature when $t=0$ must be made up by ``extra'' Gaussian area
totaling $2m\pi$, where $2m$ is the order of the boundary branch point. For boundary branch points, this argument 
implies that for small positive $t$, there cannot be a neighborhood in the parameter domain on which the Gauss map
covers more than the upper hemisphere.  
Near the boundary the normal is confined to within $O(t)$ of the $YZ$ plane,
and it seems plausible that the extra Gaussian area is contributed in the form of $m$ hemispheres.  However, that is 
difficult to prove.   The behavior of $N$ could, {\em a priori}, be quite complicated.

The main line of attack in this paper is the calculation of the eigenfunction corresponding to the first 
eigenvalue.  The starting point of this calculation is the relation between the second variations of area and of Dirichlet's integral;
namely, the condition for $\phi$ to belong to the 
kernel of the second variation of area is given by $\phi = k \cdot N$, where $k$ belongs to the kernel of 
the second variation of area.%
\footnote{Note that $\phi$ is a scalar, because with area it is enough to consider normal variations given 
by $u + \phi N$.  On the other hand $k$ is a vector since the second variation of Dirichlet's integral is defined 
as a Frechet derivative in the space of (vectors defining) surfaces $u$.  The theorem that $\phi = k \cdot N$ 
is proved in \cite{part1}.
}
The ``tangent vector'' $k = u_t$ can be shown to be a ``forced Jacobi direction'' associated
with the branch point; these are purely tangential, so as a function of the parameter $t$, the eigenfunction
goes to zero.  We calculate how the eigenfunction behaves on the $w$-plane, that is, on a region near the origin 
that shrinks to zero as $t^\gamma$.  The leading term in $t$ of this function is a harmonic function in the 
$w$-plane which inherits from the eigenfunction the property that it must have only one sign in the upper half plane.
Using this property, we are able to derive additional information about the behavior of $N$ near the origin.

Let $g(z)$ be the stereographic projection of the unit normal $N$ to $u^t$.  Then the poles of $g$, say $a_i(t)$,
are the places where $N(z)$ is the ``north pole'' $(0,0,1)$.
We analyze the way in which the $a_i(t)$   approach 
zero as $t$ goes to zero. These go to zero asymptotically in this way:
$a_i(t) = \alpha_i t^\gamma$, for some $\gamma$.  Possibly some of the $a_i$ go to zero faster than others;  
we group the $a_i$ into ``rings'',  each ring containing the roots that go to zero with the same power of $t$.
In the vicinity of  (some of) these $a_i$
the Gauss map of $u^t$ for $t > 0$ must contribute enough Gaussian area to make about $2m\pi$ ``extra'' Gaussian area,
which will disappear in the limit to account for the branch point term in the Gauss-Bonnet-Sasaki-Nitsche formula.
Consider the behavior of the normal $N$ on the boundary near the point where the branch point is when $t=0$. 
The normal must be perpendicular to the boundary, so it is confined to the vicinity of the $YZ$-plane if we take 
$\Gamma$ to be tangent to the $X$-axis.  But its behavior there might be quite uncontrollable, and moreover, several
of the $a_i$ can be associated with the same $\alpha_i$.    

In this argument there is another difficulty: the terms calculated for the eigenfunction in the 
$w$-plane might all 
cancel out, leaving the true eigenfunction hidden in the error terms.  This is the case we call ``$\HH$ constant'', since 
the eigenfunction is a power of $t$ times the imaginary part of a complex-valued function we call $\HH$.  
One of the main difficulties is to rule out this case.  We can pick any ``ring'' of roots
that go to zero as $t^{\gamma_i}$,  and look at the $w$-plane where $z=t^{\gamma_i} w$.  
We call the ``$N$-condition'' the condition that the unit normal should converge to the south 
pole $(0,0,-1)$ uniformly away from the values of $w$ where the normal takes on the north pole.
That is true on the outer ring.   One key to our argument is to show that if $\HH$ is 
constant on the outer ring, then the $N$-condition holds on the next ring,  and moreover, 
there {\em is} a next ring if $\HH$ is constant--not all the roots can go to zero as $t^{\gamma_i}$
if $\HH$ is constant on the $t^{\gamma_i}$ ring.  Once we have shown that, then on the 
innermost ring, $\HH$ cannot be constant, since there is no next ring. 
To prove this, let the $\beta_i$ be the places in the $w$-plane
where the normal 
takes on the south pole, and the $\alpha_i$ the places where it takes on the north pole.
We show that the $N$-condition on all rings follows if on each ring, there are more
$\beta_i$ than $\alpha_i$.  The explicit formula for $\HH$ that applies when $\HH$ is constant
shows first of all that there {\em is} another ring of roots going to zero faster, and second,
it
allows us to propagate the $N$-condition inwards to the next ring.

Examining the formula for $\HH$ near the origin of the $w$-plane shows that, if $\HH$ 
is not constant on the $t^\gamma$ ring, then there are no roots going to $0$ faster than $t^\gamma$, so this is the innermost ring.   Now we know that the $N$-condition holds on every ring of roots,
so the formula for the eigenfunction in the $w$-plane holds in each of the different $w$-planes
(one for each different $\gamma_i$). 
At that point the Gauss-Bonnet theorem can be used:  There are altogether $m$ roots $a_i$ where
the normal takes on the north pole.  We show (using the formula for the eigenfunction in the 
$w$-plane), applied to {\em each ring} of roots, 
 that each $a_i$ contributes at most one hemisphere of extra Gaussian area (i.e.,
the Gaussian area for $t$ positive is $2m\pi$ more than when $t=0$).  But then, there must 
really be $m$ of the $a_i$, and not some smaller number of $a_i$.  From that it follows that 
the branch points $s_i$ do not actually occur.

At that point, we know that on the innermost ring, $\HH$ is not constant;  the rest of 
the proof works only with that innermost ring.  In order to finish the proof, we 
have to analyze the behavior of $\HH$ at infinity.  Since its imaginary part $\Im\, H$
 is the limit of the first eigenfunction,
it has only one sign in the upper half-plane.  That means that many terms in the 
eigenfunction formula must cancel out.  These cancellations are very difficult to analyze 
directly.   In the last part of the proof, 
 we exhibit a differential equation that is satisfied by the function stereographic
projection of the unit normal, involving the eigenfunction $\HH$.  This equation is 
the key to the analysis of $\HH$ near infinity (in the $w$-plane).   

There is an appendix giving a ``Dictionary of Notation'', in which we list the symbols that are used in 
more than one section of the paper, and give a reminder of their definitions, as an aid to the reader.

We note in passing that as of yet, nobody has produced an example of 
a real-analytic Jordan curve bounding a minimal surface with a boundary branch point,  although Gulliver 
has exhibited in \cite{gulliver} a surface with a $C^\infty$ boundary and a boundary branch point.  The
experts I have asked all believe that they can exist, but an example is still missing.   Also, it is not 
known whether Gulliver's example is area-minimizing or not.%
\footnote{Wienholtz has used the solution of Bj\"orling's problem (\cite{hildebrandt1}, page 124) 
with boundary values $\Gamma(t^3)$ where
$\Gamma$ is a regular parametrization of a real-analytic arc, and a suitably specified normal, to produce
examples of a minimal surface partially bounded by a real-analytic arc which is not a straight line 
segment, with a branch point on the arc.}

I am grateful to Fritz Tomi for his careful reading of several early versions of this paper. 

\section{Outline of the Proof}
In this section, we list the main steps of the proof for reference. The work of B\"ohme and Tomi
reduces the task to ruling out the existence of a one-parameter family of relative minima 
of area  
bounded by a real-analytic Jordan curve $\Gamma$, running into a surface with a   
branch point.  We therefore begin with the assumption that there is a one-parameter family $u(t)$,  sometimes written ${}^t u$,
of minimal surfaces bounded by 
such that $u(0)$ has a  branch point and for $t > 0$, $u(t)$ has least eigenvalue 2 and
is immersed (both of which follow if $u(t)$ is a relative minimum for $t>0$).  The proof 
ends when this assumption is shown to lead to a contradiction.  The case of an interior 
branch point was taken care of in \cite{part1}, so  here we assume that $u$ has 
a boundary branch point when $t=0$.  The parameter domain is the upper half plane and 
the branch point is at the origin.  We assume $\Gamma$ is tangent to the $X$-axis 
at the origin and the unit normal $N$ (which extends to the branch point) is $(0,0,-1)$ 
when $t$ and $z$ are both 0.
\smallskip

(1) The Weierstrass representation has the form 
$$ u(t,z) = 
\Re \vector 
{  \onehalf \int A^2S + B^2S \,dz}
{  \iover2  \int A^2S - B^2S \,dz}
{ \int ABS \,dz}
$$
where the zeroes $a_i$, $b_i$, and $s_i$ of $A$, $B$, and $S$ go to zero as
$t^\gamma_i$, specifically as $\alpha_i t^\gamma_i$,  $\beta_i t^\gamma_i$, and 
$\zeta_i t^\gamma_i$.   The $\alpha_i$, $\beta_i$, and $\zeta_i$ are nonzero.
The order of the branch point is $2m$, its index is $k$.  The $s_i$ are 
branch points in the lower-half plane converging to 0 as $t$ goes to zero.  (Halfway
through the proof, we show they don't really exist.) 
\smallskip

(2) Fix one of the $\gamma_i$ and call it $\gamma$.  Define $w$ by $z=t^\gamma w$. 
Then in the $w$-plane,  $z-t^\gamma_j$ becomes $t^\gamma (w-\alpha_j)$ if $\gamma_j = \gamma$,
or it becomes $t^\gamma w$ if $\gamma_j > \gamma$  (a ``fast root''), or it becomes
$t^\gamma_j \alpha_j$, if $\gamma_j < \gamma$ (a ``slow root'').  Slow roots give powers of $t$,
fast roots give powers of $w$.
\smallskip

(3) In the $w$-plane, $A$, $B$, and $S$ become $\AA$, $\BB$, and $\SS$ times a power 
$t^n$  (plus terms with higher powers of $t$).   $\AA$,  $\BB$, and $\SS$ are polynomials
in $w$.   The degree of $\AA$ is the number of $\alpha_i$ that go to zero at the rate
$t^\gamma$.  There is a different $w$-plane for each distinct choice of $\gamma$ among the 
$\gamma_i$.  We refer to the $\alpha_i$, $\beta_i$, and $\gamma_i$ that go to zero 
as $t^\gamma$ as a ``ring of roots.''
\smallskip

(4) $B/A$, which is $g$ in the Weierstrass representation, is the stereographic projection
of the unit normal $N$.  When we pass to the $w$-plane, there is a number $\K$ such that
  $$ \frac B A = t^\K \frac \BB \AA  + O(t^{\K + 1}) $$
$\K$ depends on the ring, i.e. there may be a different $\K$ for each ring.
If (and only if) $\K > 0$ then we have the ``$N$-condition'', which means that as $t$ goes
to zero,  $N$ converges to the south pole $(0,0,-1)$, uniformly on compact subsets of 
the $w$-plane away from the $\alpha_i$.  On the slowest (outermost) ring, we have 
$\K = k\gamma$, where $k$ is the index of the branch point, so the $N$-condition holds there.
\smallskip

(5) There is a number $\mm$ (for each ring) such that 
$$ \int A^2 S\, dz = t^{\mm} \int \AA^2 \SS \, dw + O(t^{\mm+1}).$$
On the outermost ring, $\mm = (2m+1)\gamma$, so $\K/\mm = k/(2m+1)$.
\smallskip

(6) The eigenfunction for the least eigenvalue is given by $\phi = u_t \cdot N$.
The eigenfunction has only one sign.   Our plan is to calculate the eigenfunction and 
derive more and more information, and eventually a contradiction, from the fact that the
eigenfunction has only one sign.
\smallskip

(7) If the $N$-condition is satisfied, then $\AA$, $\BB$, and $\SS$ are real on the real axis 
(Lemma~\ref{lemma:realAA}).
\smallskip

(8) If the $N$-condition is satisfied, then the eigenfunction is given (Lemma~\ref{lemma:eigenfunctionCompact}) by  
\begin{eqnarray*} 
\frac {-1 } { t^{\mm + \K-1}}    \phi 
&=&  \Im\, \HH + O(t) 
\end{eqnarray*}
where 
\begin{eqnarray*}
\HH &=& \frac \BB \AA \, \bigg(\mm\int_0^\w \AA^2\SS\,d\w \bigg) 
 - (\mm + \K)\int_0^\w \AA\BB\SS \,d\w     
\end{eqnarray*}
Here $\AA$ and $\BB$ are polynomials in $w$, so $\HH$ is a rational function.
The principal difficulty in the proof is that $\HH$ might be constant.
\smallskip

(9) If $\HH$ is constant and the $N$-condition holds, 
then by setting the derivative $\HH_w$ equal to zero we 
find (Lemma~\ref{lemma:Hconstant}) 
$$ \frac \BB \AA = C \bigg(\int\,\AA^2\SS\,d\w\bigg)^{ \K/\mm}$$
Then $\BB/\AA$ is a polynomial, vanishing at the origin, so there are some fast roots,
hence this is not the innermost ring.  
\smallskip

(10) If $\HH$ is not constant and the $N$-condition holds, 
then there are no fast roots, since if there were, then
$\Im\, \HH$ would take two signs in the upper half plane near the origin.
Hence such a ring must be the innermost ring.
\smallskip

(11)   The $N$-condition propagates from a ring where $\HH$ is constant and the 
$N$-condition holds,  to the next ring down (Lemma~\ref{lemma:descendingrings}).  The reason for this is that, since $\BB/\AA$ is a 
polynomial,  there are at least as many $\beta_i$ as $\alpha_i$, on a ring where $\HH$ is 
constant.  Passing to a lower ring, we pick up powers of $t$ from each $\beta_i$ and $\alpha_i$,
but at least as many enter the numerator of $B/A$ as enter the denominator, so the 
power of $t$ in $B/A$ is still positive on the next ring down.
 Hence the $N$-condition holds on all rings, 
and $\HH$ is constant on all but the last ring, and $\HH$ is not constant on the last ring.
\smallskip

(12) The value of $\K/\mm$ is the same on all the rings, namely $k/(2m+1)$
 (Lemma \ref{lemma:Hnotconstant}.)
The reason for this is that, by the equation for $\HH$,  $\K/\mm$ is the ratio of the 
degrees of $\BB/\AA$ and $\int \AA^2 \SS\, dw$, and looking at the behavior near the 
origin on one plane gives us the degrees of those polynomials on the next ring down.
\smallskip

(13) Each $a_i$ contributes one hemisphere of Gaussian area (not a whole sphere).  
Since the normal is confined within $O(t)$ of the $YZ$-plane, it suffices (Lemma \ref{lemma:GaussAreaPrincipal}) to show that 
near $\alpha_i$ in the $w$-plane,  when $t$ is positive, the normal cannot take on both the 
``west pole''  $(-1,0,0)$ and the ``east pole'' $(1,0,0)$   (Lemma~\ref{lemma:hemisphericcovering}.)
The reason for this is that the eigenfunction is given by $u_t \cdot N$, and the first
component of $u_t$ near $\alpha_i$
is approximately (a power of $t$ times) $\mm \int_0^{\alpha_i}\, \AA^2\SS\,dw$,  which is positive or negative
according as $\alpha_i$ is positive or negative; hence if both the east and west poles
were taken on, $\phi$ would take two signs.
\smallskip

(14) The total extra Gaussian area is therefore $2\pi$ times the total number of $a_i$.
But by the Gauss-Bonnet theorem, it must be $2\pi m$.  Therefore the total number of $a_i$ 
is $m$. But the degree of $\AA^2\SS$ is $2m$; hence $S=1$ and the branch points $s_i$ do not exist.
\smallskip

(15) The rest of the proof focusses only on the innermost ring, where $\HH$ is not constant.
 Recall $\sigma = \int_0^w \AA^2 \, dw$.   We have (Lemma~\ref{lemma:exactB/A})
\begin{eqnarray*}
\frac {\BB}{\AA} &=& c\sigma^{\K/\mm} \exp \bigg(\int \frac {\HH_w}{\mm \sigma \BB /\AA} \, d\w\bigg)
\end{eqnarray*}
The integral is a complex path integral.  This generalizes the equation given above
for the case of constant $\HH$, since in that case $\HH_w = 0$ and the exponential factor is 1.
\smallskip

(16) $\K/\mm$ (which is equal to $k/(2m+1)$)  is an integer (Lemma~\ref{lemma:K/minteger}).
The reason why this result is important is that it makes $\sigma^{\K/\mm}$ a polynomial, so 
$\BB/\AA$ is a polynomial times the exponential factor.

The idea of the proof is that, as $w$ passes through the origin on the real $w$-axis, $\BB/\AA$ 
remains real, and does
not change signs (since there are no fast roots).  Hence the right side
must also remain real.  
Calculation shows $\HH_w \neq 0$ at the origin, so the integrand has a simple pole.
 We calculate the residue and find it to be $-\K/\mm$.  This must
be an integer in order that the right hand side remain real on both sides of the origin. 
\smallskip

(17) The non-real zeroes of $\sigma$ are also zeroes of $\BB$, and with 
the same multiplicity as zeroes of $\BB$ and zeroes of $\sigma^{\K/\mm}$ (Lemma~\ref{lemma:nonrealzeroes}).
\smallskip

(18)  Let $R$ be the degree of $\BB$ and $Q$ the degree of $\AA$.  We consider the 
Laurent expansion
$$ \sigma^{-k/\mm} \bigg( \frac \BB \AA \bigg) = 1 + \mu w^{-J} + O(w^{-J-1}).$$
On the left side we have a rational function, since $\K/\mm$ is an integer. 
If we know the degrees of the numerator and denominator of a 
rational function, we can bound the exponent of the first term in its Laurent expansion
(Lemma~\ref{lemma:Laurent}.)  On the face of it the numerator is $\BB$ and the denominator 
is $\AA \sigma^{\K/\mm}$.  But by (17), there is cancellation: all the non-real zeroes of $\BB$ 
  cancel out of both numerator and denominator.  After that cancellation the general bound 
in Lemma~\ref{lemma:Laurent} gives us $J \le R-2Q$  
 (Lemma~\ref{lemma:Laurent2}). 
 \smallskip
 
 (19) Now we come to the final steps, given in the proof of Theorem~\ref{theorem:main1}.
First we examine the integrand in 
$$
 \int \frac {\HH_w}{\mm \sigma \BB /\AA} \, d\w
$$
Using the fact that $H_w$ is asymptotic to a constant or to $w^{-2}$, we find that the 
integrand is asymptotic to $w^{-2+Q+R\pm 1}$.  Recall that $Q$ and $R$ are the degrees of $\AA$ and 
$\BB$, so $Q+R \ge 1$, which means the exponent is negative with magnitude at least 2.  Hence
the integral has a finite limit, say $L$, as $w$ goes to infinity.
\smallskip

(20) So we have  
\begin{eqnarray*}
\frac {\BB}{\AA \sigma^{\K/\mm}} 
 &=& c \, \exp \bigg(\int  \frac {\HH_w \,dw}{\mm \sigma \BB/\AA}  \bigg).  
\end{eqnarray*}
On the left side, the first term $w^{-J}$ in the Laurent expansion has $J \le R-2Q$.
But on the right, because $H_w$ is asymptotically constant or $w^{-2}$, the integrand
is asymptotically $w^{-(2+Q+R\pm 1)}$. If we differentiate, we get on the left 
a term in $w^{-(J+1)}$ and on the right a term in $w^{-(2+Q+R \pm 1)}$.  There is a tension 
here:  on the left, the first term in the Laurent series is not too far out, because 
$\K/\mm$ is an integer.  On the right, it has to be at least so far out, because $\HH_w$
is constant or $w^{-2}$.  Calculation shows  that these results cannot both 
be true unless $Q=0$, i.e. there are no $\alpha_i$; but that is also
impossible (Lemma~\ref{lemma:someA}).  That is the final contradiction.

\section{Preliminaries}

\subsection{ Analyticity at the boundary}
One of the reasons we need to work with a real-analytic boundary is that we need to know that the 
minimal surface can be analytically extended across the boundary, so that it is defined in some neighborhood 
of each boundary point and given by a power series there.  This well-known result is due to Lewy, and a proof can 
 be found in \cite{tromba-hildebrandt}, p.~107, Theorem 3.  It will be taken for granted in the rest of 
this paper and not cited explicitly.  

\subsection{ Boundary parametrization}
By a real-analytic arc, or real-analytic Jordan curve, we mean an arc or curve that can be parametrized as a real-analytic function 
of arc length.  We will consider surfaces bounded by a real-analytic Jordan curve $\Gamma$.
We suppose that $\Gamma$ passes through the origin tangent to the $X$-axis, and that $u$ takes the 
portion of the real axis near origin onto $\Gamma$, with $u(0) = 0$.  
We still are free to orient the $Y$ and $Z$ axes.  We do this in such 
a way that the normal at the branch point (which is well-defined)  points in the positive
$Z$-direction.  With $\Gamma$ oriented in this way,  there will be two positive integers
$p$ and $q$ such that  $\Gamma$ has a parametrization in the form
$$ \Gamma(\tau) = \vector
{ \tau + O(\tau^2) }
{ \frac{C_1 \tau^{q+1}}{q+1} + O(\tau^{q+2})}
{ \frac{C_2 \tau^{p+1}}{p+1} + O(\tau^{p+2})}
$$
for some nonzero real constants $C_1$ and $C_2$.   But it is possible to choose $\tau$ more 
carefully so that we do not have the $O(\tau^2)$ term in the first coordinate. 
Let $\tau(z) = X(z) = \Re \onehalf \int f - fg^2 dz $.  Then on the boundary we have 
\begin{eqnarray*}
u(z) &=& \Gamma(\tau(z)) 
\end{eqnarray*}
This parametrization and function $\tau(z)$ were inspired by Lewy's equation (see \cite{tromba-hildebrandt}, p.~108).

We therefore have
\begin{eqnarray}
 \Gamma^{\prime}(\tau) &=& \vector
{1}
{C_1 \tau^q + O(\tau^{q+1})}
{C_2 \tau^p + O(\tau^{p+1})} \label{eq:taudef}
\end{eqnarray}

\subsection{ Order and index of a branch point}
We write $u(z) = (X(z),Y(z),Z(z))$.
We make use of the Enneper-Weierstrass representation of $u$ (see e.g. \cite{hildebrandt1}, p. 112)
$$ u(z) = 
Re \vector 
{ \onehalf \int f - fg^2 \,dz}
{ \iover2  \int f + fg^2 \,dz}
{ \int fg \,dz}
$$
where $f$ is analytic and $g$ is meromorphic in the upper half-disk.

\begin{definition} The {\em order} of the branch point is the order of the zero of $f$.  
The {\em index} of the branch point is the order of the zero of $g$.
\end{definition}

\noindent{\em Remark}.  We follow \cite{nitsche}, p.~315 in this definition of ``index''.
Tromba \cite{tromba-hildebrandt} defines ``index''  to mean the sum of the orders of $f$ and $g$,
i.e. Tromba's index is Nitsche's index plus the order.

The order of a boundary branch point of a solution of Plateau's problem must be even 
(since the boundary is taken on monotonically).
It is customary to write it as $2m$, and to use the letter $k$ for the index.  Thus 
$f(z) = z^{2m} + O(z^{2m+1})$ and $g(z) = cz^k + O(z^{k+1})$ for some constant $c$.

\subsection{Gauss-Bonnet theorem for branched minimal surfaces} \label{section:Gauss-Bonnet}

The ``geodesic curvature'' $\kappa_g$ of a boundary curve $\Gamma$ bounding a surface $u$ 
is the component of the curvature vector of $\Gamma$ in the plane of $u$.  It is therefore
bounded by the magnitude of the curvature vector of   $\Gamma$, which is defined independently
of any surface. 

The Gauss-Bonnet formula says that for regular surfaces  (minimal or not)
$$  \int_\Gamma \kappa_g = 2\pi -  \int KW \, dx\, dy $$
Note that for minimal surfaces, $KW$ is negative, so both terms on the right are positive.
For minimal surfaces with branch points, there is another term in the Gauss-Bonnet formula:
$$ \int_\Gamma \kappa_g = 2\pi - \int KW \, dx\, dy  + 2\pi M$$
where $M$ is the sum of the orders of the interior branch points and half   the 
orders of the boundary branch points. See \cite{tromba-hildebrandt}, p. 195,
or \cite{nitsche}, p. 331, equation (155).  

The ``total curvature'' of a minimal surface is $-\int KW\, dx\,dy$.  This is the 
``Gaussian area'', or the area of the image on the Riemann sphere of the unit normal $N$ to the
surface.  The expression $\vert KW \vert$ is the Jacobian of $N$.  
One way to remember the Gauss-Bonnet formula is that each interior branch of order $2m$
point contributes as much $m$ spheres of Gaussian area, and 
each boundary branch point of order $2m$
contributes as much as $m$ hemispheres of Gaussian area.  This way of looking at the 
formula will be helpful when we consider a family of minimal surfaces without branch 
points for positive values of the parameter $t$, but having a branch point when $t$ is zero.
So for positive values of $t$, there must be compensating ``bubbles'' of Gaussian area.

Somewhat confusingly, the phrase ``total curvature'' is applied both to surfaces and boundary
curves.  Applied to a boundary curve, it means the integral of the magnitude of the curvature 
vector over the whole curve.  The total curvature of $\Gamma$ is thus greater than the 
left side of the Gauss-Bonnet theorem above, stated using the geodesic curvature. 

\subsection{On the zero set of real-analytic functions}
\label{section:analyticity}

A subset $U$ of  $R^n$ is a called an {\em analytic set}  if locally it is the set of simultaneous zeroes of a finite number of 
real-analytic functions from  $R^n$ to some $R^m$.  There is a classical theorem about the structure of analytic sets:

\begin{lemma} \label{lemma:analyticity}  Analytic sets are analytically triangulable.  That means that each 
analytic set is the union of a finite number of real-analytic homeomorphic images of closed simplexes,
meeting only at their boundaries.
\end{lemma}

Six references for this theorem are given in \cite{part1}, where the theorem is discussed on page 117. 
The theorem refers to {\em closed} simplexes, which means that the homeomorphisms
in the conclusion are analytic even at the endpoints of intervals or boundary points of 2-simplexes.

Several times in this paper we have reason to consider the zero set of a function that is analytic 
in two variables (one real and one complex).  When we refer to a function of a real variable $t$
and a complex variable $z$ as ``analytic'', we mean it is given by a power series in $z$ and $t$,
so in particular it is real-analytic in $t$ for fixed $z$ and complex-analytic in $z$ for fixed $t$.
The following lemma describes the zero sets of such functions. 

\begin{corollary} \label{corollary:analyticity}

 Let $A$ be an open subset of $\R \times \C$ containing
the point $(0,p)$. Let $f$ be a complex-valued   analytic function defined on $A$.
Suppose $f(0,p) = 0$ and that for each sufficiently small $t > 0$, $f(t,\cdot)$ is not constant. 
Then there exists a neighborhood $V$ of $(0,p)$ and finitely many
functions $c_i$ such that $f(t,c_i(t)) = 0$,  every zero of $f$ in $V$ for $t > 0$ 
has the form $(t,c_i(t))$,
and the $c_i(t)$ are analytic in some rational power of $t$, with $c_i(0) = p$.  
\end{corollary}

\noindent{\em Proof}. 
If $f$ does not depend on $t$ at all, the theorem is trivial. Otherwise, if 
$f$ is divisible by some power of $t$, we can divide that power out without changing the zero set of $f$ for $t > 0$, 
so we can assume that $f(t, \cdot)$ is not constant for any sufficiently small $t$, including $t=0$.
Let $S$ be the zero set of $f$. 
Then the dimension of $S$ is one, and by Lemma \ref{lemma:analyticity}, 
$S$ is composed of finitely many analytic 1-simplexes, i.e.
analytic paths given by $t = d_i(\tau)$, $w =c_i(\tau)$, 
where $\tau$ is the parameter in which $c_i$ is analytic. 
If any of these simplexes do not pass through the $(0,p)$, decrease $V$ to exclude them.  The 
ones that do pass through the $(0,p)$ should be divided into two simplexes, each ending at $(0,p)$.
Then we can assume that $\tau = 0$ corresponds to $(0,p)$, i.e. $d_i(0) = 0$ and $c_i(0) = p$.
We have
$$f(d_i(\tau),c_i(\tau))  = 0.$$
Since $d_i$ is analytic and not constant, it has only finitely many critical points.
Decrease the size of $V$ if necessary to exclude all nonzero critical points of $d_i$. Then 
either $d_i$ is increasing in some interval $[0,b]$ or decreasing in some interval $[0,b]$.
In either case it has an inverse function $\phi = d_i^{-1}$, so $d_i(\phi(t)) = t$.
If the leading term of $d_i$ is $\tau^n$, then $\phi$ is real-analytic in $t^{1/n}$. We 
can parametrize the zero set of $f$ by $\tilde c_i(t) = c_i(\phi(t))$, which is real-analytic
in $t^{1/n}$.  That completes the proof.

{\em Remark:}  In applications it will usually be possible to replace the original parameter $t$ by 
$t^{1/n}$,  enabling us to assume that the zero set is analytic in $t$.  

\subsection{Regularity results for relative minima of area} \label{section:regularity}
We need the following theorem:

\begin{theorem} \label{theorem:regularity} Let $u$ be a minimal surface bounded by a real-analytic Jordan curve $\Gamma$ and 
furnishing a relative minimum of area in the $C^0$ metric on the closed disk.  Then $u$ has neither 
interior nor boundary branch points.
\end{theorem}

It is an open problem to replace $C^0$ by $C^n$ in this theorem.  We now summarize the relevant results
from the literature.

In \cite{tromba-hildebrandt}, pp. 554-560, there is an extensive {\em Scholia}  reviewing the various 
results and proofs of regularity for interior and boundary branch points as matters stood in 2010.
This {\em Scholia} is reprinted and updated to 2012 in \cite{tromba-book}, pp.~169--175.
We are here concerned only with the case of real-analytic boundary,  but we need the regularity
of {\em relative} minimizers (of area or Dirichlet energy), not just absolute minimizers, and 
both interior and boundary regularity.

In the theory of branch points, one distinguishes between true branch points and false 
branch points (see p.~58 of \cite{tromba-hildebrandt} for the definition).  False branch points
cannot exist for any solution of Plateau's problem, minimizing or not; see the discussion 
in the {\em Scholia} ({\em op.~cit.}) for references.   In the rest of this discussion, we
consider true branch points.  

First we take up interior regularity.  It seems that \cite{beeson-regularity}
was the first to claim interior regularity for relative minima, as opposed to absolute minima, 
but Osserman's proof also 
yields such a result.  Osserman's proof works for $C^0$ relative minima, but not for $C^1$ relative
minima, because the proof involves ``smoothing off'' a wedge.  
The argument in \cite{beeson-regularity} is (also) a local 
argument, but does not involve cut-and-paste: it shows how to decrease the area (in a $C^n$ smooth way) in a neighborhood of a branch point.  To obtain a global result 
one must supplement that argument by a smoothing argument to remove the ``crease'' introduced by decreasing 
area locally.  Wienholtz has correctly pointed out that smoothing argument given in \cite{beeson-regularity}
 works only for $C^1$, although 
$C^n$ is claimed.  Nevertheless the result is correct for $C^1$, as confirmed in the {\em Scholia} just mentioned.
So as far as interior branch points go,  we could replace $C^0$ by $C^1$ in the theorem.

In \cite{tromba-book},  Tromba made direct calculations attempting to decrease the Dirichlet integral (globally) for a 
minimal surface with an interior branch point.  While his method was successful in most cases,  in certain cases he was still forced to make 
a local argument, so that $C^1$ can still not be replaced with $C^n$ in the theorem; but Tromba did provide 
an independent proof of the interior regularity for $C^1$ relative minima.  

Now we consider boundary regularity.  In the case of real-analytic boundary curves, area and energy
can be locally decreased near a true boundary branch point, as proved independently in 
  \cite{gulliver-lesley} and \cite{white}.  Both proofs rely 
 on an Osserman-style cut-and-paste.
  Hence they work for $C^0$ relative minima, although the papers claim the theorem only for 
absolute minima, and the {\em Scholia} cited above also does not state that the theorem works
for $C^0$ relative minima. 

However, we do not have to rely on my claim that these two proofs work for $C^0$ relative minima. 
Instead we can rely on the work of
 Tromba \cite{tromba-book}, who also made calculations for a minimal surface with a boundary 
branch point.   While serious difficulties remain in the case of boundaries that are only $C^n$,
Tromba's method does work
in the case of a real-analytic boundary curve.   But, as mentioned
above, Tromba also needed a local argument in some cases, so his proof also works only for 
$C^0$ relative minima.  For Tromba's explicit statement that his method works in the 
real-analytic boundary case, see p.~168 of \cite{tromba-book}.  It is the last sentence in the book, not counting
the {\em Scholia}.

\subsection{ Dirichlet integral and area}
Let $D$ be a plane domain; for our purposes, we may suppose $D$ is a disk or a half-plane, so we will not 
worry about the conditions on the boundary  $\partial D$ of $D$.  
The Dirichlet integral, or ``energy'', of a surface $u$ defined in $D$ is defined by 
$$ E(u) = \onehalf \int_D u_x^2 + u_y^2 \, dx\,dy $$
The area functional is defined by 
$$ A(u) = \int_D \vert u_x \times u_y \vert \, dx\, dy$$
These functionals are defined on various spaces of functions on $D$, for example $C^{n,\alpha}$ or the 
Sobolev spaces $W^{k,p}$.  Since we will be restricting attention to surfaces bounded by a fixed real-analytic Jordan
curve $\Gamma$, we are interested in a space of functions defining surfaces bounded by $\Gamma$.  
These function spaces can be considered as Hilbert manifolds, and the tangent space at a surface $u$ 
is a vector space of functions defined on $D$ whose values are tangent to $\Gamma$ on the boundary of $D$.
These functionals have Frechet derivatives, which we denote
by $DE(u)$ and $DA(u)$.   These are linear mappings on the tangent space.  

Minimal surfaces are defined as critical points of $E$,  but it is well known that they 
are also critical points of $A$.  In connection with $A$, it is common to consider only ``normal variations'', 
i.e. to restrict $DA(u)$ to the subspace of the tangent space  consisting of those $k$  such that $k(x,y)$
is normal to $u$.  It can be shown that if $DA[u]$ vanishes on this subspace then $u$ is minimal.
(Details can be found in \cite{nitsche}, p. 94).   In connection with $E$, it is common to consider 
only harmonic surfaces; or equivalently, spaces of functions defined on the boundary of the parameter domain
and mapping it to $\Gamma$ in a way homotopic to the identity (since such functions have a unique harmonic extension 
to the interior).

\subsection{The second variation of Dirichlet's integral and the forced Jacobi fields }

At a minimal surface $u$,  we can consider the 
second variation $D^2 E(u)$.  This is a bilinear mapping on the tangent space at $u$.  The members of this 
tangent space are ``tangent vectors'' $k$, i.e. functions defined on $\partial D$ such that $k(\xi)$ is tangent
to $u(\xi)$ for $\xi$ on $\partial D$.
The kernel of $D^2 E(u)$ consists of tangent vectors $k$ such that $D^2 E(u)[k,j] = 0$ for all 
tangent vectors $k$.  We write $D^2 E(u)[k]$ for $D^2 E(u)[k,k]$.  By 
diagonalizing the bilinear form, one can prove that the kernel of $D^2 E(u)$ consists of those $k$ for which 
$D^2 E(u)[k] = 0$.  
 
\begin{theorem}[Tromba] The second variation of Dirichlet's integral is given by 
$$ D^2E[u](h,k) = \int k(h_r - \tilde h_\theta)\, d\theta$$
where $k = \lambda u_\theta$ and $h = \eta u_\theta$ and  $\tilde h = \eta u_r$ and $\tilde k = \lambda u_r$.
The tangent vector $k$ to the minimal surface $u$ belongs
to Ker $D^2 E[u]$ if and only if 
$$ u_\theta(k_r - \tilde k_\theta) = 0$$
or equivalently
$$ k_z u_z = 0$$
\end{theorem} 

\noindent{\em Proof}.  The first formula is equivalent to the one 
given in \cite{tromba} (bottom of p.~53, with $h$ and $k$ 
interchanged); for a stand-alone one-page proof of it by direct calculation, and the easy derivation 
of the second two formulas from the first, see \cite{beeson-notes}.
 
Consider the kernel equation $u_\theta(k_r - \tilde k_\theta) = 0$.  One way in which this could 
be satisfied is if $k_r -\tilde k_\theta = 0$; vectors $k$ satisfying this condition and not 
induced by the conformal group are called ``forced Jacobi fields'' or ``forced Jacobi directions''.
Tromba proved that they do not occur in the absence of branch points, and that in the presence of branch
points there are two for each interior branch point (counting multiplicities) and one for each 
boundary branch point, so that the space of forced Jacobi fields is finite dimensional.  (There can 
be at most finitely many branch points, even if the boundary is not real-analytic, as long as the 
total curvature of the surface is finite, thanks to the Gauss-Bonnet formula for branched minimal 
surfaces.)  The forced Jacobi directions are just the directions $k$ such that the function 
$K = k + i \tilde k$ is complex analytic, i.e. such that $\tilde k$ is the conjugate harmonic 
function of $k$. 

Another important characterization of the forced Jacobi fields is this: they are exactly the 
tangent vectors of the form 
$$ k = \Re\, (i\omega z \, u_z)$$
where $i \omega z$ is a function meromorphic in the parameter domain, and having a pole of 
order at most $m$ at each branch point of order $m$.  Any function $\omega$ with suitable behavior
on the boundary, and poles of the right orders at the branch poitns, will produce a tangent vector
by this equation.  The reason for writing the equation with $\omega z$ instead of with $\omega$
is that in case the parameter domain is the unit disk, the appropriate boundary condition is that 
$\omega$ be real on $S^1$.  In case the parameter domain is the upper half plane, the condition 
is that $i \omega z$ be real on the $x$-axis.  The Appendix of \cite{bohme-tromba} contains Tromba's
treatment of the forced Jacobi fields.

\begin{lemma}[Tromba \cite{tromba}] \label{lemma:Tromba} Suppose $u$ is a minimal surface, and $k$ is a tangent vector 
belong to $Ker D^2E[u]$ whose harmonic extension is everywhere tangent to $u$.  Then $k$ is 
a forced Jacobi direction or a direction induced by the conformal group.
\end{lemma}

\subsection{The second variation of area}
\begin{lemma}  [Laplacian of the Gauss map]\label{lemma:laplacianN} Let $u$ be a minimal surface.
Then the Laplacian of its unit normal is given by  
$\Delta N = 2KWN$.
\end{lemma}

\noindent{\em Proof}.
 To prove this elegantly, we make use of the general fact that the Laplace-Beltrami operator of any 
surface $S$, applied to the position vector of $S$, is exactly twice the mean curvature of $S$.
(For a proof of that fact, see p. 45 of \cite{hildebrandt1}, formula (30), with $f$ the position vector of $X$
so that $H_f$ in formula (27) is the first fundamental form of $X$.)
Apply this fact to the Riemann sphere, whose position vector $h(w)$ coincides with its 
unit normal. Thus $\Delta h = -2h$.  Next, note that the map from the disk $D$ to the Riemann
sphere induced by the Gauss map of $u$ is a conformal map with Jacobian $-KW$.  Under a conformal 
map, the Laplace-Beltrami operator changes to the Laplace-Beltrami operator on the range 
surface, multiplied by the Jacobian of the mapping.  Hence $\triangle N = 2KW N$, and the 
lemma is proved. 

One can define the second variation of area, $D^2 A(u)$, on the original tangent space, or on the 
subspace of the tangent space corresponding to normal variations.  Nitsche tells us (\cite{nitsche}, p. 95) that 
we get the same result, but that the calculation is too long to put in his 562 page book.  It is customary 
to consider the second variation of area on normal variations only.  The second
variation $D^2A(u)\phi$ is a bilinear functional on normal variations; technically 
normal variations are ``tangent vectors'' to a Hilbert manifold of surfaces, but of course
they are not tangent to the surfaces themselves.  
  Abusing notation slightly, we 
write $D^2 A(u)[\phi,\psi]$ to stand for 
$D^2 A(u)[k,j]$ where $k$ is the tangent vector $\phi \cdot N$ and $j=\psi\cdot N$ (and $N$ is the unit normal to $u$).
Thus $D^2 A(u)[\phi,\psi]$
is, in classical terms, 
$$ \frac {\partial^2}{\partial s \partial t} A(u + s\phi + t\psi) $$
evaluated at $t=0$ and $s=0$.  This derivative is only used when $\phi$ and $\psi$
are in the kernel of the first variation $DA$ (which of course holds when $u$ is minimal).
In case $\phi = \psi$ we recover the classical second variation of area:
$$ D^2 A(u)[\phi,\phi] = \frac {\partial^2}{\partial t^2} A(u + t\phi),$$
which by a further abuse of notation is usually written $$D^2A(u)[\phi].$$
 One can calculate that
\begin{equation} 
 D^2A(u)[\phi,\psi] = \int_D \psi (-\Delta \phi +2KW\phi) \, dx\, dy, 
 \end{equation}
 where $K$ is the Gauss curvature and $W = \vert u_x \times u_y \vert$.  See for example
 \cite{nitsche}, p. 96, where the case $\psi = \phi$ is calculated, but the calculation easily adapts to the 
 bilinear second variation.

\begin{lemma} \label{lemma:KerD2A}
$\phi$ belongs to the kernel of $D^2A(u)$ (for normal variations) just in case $\phi$ is zero on the 
boundary of the parameter domain and satisfies
$$ \Delta \phi = 2KW\phi$$
in the interior.
\end{lemma}

\noindent{\em Proof.}  Apply the fundamental lemma of the calculus of variations to the preceding formula.

\subsection{Connections between $D^2A$ and $D^2E$}

\begin{theorem} \label{theorem:EtoA}
Let $u$ be a minimal surface in $R^3$ with $C^n$ boundary, and unit normal $N$.
Let $k$ be in $Ker D^2E[u]$. Then $\phi = k \cdot N$ belongs to $Ker D^2 A[u]$.
\end{theorem}
\begin{corollary} If $Ker D^2 A[u]$ has no kernel among normal variations, $Ker D^2 E[u]$ 
contains only the conformal and forced Jacobi directions.
\end{corollary}

\noindent{\em Proof}. The Corollary follows immediately from the theorem and Tromba's
lemma (Lemma \ref{lemma:Tromba}).   We now prove the theorem.  Suppose $k$ is in $Ker D^2 E[u]$.  By 
Lemma \ref{lemma:KerD2A}, it will suffice to  show that
$\phi =  k\cdot N$ satisfies $\triangle \phi - 2KW \phi = 0$. We have 
$$\triangle \phi =  \triangle(k\cdot N) = (\triangle k) \cdot N + 2\nabla k \nabla N + k \triangle N.$$
The first term vanishes because $k$ is harmonic.  We claim the second term vanishes also.
To prove this, fix a point $z_0$ in the unit disk, and choose coordinates $a$ and $b$
in a neighborhood of $z_0$ that diagonalize the first fundamental form at $z_0$, so that 
$N_a = \kappa_1 u_a$ and $N_b = \kappa_2 u_b$, where $\kappa_1$ and $\kappa_2$ are the 
principal curvatures of $u$ at $z_0$.  If these equations hold in a whole 
neighborhood, then $a$ and $b$ are called ``local curvature coordinates''; it costs some 
trouble to prove they exist, and we do not need them; we need the first fundamental form 
to be diagonalized at one point $z_0$ only.  To do this, 
we take $a$ and $b$ to be a 
certain linear combination of $x$ and $y$.  
Let  $\nu$ be the angle between the positive $x$-direction and the positive $a$-direction
(so $\nu$ is a function of $z_0$ but not of $z$).  Then we define $w = e^{i\nu}(z-z_0)$,
and define $a$ and $b$ by $w = a+ib$, so $a$ and $b$ are coordinates rotated by $\nu$ 
from $(x,y)$. 

Because $u$ is a minimal surface, we have 
$\kappa_1 = - \kappa_2$.  Then at $z_0$ we have
\begin{eqnarray*}
\nabla k \cdot \nabla N &=& k_a N_a + k_b N_b \\
&=& k_a \kappa_1 u_a + k_b \kappa_2 u_b \\
&=& \kappa_1 (k_au_a -k_bu_b)
\end{eqnarray*}
We have $u_w = u_a-iu_b$ and $k_w = k_a -ik_b$, so 
\begin{eqnarray*}
\nabla k \cdot \nabla N &=& \kappa_1 \Re\,(k_w \cdot u_w)  \\
&=& \kappa_1 \Re\,((k_z z_w) \cdot (u_z z_w)) \\
&=& \kappa_1 \Re\,(z_w (k_z  \cdot u_z)) \\
&=& \kappa_1 \Re\,(e^{-i\nu} (k_z \cdot u_z))
\end{eqnarray*}
since $z_w = e^{-i\nu}$.
  Since $k$ is assumed to be in $Ker D^2 E[u]$, we have
$k_z u_z = 0$.  Hence the term $\nabla k \nabla N$ vanishes at $z_0$.
But $z_0$ was arbitrary; hence $\nabla k \nabla N$ vanishes everywhere,
 and we have proved $\triangle \phi = k \cdot N$.

The proof of the theorem is thus reduced to proving $\triangle N = 2KWN$.  But this is 
Lemma \ref{lemma:laplacianN}.  That completes the proof.

The converse of this theorem is also true (but more difficult):
  every $\phi$ in the kernel of $D^2 A(u)$ arises as 
$k \cdot N$ for some $k$ in the kernel of $D^2 E(u)$.  We do not need this result in this paper, but the curious
can find a proof in \cite{beeson-notes} or \cite{beeson-Bonn}.

\subsection{Stereographic projection}
We take the ``Riemann sphere'' to be the unit sphere $\{(x,y,z): x^2 +y^2 + z^2 = 1\}$.
Stereographic projection ${\bf St}$ is defined by 
$$ {\bf St}((x_1,x_2,x_3)) = \frac{ x_1+ix_2}{1-x_3}$$
and its inverse is given by 
$$ {\bf St}^{-1}(z) = \frac 1 {1 + \vert z\vert^2} \vector { 2 \,\Re\ z} {2 \,\Im\ z}{\vert z\vert^2 -1}.$$
Thus the equator projects onto the unit circle, and the ``north pole'' $(0,0,1)$ projects onto $\infty$,
while the southern hemisphere projects onto the unit disk, with the south pole $(0,0,-1)$ going to origin.  
The picture then has the plane passing through the equator of the sphere.   
We follow \cite{osserman}, p. 46, in these details, and we mention them because some other authors use 
a sphere of radius $1/2$, with the picture having the sphere entirely above the plane, tangent to 
the plane where the south pole of the sphere touches the origin of the plane.

\subsection{The eigenvalue problem associated with the second variation of area}

Associated with a conformal map $N$ defined on a region $\Omega$ in the plane and taking values in the Riemann sphere,
there is a natural eigenvalue problem:
\begin{eqnarray*}
\Delta \phi -  \onehalf \lambda \vert \nabla N \vert^2 \phi &=& 0 \mbox{\qquad in $\Omega$ } \\
\phi &=& 0 \mbox{\qquad on $\partial \Omega$}
\end{eqnarray*}
The Jacobian of $N$ is $\onehalf \vert \nabla N \vert^2$.
In case $N$ is the Gauss map of a minimal surface,  the Jacobian is also  $-KW$, so 
$ \onehalf \vert \nabla N \vert^2 = - KW$ and 
the eigenvalue equation becomes
$$ \Delta \phi + \lambda KW \phi = 0$$
As proved above in Lemma \ref{lemma:KerD2A}, if $\phi$ is in the kernel of $D^2A(u)$, 
then $\phi$ is an eigenfunction of this equation for $\lambda = 2$.  But for purposes of 
this section, $N$ can be any map from the disk to the Riemann sphere.%
\footnote{Some readers may be familiar with another 
form of the eigenvalue equation for which the critical eigenvalue is zero rather than 2,  or with this form 
but with a factor of 2 inserted so that the critical eigenvalue is 1 instead of 2; both 
forms are discussed in \cite{nitsche}, p. 103, cf. equations (62) and (62${}^\prime$).  }

We denote the least eigenvalue $\lambda$ of this problem by $\lambda_{N,\Omega}$, or simply 
by $\lambda_\Omega$ when $N$ is clear from the context.  Sometimes we use the notation $\lambda_{\min}$.
Sometimes, for a minimal surface $u$, we speak of the ``least eigenvalue of $u$''  rather than the 
``least eigenvalue of the second variation of $u$'' or ``the least eigenvalue of the eigenvalue problem 
associated with the second variation of $u$.''
 
The least eigenvalue is well-known to be equal to the infimum of the 
Rayleigh quotient 
$$ R[\phi] = \frac { \int \int_\Omega \vert \nabla \phi \vert ^2 \,dx\, dy}
                   { \int \int_\Omega \onehalf \vert \nabla N \vert^2 \phi^2 \,dx\,dy}.
$$ 
When we speak of the least eigenvalue $\lambda_\Omega$ of a region $\Omega$ on the Riemann sphere,
we mean the following:   Let $\Delta$ be the stereographic projection of $\Omega$ and $N$ the 
inverse of stereographic projection.  Then $\lambda_\Omega := \lambda_{N,\Delta}$.  The eigenvalue
problem $\Delta \phi - \frac 1 2 \lambda \vert \nabla N \vert^2 \phi$ on $\Delta$ is equivalent
to the problem $\Delta \phi = \lambda \phi$ on $\Omega$, where now $\Delta$ is the Laplace-Beltrami
operator on the sphere.  If $\Omega$ contains the north pole, we should use stereographic projection 
from some point not contained in $\Omega$.  We do not need to discuss the case when $\Omega$ is the 
entire sphere.

{\em Example}. 
We compute the least eigenvalue  when $N(\Omega)$ is a hemisphere.
In this case the eigenfunction in the lower hemisphere is minus the 
$Z$-component of $N$.  For example with $\Omega$ equal 
to the unit disk and $g(z) = z$, we have $N(z)$ the inverse of stereographic 
projection.  With $\vert z \vert = r$ and $z = x + iy$ we have
$$N(z)= \frac 1 {1+r^2} \vector{2 x} {2y }{r^2 -1 }.$$
The eigenfunction $\phi$ is given, with $\vert z \vert = r$, by 
$$ \phi(z) = \frac{ 1-r^2}{1+r^2}.$$
A few lines of elementary computations (or a couple of commands to a computer algebra program) show that
\begin{eqnarray*}
\triangle \phi &=&  \phi_{rr} + \frac 1 {r} \phi_r  \\
   &=&  \frac {8(r^2-1)} {(1+r^2)^3} \\
\end{eqnarray*}
and
\begin{eqnarray*}
  \vert \nabla N \vert^2 &=& N_x^2 + N_y^2 = \frac{ 8r}{(1+r^2)^2} \\
  \vert \nabla N \vert^2 \phi &=&\frac {8(1-r^2)} {(1+r^2)^3} 
\end{eqnarray*}
Hence, $\Delta \phi - \vert \nabla N \vert^2 \phi = 0$, which means the eigenvalue of a hemisphere is 2.

\begin{lemma} \label{lemma:lambda-monotonicity} [Passing to the Riemann sphere does not decrease the eigenvalue] Let $\Omega$ be a connected open set on the sphere
with least eigenvalue $\lambda_\Omega$.  Suppose $\Delta$ is a region (open set) in the 
plane and $\Omega \subseteq N(\Delta)$ and $N(\partial \Delta) \cap \Omega = \phi$.  Suppose that
boundary of $\Omega$ is $C^2$ and the boundary of $\Delta$ is piecewise $C^2$ with the pieces meeting
at positive angles. 
 Then $\lambda_{N,\Delta} \le \lambda_\Omega$,
and strict inequality holds if $N^{-1}\partial\Omega$ contains an interior point of $\Delta$.
\end{lemma}

\noindent{\em Remark}. Regarding the assumptions on the boundaries,  the proof requires that 
the least eigenvalue be the minimum of the Rayleigh quotient, and that the gradient of the 
least eigenfunction of $\Omega$ not vanish at any boundary point.  The hypotheses given 
imply these conditions but still allow $\Delta$ to be a half-disk. See \cite{gilbarg}.  

\noindent {\em Proof}.
Let $\phi$ be the least eigenfunction of $\Omega$ and define $\psi$ on $\Delta$ by 
setting $\psi(z) = \phi(N(z))$ if $N(z) \in \Omega$ else $\psi(z) = 0$.  Then 
$\psi$ is admissible in the Rayleigh quotient for $\Delta$, since $\phi$ is zero outside 
$\Omega$ and $N(\partial \Delta) \cap \Omega = \phi$.
 The Rayleigh quotient in question is 
$$ \frac { \int \int_\Delta \vert \nabla \psi \vert ^2 dx\, dy}
         { \int \int_\Delta \frac 1 2 \vert \nabla N \vert ^2 \psi^2}
$$
(The factor $1/2$ was explained above.) 
Since $\Omega \subseteq N(\Delta)$, on the support of $\psi$, $N$ is a covering map (i.e., locally a 
homeomorphism), except 
at the points of ramification of $N$, which are isolated.  Since $\Omega$ is connected 
and $N(\partial \Delta) \cap \Omega = \phi$,
the number of sheets over (cardinality of the pre-image of) $N$ of each non-ramification point is the same.   
Hence each of the two
integrals in the (numerator and denominator of the) Rayleigh quotient is the number of sheets times the corresponding integral
on the Riemann sphere, with $\phi$ in place of $\psi$.  That is, the Rayleigh quotient for 
$\psi$ on $\Delta$ equals the Rayleigh quotient for $\phi$ on $\Omega$, which is $\lambda_\Omega$.
Since $\lambda_\Delta$ is the minimum of such Rayleigh quotients, $\lambda_\Delta \le \lambda_\Omega$.  
Now suppose there is an interior point $p$ of $\Delta$ in $N^{-1} (\partial\Omega)$.
By Hopf's lemma, $\nabla \phi$ is never zero at a 
point on $\partial \Omega$ where the boundary is $C^2$.
 By analyticity, $\nabla N$ is zero only at isolated points; so 
there is a point $q$ near $p$ which is still on $N^{-1} (\partial\Omega)$ at which $\nabla N$ 
is not zero and $\nabla \psi$ is not zero.  Hence near $q$, the set $N^{-1}(\partial\Omega)$ is a smooth arc,
and  $\nabla \psi$ is zero on one side of it and bounded away from zero on the other side.
Hence we can ``smooth out the edge'' near $q$ to obtain a function $\psi^\prime$ which is 
admissible for the Rayleigh quotient and has smaller Rayleigh quotient than $\psi$.  Hence
$\lambda_\Delta < \lambda_\Omega$.  This completes the proof of the lemma.

\begin{corollary} 
If the Gaussian image of a minimal surface defined in (the open set) 
$\Delta$  contains a hemisphere, and
at least one boundary point of the hemisphere, 
then the eigenvalue $\lambda_{N,\Delta}$ is less than 2.  
\end{corollary}

\noindent{\em Proof}.  The eigenvalue of a hemisphere is 2.

\subsection{ Dependence of the Gauss map on a parameter}

In this section we present a result which will be applied to the Gauss map of a 
one-parameter family of minimal surfaces.  However, we present it here as a
result about conformal mappings from a plane domain to the Riemann sphere.

Let  $D^+$ be the upper half of the unit disk.  We will consider a one-parameter
family of analytic mappings from $D^+$ to the Riemann sphere.  (These arise
in our work as the Gauss maps of a one-parameter family of minimal surfaces.)  
 We write $t$ for the parameter and 
$N(t,z)$ for the value of the mapping.  For brevity we often omit the explicit
$t$-dependence and write $N(z)$, or $N^t(z)$.   Thus $N^t$ is the map $N(t,\cdot)$
from $D^+$ to the Riemann sphere, for a fixed $t$.  We define $g = {\bf St} \circ N$.  The map
$N$ is thus conformal except at the zeros of $g^{\prime}$ and poles of $g$ of order
more than one.  (These are sometimes known as {\em ramification points} of $g$.)

We will suppose that for sufficiently small $t$, the map $g$ is a quotient of 
functions $B/A$, where $B$ and $A$ are analytic jointly in $t$ and $z$ (in the sense 
explained in  section~\ref{section:analyticity}), and that when 
$t=0$, $A$ has a zero of order $m$, and $B$ has a zero of order $m+k$, for some positive integers
$m$ and $k$.  (In our work, these arise from a boundary branch point of order $2m$ and 
index $k$, but here it is not necessary to be so specific.)

By Corollary \ref{corollary:analyticity}, 
near the origin there
exist paths $a_i(t)$ and $b_i(t)$  describing the poles and zeroes of $g$, respectively,
so that each $a_i$ and $b_i$ is analytic in a rational power of $t$.   Let $s$ be the
greatest common divisor of the denominators of these rational powers, and replace
$t$ by $t^{1/s}$.   Then $a_i$ and $b_i$ will be analytic in $t$.   Thus we may suppose that
the parameter $t$ has been chosen such that $a_i$ and $b_i$ are analytic in $t$.  We
also may suppose without loss of generality that for all sufficiently small positive $t$,
we have $a_i(t) \neq b_j(t)$
for all $i$, $j$, since otherwise by analyticity we would have, 
for some $i$ and $j$,  $a_i(t) = b_j(t)$
for all sufficiently small $t$, and then the factor $z-a_i(t)$ could be cancelled out of both
$A$ and $B$.  The partial derivatives of a quotient of real analytic functions are 
again quotients of real-analytic functions, so the structure of the set of critical points
of the real and imaginary parts of such a function is also known:  it is the union of a 
finite set of real-analytic arcs.

%
\subsection{Dependence of eigenvalues on a parameter}
 Let $D$ be the unit disk and let $0 < \alpha < 1$.  Let $F:[0,1] \times D \rightarrow \R$, and suppose that 
$F(t,x,y) \ge 0$ for all $t$, $x$, and $y$.  For each fixed $t$ let $F^t$ be the function of $x,y$ given by $F^t(x,y) = F(t,x,y)$.
Then the eigenvalue problem $\phi = \lambda F^t \phi$ (with $\phi = 0$ on the boundary)  has eigenvalues $\lambda_1^t,\lambda_2^t, \cdots \lambda_n^t$, given 
in non-decreasing order.  For $n$ such that the eigenspace of $\lambda_n$ has dimension greater than 1,  that leaves some
ambiguity about the indexing of the eigenvalues.  Suppose that $F$ depends in some smooth way on $t$, perhaps even 
real-analytically.  Then what can be said about the dependence of $\lambda_n$ on $t$?  At least this much:

\begin{lemma}[Continuous dependence of eigenvalues] \label{lemma:continuouseigenvalues}
Consider the eigenvalue problem $\Delta \phi = \lambda F_t \phi$, with $\phi = 0$ on the boundary.
For each fixed
value $t_0$ of $t$, there are continuous functions $\lambda_n(t)$ defined in a neighborhood of $t_0$ giving the eigenvalues.
\end{lemma}

\noindent{\em Proof}. The proof can be found on page 419 of \cite{courant-hilbert}.  

 In fact more is true.  If $F$ depends real-analytically on $t$,  then it is possible to number the eigenvalues in such a way that 
 they are given by real-analytic functions of $t$, as shown in \cite{kato}, p. 370, p. 387; in general this numbering will not be in order 
 of size, since two eigenvalues, for example $\lambda_2$
and $\lambda_3$,  might ``cross'' at a certain value of $t$.  Of course no other eigenvalue can cross or even touch $\lambda_1$ since
its eigenspace is always one-dimensional.  For our purposes simple continuity of $\lambda_1$ and $\lambda_2$ suffices.

\subsection{Tomi's no-immersed-loops theorem}

We state and sketch the proof of a theorem due to Tomi \cite{tomi1}, or maybe it is due to Tomi and B\"ohme \cite{bohme-tomi}.  It is difficult to give an exact reference for this 
theorem as the paper where it is stated \cite{tomi1}  contains deeper theorems about the structure of the solution set of Plateau's problem,
and the calculation needed for this proof is referenced to \cite{bohme-tomi} where, at the crucial point, the paper says
``Durch eine elementare aber etwas m\"uhsame Rechnung finder man$\ldots$''  (by an elementary but somewhat tiresome computation one finds),
and does not give the computation.   Also, they required the boundary to 
be $C^{4,\alpha}$,  but that was for other reasons in their paper; $C^2$ is enough for the theorem stated here (although in our paper,
the boundary is always real-analytic anyway).

\begin{theorem} [Tomi] Let $\Gamma$ be a $C^{2}$ Jordan curve and suppose $u = u(t)$ is a periodic one-parameter family of 
minimal surfaces, $C^2$ as a function of $z$ and $C^1$ in $t$, bounded by $\Gamma$ for each $t$, and satisfying a three-point condition.  Suppose that $u_t$ is not identically zero as a function of $z$ for any $t$, and that each $u(t)$ has $\lambda_{\min} = 2$.
Then some $u(t)$ has a branch point, either in the 
interior or on the boundary.
\end{theorem}

\noindent{\em Remarks}.  By a ``periodic family'', we mean that $u(t + 2\pi)= u(t)$; the exact period is not relevant.  The condition 
$\lambda_{\min} = 2$ will be fulfilled if $u(t)$ is a relative minimum of area, but it is a more general condition.  The condition that 
the $u(t)$ satisfy a three point condition is only needed to guarantee that $u_t$ is not a conformal direction. 
\smallskip

\noindent{\em Proof}. Suppose, for proof by contradiction, that $u(t)$ has no branch point.
 Since each $u(t)$ is a minimal surface, the first variation of Dirichlet's integral $E$
 is zero, so $E(t) = E[u(t)]$ is 
constant.  Hence the second derivative $\partial^2 E/\partial t^2 = 0$.  Hence the second variation of $E$ is zero
in the direction $u_t$.  That is, 
$$ D^2 E[u](u_t) = 0$$
Define 
$$ \phi := u_t \cdot N$$
where $N$ is the unit normal to $u(t)$.
(We suppress the $t$-dependence in our notation, writing $u$ instead of $u(t)$ and not indicating the $t$-dependence of $\phi$ and $u_t$.)
Because there is no branch point, there are no forced Jacobi directions.  Because of the three-point condition, $u_t$ is not a 
conformal direction. 
Then $D^2[A](\phi) = 0$ as shown in the previous section.   Since $u$ has no branch points, $\phi$ is not identically zero.
By Theorem \ref{theorem:EtoA}, $\phi$ is an eigenfunction of 
$$ \Delta \phi = 2 KW\phi$$
over the parameter domain $D$, with $\phi = 0$ on the boundary $\partial D$.
We define the ``volume integral'' to be 
$$ V(t) := \int_D u \cdot ( u_x \times u_y) \, dx dy \ = \ \int_D u \cdot N\, W\, dx dy$$
where $N$ is the unit normal to $u$.  
 The key to the proof is the ``m\"uhsame Rechnung'' that 
$$ \frac {\partial  V} {\partial T} = V_t =  \int \phi W \,  dx dy $$
That computation, not found in either \cite{tomi1} or \cite{bohme-tomi}, is written out in detail in \cite{beeson-notes}.  Since by hypothesis $\lambda_{\min} = 2$,  for each $t$ the function $\phi = \phi(t)$ has only 
one sign in the interior of the parameter domain.  Since $\phi$ is continuous in $t$ (because $u$ is $C^1$ in $t$),
 that sign is the same for all $t$ in 
$[0,2\pi]$.  Now choose $t_0$ at which $V(t)$ has its minimum value.   Then we have
\begin{eqnarray*}
0 &=& V_t(t_0) \\
&=& \int u_t \cdot NW\,  dx dy \\
&=& \int \phi W\, dx dy
\end{eqnarray*}
But since $\phi$ is not identically zero and has one sign, this is a contradiction.

\subsection{B\"ohme and Tomi's structure theorem}

Since the 1930s it has been known that, if $\Gamma$ bounds
infinitely many minimal surfaces, there is a sequence $u_n$
of them converging to a minimal surface $u$ bounded by $\Gamma$.
We need an ``analytic compactness'' result that will give us 
a one-parameter family of minimal surfaces instead of just a sequence.
B\"ohme and Tomi \cite{bohme-tomi, bohme} and Tromba \cite{tromba} applied
nonlinear global analysis to the theory of minimal surfaces.  
 We will state and use 
the structure theorem of B\"ohme and Tomi.  This theorem is Satz 2.4, page 15 in \cite{bohme-tomi}.
That theorem is stated more generally, to cover the case of constant mean curvature as well as 
the minimal-surface case;  here we state it only for minimal surfaces.

 \def\F{\mathcal F}
 
\begin{theorem} [B\"ohme-Tomi structure theorem] Let $u$ be a minimal surface bounded by a real-analytic
Jordan curve in $R^3$.  Let $0 < \alpha < 1$. 
 Let $\F$ be the space of $C^{2+\alpha}$ maps from $D$ to $\R^3$ satisfying a three-point condition 
 and whose restrictions to the boundary map $S^1$ onto the range of $\Gamma$ and 
 are homotopic to $\Gamma$.   Let ${{\mathcal M^\prime}(\Gamma)}$ be the set of minimal surfaces in $\F$.
 Then there is a neighborhood $W$ of $u$ in the space $C^{2,\alpha}$, such that 
 $W \cap {\mathcal M^\prime}(\Gamma)$
 is an analytic subset of a finite-dimensional analytic submanifold of $\F$.
 \end{theorem}

\noindent{\em Warning}.  The surfaces in $W$ are {\em not} required to map the boundary of the parameter domain
onto $\Gamma$  monotonically.  That is why we do not simply say that surfaces in $\F$ are ``bounded by $\Gamma$.''
The notation ${\mathcal M^\prime}$ instead of ${\mathcal M}$ means that monotonicity is not required.
\smallskip

\noindent{\em Remarks}.  The proof goes by introducing a map ${\mathcal L}$ from $W$ to $C^{2+\alpha}(D) \times C^{2+\alpha}{S^1}$,
defined by
$${\mathcal L}(u)  = (\Delta u, u_r \cdot u_\theta)$$
Since a harmonic function conformal on the boundary is also conformal in the interior, the zeroes of ${\mathcal L}$ are 
exactly the minimal surfaces.  B\"ohme and Tomi show that ${\mathcal L}$ is a Fredholm map, which means that 
the Frechet derivative $D{\mathcal L}(u)$ has finite-dimensional kernel and cokernel.  The finite-dimensional 
manifold in the theorem is parametrized by the kernel of ${\mathcal L}$.
 The theorem essentially boils down to an application of the implicit function theorem; 
if $D{\mathcal L}$ is zero in the $n$ directions $t = (t_1,\ldots,t_n)$, and nonzero in other directions, 
then the non-isolated zeroes of ${\mathcal L}$ are locally analytic functions of $t$.  That remark, of course,
is not a proof--it is intended only to convey the idea of the proof in \cite{bohme-tomi}.

The kernel of $D{\mathcal L}$ is essentially the kernel of $D^2 E$, the second variation 
of Dirichlet's integral.  We now know that kernel consists of the forced Jacobi directions (in case $u$ has branch points),
plus directions with non-zero normal components.  

In \cite{bohme}, Satz 6 and Satz 11,  B\"ohme spells out the consequences of the structure theorem more explicitly.
With $\Gamma$, $u$, and $\F$ be as in the previous theorem, let ${\mathcal M}^\prime (\Gamma)$  be the set of minimal surfaces in $\F$. 
The notation ${\mathcal M}^\prime(\Gamma)$ is B\"ohme's.  It contrasts with ${\mathcal M}(\Gamma)$, in that 
members of ${\mathcal M}(\Gamma)$ must take on $\Gamma$ monotonically (that is, topologically);  both sets impose
a three-point condition.  
Let ${\mathcal L}$ be as defined above, so ${\mathcal L}(u) $ is zero just when $u$ is a minimal surface in $\F$.  Then
the kernel of $D{\mathcal L}(u)$ is a finite-dimensional vector space $V$.    B\"ohme's Satz 11 uses terminology from his 
Satz 6; putting in that terminology, Satz 11 becomes:

\begin{theorem}[B\"ohme analytic structure theorem]   
There is a neighborhood W of $u$ and an analytic map $j$ defined on the unit ball $U$ in $R^n$, where $n$ is the 
dimension of the kernel of $D^2 E$, such that
 $ {\mathcal A} := j^{-1} {\mathcal M}^\prime(\Gamma)\cap W)$ is an analytic subset of $U$.  
\end{theorem}

B\"ohme does not define ``analytic subset'',  and we wish to make sure there is no ambiguity about it.  It means,
the set of simultaneous zeroes of a finite number of real-analytic functions from $U$ to some $R^m$.  
 Explicitly, we have 
$$ {\mathcal M}(\Gamma)\cap W) = \{ v \in j(U) : {\mathcal L}(v) = 0$$
Then $j(t)$ is a function of $z$, and $j(t)(x,y)$, which we can write as $j(t,x,y)$, is real-analytic in
all three variables, as B\"ohme points out.   We have 
$$ {\mathcal A} = \{ t \in V \ : \ ({\mathcal L} \circ j)(t)) = 0 \} $$
This ``analytic set'' is the set of zeroes of a real-analytic function from a finite-dimensional space to 
an infinite-dimensional space, but since $U$ is finite-dimensional and ${\mathcal L} \circ j$ is analytic,
${\mathcal L} \circ j(U)$ is also a finite-dimensional 
manifold.  

 The structure of the set of minimal surfaces ${\mathcal M}^\prime(\Gamma)$ bounded by $\Gamma$ in a neighborhood of $u$ is thus 
reduced to the study of the local structure of analytic subsets of $R^n$.  If we want to also require monotonicity
on the boundary, the set ${\mathcal M}(\Gamma)$ becomes {\em semianalytic}; that is, defined by a Boolean combination
of analytic sets.  For our purposes, it is enough to consider ${\mathcal M}^\prime(\Gamma)$.

\begin{corollary} Let  $\Gamma$ be a real-analytic Jordan curve
in $\R^3$. Then the set of minimal surfaces bounded by $\Gamma$ is locally analytically triangulable.
 That is,  if $u$ is a minimal surface bounded by $\Gamma$,
and $\F$ is the space of minimal surfaces bounded by $\Gamma$ and satisfying a three-point condition,
then there is   a neighborhood $W$ of $u$ in $\F$ such that the set of minimal 
surfaces in $W$ is a union of analytic images of simplexes in $R^n$.
\end{corollary}

\noindent{\em Proof}. By B\"ohme's structure theorem and Lemma \ref{lemma:analyticity}.
\smallskip

B\"ohme's theorem leaves open the possibility that the analytic simplexes have their vertices at $u$.
If $u$ is a relative minimum of area,  more can be said:  the space of minimal surfaces is one-dimensional
and a 1-simplex passes right {\em through} $u$, not just up {\em to} it:

\begin{lemma}[Tomi] \label{lemma:TomiStructureLemma}  Let $u$ be a relative minimum of area bounded by a real-analytic Jordan curve $\Gamma$.  
If $u$ is not isolated, then the set of minimal surfaces bounded by $\Gamma$ in some neighborhood of $u$
is exactly a one-parameter family $u^t$, with $u^0 = u$, defined for $t$ in some open interval containing 0.
\end{lemma}

\noindent{\em Proof}.  See \cite{tomi2}.

\section{ Basic setup and Weierstrass representation}

\subsection{A one-parameter family terminating in a branched surface}

The starting point of our work is the following theorem.  The statement and proof can be sketched as follows:
if $\Gamma$ bounds infinitely many relative minima of area, then there is a one-parameter family of minimal surfaces
$u^t$ bounded by $\Gamma$.  If they are all without branch points, then they all have $\lambda_{\min} = 2$, and
by compactness the family must loop, but Tomi's theorem prevents that, so they must run into a branched minimal surface.  
To make this ``run into'' precise we need B\"ohme's structure theorem; and to guarantee that they remain relative 
minima until they run into the branched surface, we need the continuity of the first two eigenvalues.   Here are 
the details:

\begin{theorem} \label{theorem:one-parameter}
Let $\Gamma$ be a real-analytic Jordan curve bounding infinitely many minimal surfaces without branch 
points (satisfying a three-point condition) and with $\lambda_{\min} \ge 2$.
 Then there exists a one-parameter family of minimal 
surfaces $u^t$, bounded by $\Gamma$,  defined for $t$ in some closed interval containing 0, 
such that for $t > 0$, $u^t$ is a minimal surface without branch points and with $\lambda_{\min} \ge 2$, and $u^0$ has 
an (interior or boundary) branch point.   
\end{theorem}

\noindent{\em Proof}.  Suppose that $\Gamma$ bounds infinitely many minimal surfaces with $\lambda_{\min} \ge 2$ and
without branch points.   
Then by compactness, there is a surface $u$ that is a limit of a sequence $u_n$ of such surfaces, satisfying the same 
three-point condition as the $u_n$.  
Since the limit of minimal surfaces is minimal, this surface is also minimal.    Let $\F$ be the set of minimal surfaces
satisfying the same three-point condition and taking $S^1$ to $\Gamma$ with winding number 1, but not necessarily monotonically.
By (the corollary to) B\"ohme's structure theorem, the set of minimal surfaces in $\F$ near $u$
is analytically triangulable.   Pick $n$ large enough that $u_n$ lies on the triangulation, 
i.e.,  $u_n$ is the image under some analytic map $j$ defined on a $k$-simplex $S$ in $R^n$ taking
values in $\F$.    Since $u_n$ has no branch points, it takes the boundary monotonically.
By   Lemma \ref{lemma:TomiStructureLemma}, the set of minimal surfaces bounded by $\Gamma$ near $u_n$ is one-dimensional,
so the simplex $S$ is a $1$-simplex, and there is a one-parameter family of minimal surfaces $u^t$ defined for
$t$ in $[0,1]$, connecting
$u^0 = u$ to $u^1 = u_n$.%
\footnote{Actually, B\"ohme's structure theorem suffices for this proof, since according to that theorem there
is {\em some} simplex containing $u_n$, possibly of dimension more than 1,  that extends analytically to $u^0$.
In that simplex we can find some analytic path.}   
Let $q$ be the infimum of the set of $t \le 1$ such that $u^t$ has no boundary or interior branch points.  This set
is an open set, since it is the set where $ \vert \nabla u \vert^2 > 0$ on the closed unit disk;
 hence $q < 1$.  
Now, we just forget about the $u^t$ for $t < q$.  That is,
replacing $t$ by $t + q(1-t)$,  we can assume without loss of generality that $q=0$; that is, for $t > 0$ the surface
$u^t$ has no  branch points and, in view of the fact that there are no boundary branch points, each $u^t$
 takes the boundary monotonically. 
 Let $\phi = u_t \cdot N$.  Since for $t > 0$, $u^t$ has no branch points,
  Theorem \ref{theorem:EtoA} implies that
 for $t > 0$, $\phi$ is an eigenfunction of 
 $\Delta \phi + 2KW\phi$.  Hence  2 is an eigenvalue of $\Delta \phi + \lambda KW \phi$ for each $t > 0$.
 For $t = 1$ we have $\lambda_1 = 2 $ and $\lambda_2 > \lambda_1$, since the least eigenvalue
 has a one-dimensional eigenspace. By Lemma \ref{lemma:continuouseigenvalues}, $\lambda_1$ and $\lambda_2$ depend continuously on $t$ for 
 $t >0$, since for $t>0$ we have $KW >0$.  Because the eigenspace of $\lambda_1$ is one-dimensional, we cannot 
 have $\lambda_2 = \lambda_1$ for any value of $t$;  hence the two continuous functions $\lambda_1(t)$ and 
 $\lambda_2(t)$ do not cross (or even touch), and we have $\lambda_1 = 2$ for all $t$ in $(0,1]$.    That completes the proof.
 
 \begin{corollary} Let $\Gamma$ be a real-analytic Jordan curve bounding infinitely many relative minima 
 of area in the $C^0$ metric.  
relative minima of area.  There there exists a one-parameter family of minimal 
surfaces $u^t$, bounded by $\Gamma$,  defined for $t$ in some closed interval containing 0, 
such that for $t > 0$, $u^t$ is a minimal surface without branch points and with $\lambda_{\min} \ge 2$, and $u^0$ has 
an (interior or boundary) branch point.   
\end{corollary}

\noindent{\em Proof}.  This would be true with $C^n$ in place of $C^0$, if we knew that $C^n$ relative minima have no 
branch points.  See the discussion in Section \ref{section:regularity} above.
\smallskip

\noindent{\em Remark}.  The conclusion of the corollary can be strengthened to claim that each $u^t$ for $t > 0$ not 
only has $\lambda_{\min} = 2$ but is also a relative minimum of area.  We do not need this result, and it is more 
complicated to prove than the lemma above; in \cite{part1} it is proved by reference to a normal form theorem for 
the Dirichlet integral due to Tromba.

\subsection{ Weierstrass Representation}

We suppose we are in the following situation: a real-analytic Jordan curve, not lying in a plane, bounds a 
one-parameter family of minimal surfaces $u(t,z)$. The surfaces are parametrized by $z$ in the upper half plane, 
so that on the real line each $u(t,\cdot)$ is a reparametrization of $\Gamma$.  The surface $u(0,\cdot)$ has
a  boundary branch point of order $2m$ and index $k$ at the origin.  The surfaces $u(t,\cdot)$ for $t > 0$
have no branch points, and they have $\lambda_{\min} = 2$.  
All our arguments will be local; we shall only be concerned
with what happens near the origin.\footnote{
We do not need to assume that the branch point is a true branch point.  We can rule out the existence of a 
one-parameter family of the type considered without that assumption.  This does not, however, provide a new proof
of the non-existence of false branch points which are not the terminus of a one-parameter family of relative
minima of area.}

We suppose that $\Gamma$ passes through the origin tangent to the $X$-axis.  Then since $\Gamma$ is 
a Jordan curve, it is not contained in a line.   We still are free to orient the $Y$ and $Z$ axes.  We do this in such 
a way that the normal to $u(0, \cdot)$ at the branch point (which is well-defined)  points in the positive
$Z$-direction.  

Finally,  we can assume that for all $t$ we have $u(0,0) = 0$.   If this is not already the case,
then we can make it so by applying (for each $t > 0$) a conformal transformation that moves $u(0,0)$ to the origin.
If we do this while using the upper half-plane for a parameter domain, then the transformation in question is 
simply the map $z \mapsto z- \xi(t)$, where $u(\xi(t)) = 0$.  $\xi(t)$ can be found analytically in a rational power of $t$,
and we can replace $t$ by a new parameter in which the family will be analytic.   This assumption will be used (for example)
to know that $u$ can be found on the boundary by integrating $u_x$ starting from 0.

The surfaces $u(t,z)$ have an Enneper-Weierstrass representation
$$ u(t,z) = 
\Re \vector 
{  \onehalf \int f - fg^2 \,dz}
{  \iover2  \int f + fg^2 \,dz}
{ \int fg \,dz}
$$
The branch points of $u(t,\cdot)$ are the places where $f$ and $fg^2$ vanish simultaneously.  
In \cite{part1} we were concerned with interior branch points only.  In that case we can 
argue that since there are no branch points for $t > 0$,  the zeroes of $f$ are of even 
order and coincide with the poles of $g$, so that $f$ and $g$ have the forms stated in 
Lemma 3.3 of \cite{part1}, namely:
\begin{eqnarray*}
 f(z) &=& f_0(z) \prod_{i=1}^m(z-a_i)^2\\
 g(z) &=& \frac{\displaystyle h(z) \prod_{i=1}^{m+k}(z-b_i)}{\displaystyle \prod_{i=1}^m(z-a_i)}
\end{eqnarray*}

In this paper, boundary branch points are considered.   Hence, we must consider
 the possibility that $f$ and $fg^2$ may have some common zeroes $s_i(t)$, corresponding 
to a branch point (or branch points) in the lower half plane that converge to the origin as $t$
goes to zero.
 These common zeroes of $f$ and $fg^2$ 
could occur with different multiplicities in $f$ and $fg^2$.  Say that $s_i$ occurs
with multiplicity $m_i$ in $f$ and $n_i$ in $fg^2$.
We must have $\Im(s_i(t)) < 0$ for $t > 0$, since there will be branch points at the $s_i$, 
and $s_i(t)$,  like $a_i(t)$ and $b_i(t)$, depends analytically on some power of $t$
called $t^\gamma$ and converging to 0 as $t$ approaches 0.  

We define 
$$S(z) = \prod_{i=1}(z-s_i)^{\min(n_i,m_i)} $$
Then all zeroes of $S$ are zeroes of $f$ and$fg^2$ with at least the multiplicity they have in $S$.
If $s_i$ occurs with a greater multiplicity in $fg^2$ than in $f$, the extra occurrences are zeroes
of $g^2$, and hence of $g$, so $n_i - m_i$ is even.  Similarly, if $s_i$ occurs with a greater
multiplicity in $f$ than in $fg^2$, the extra occurrences are poles of $g^2$, and hence poles of $g$,
so $m_i-n_i$ is even.   We have 
$$ \frac {fg^2} f  = g^2 $$
so the remaining zeroes of $f$ and $fg^2$ are double.  We can therefore define functions $A$ and 
$B$, real-analytic in $t$ and complex-analytic in $z$, such that 
\begin{eqnarray*}
 f(z) &=&  A^2 S \\
 f(z)g^2(z) &=&  -  B^2 S\\
 g(z) &=& \frac{fg^2}{f} = \frac{i  B}{A} \\
f(z)g(z) &=& i A B S
\end{eqnarray*}
We already defined the $a_i(t)$ and $b_i(t)$ to be the poles and zeroes of $g$, respectively.
Let $2N$ be the number of zeros of $S$. (It must be even since the number of zeroes of $f$ 
is $2m$ and all the zeroes of $f$ besides those of $S$ are double.)
Then we can find functions $A_0$ and $B_0$, real-analytic in $t$ and complex-analytic in $z$,
not vanishing at the origin, such that
\begin{eqnarray*}
A(t,z) = A_0(t,z) \prod_{i=1}^{m-N}(z-a_i) \\
B(t,z) = i \lambda B_0(t,z) \prod_{i=1}^{m+k-N}(z-b_i) \\
\end{eqnarray*}
and $A_0(0,0) =  1$, $B_0(0,0)$ is real, and $\vert \lambda \vert = 1$ with $\lambda \neq -1$.
Note that the common zeroes $s_i$ do not enter into the formula for $g$.   
We can assume $A_0(0,0) = 1$,
since we can always change $z$ to a constant times $z$.   The leading term of $B(z)$, namely
$i\lambda B_0(0,0)$ is some nonzero complex constant, but because of the presence of $\lambda$
in the formula, we can assume $B_0(0,0)$ is real.

Note that with this definition of $A$ and $B$, some of the $a_i$ and/or some of the $b_i$ may 
have been defined to be equal to some of the $s_i$, in case $s_i$ occurs with a greater multiplicity
in $f$ than in $fg^2$ or vice-versa.  This is possible since these extra occurrences are double,
since they occur as roots or poles of $g$.  Each branch point $s_i$ will either occur as an $a_i$,
or as a $b_i$, or neither one (in case its multiplicities in $f$ and $fg^2$ are equal);  so it will 
still be true that $a_i(t)$ is not identically equal to $b_j(t)$ as $t$ varies, even if the $s_i$ 
occur and are treated this way.

As we have already remarked, it follows from Corollary \ref{corollary:analyticity} that the $a_i$, $b_i$, 
and $s_i$ depend analytically on $t$.
Although they all approach zero as $t$ goes to zero,
some may go to zero faster than others.  The roots $a_i$, $b_i$, and 
$s_i$ occur in ``rings'',  going to zero at different speeds.  There is a (finite) sequence of powers $\gamma_1 < \gamma_2 < \ldots $
such that every root goes to zero as $t^{\gamma_n}$ for some $n$; we say those roots 
belong to the $n$-th ring.  Specifically, the roots in the $n$-th ring are those for which 
there exist nonzero complex constants $\alpha_i$, $\beta_i$, and $\zeta_i$ such that 
\begin{eqnarray*}
a_i &=& \alpha_i t^{\gamma_n} + O(t^{\gamma_n + 1}) \\
b_i &=& \beta_i t^{\gamma_n} + O(t^{\gamma_n + 1}) \\
s_i &=& \zeta_i t^{\gamma_n} + O(t^{\gamma_n + 1}) 
\end{eqnarray*}

\section{The Gauss map and the $\w$-plane} \label{section:wplane}
 Suppose that $\gamma$ is one of the 
numbers $\gamma_n$; that is, $\gamma$  is a number 
such that some of the roots $a_i$, $b_i$, or $s_i$ go to zero as $t^\gamma$, but not necessarily the least such number.
 Roots that go to zero as some higher
power of $t$ are called ``fast roots'', and ones that go to zero as a lower power of $t$ are called ``slow roots''.
\begin{eqnarray}
&& \mbox{$Q$ is the number of $a_i$ that go to zero as $t^\gamma$ or faster.}  \label{eq:Qdef} \\
&& \mbox{$S$ is the number of $s_i$ that go to zero as $t^\gamma$ or faster.} \label{eq:Sdef} \\
&& \mbox{$R$ is the number of $b_i$ that go to zero as $t^\gamma$ or faster.}\label{eq:Rdef} \\
&& \mbox{$P$ is the total number of fast roots.}\label{eq:Pdef}
\end{eqnarray}  
Thus if $\gamma = \gamma_1$, we have $Q = m$, $S = 2N$,  and $R = m+k$.
Recall that $z = a_i(t)$ are the zeroes of $A$ and $z = b_i(t)$ are the zeroes of $B$.
We will generally be considering only one ``ring'' of roots at a time, and when we speak of 
$\alpha_i$, $\beta_i$, and $\zeta_i$, we have in mind those roots that go to zero as $t^\gamma$:
\begin{eqnarray*}
a_i &=& \alpha_i t^\gamma(1 + O(t)) \\
b_i &=& \beta_i t^\gamma(1 + O(t)) \\
a_i &=& \zeta_i t^\gamma(1 + O(t)) \\
\end{eqnarray*}
We have 
$$A = A_0(t,z) \prod_{i=1}^{m-N}(z-a_i(t))$$
We introduce a new variable $\w$ by $z = t^\gamma \w$.  Then a factor
$z-\alpha_i t^\gamma$ becomes $t^\gamma (\w - \alpha_i) $,  and a factor $z - \alpha_i t^{\gamma_i}$ for $\gamma_i < \gamma$
becomes $$t^{\gamma_i}(z t^{\gamma - \gamma_i} - \alpha_i) = t^{\gamma_i}(-\alpha_i + O(t)).$$  
Thus the slow roots introduce constants and powers of $t$, and the roots that go to zero as $t^\gamma$ 
introduce factors containing $w-\alpha_i$, and the roots that go to zero faster than $t^\gamma$ introduce powers of $w$ times $(1+O(t))$.
 All products without explicit indices will be understood to be over those roots which are $O(t^\gamma)$,
and in such products (but not elsewhere) we will assume that the $\alpha_i$ corresponding to the fast roots are zero.  Thus
$\prod (w-\alpha_i)$  contains a factor of $w$ for each fast root, and the total number of factors is the number of $a_i$
that go to zero as $t^\gamma$ or faster.  Define, for the roots that go to zero as $t^\gamma$ or faster,
\begin{eqnarray*}
\alpha &:=& \mbox{the product of all the $-\alpha_i$ over the slow $a_i$ } \\
\beta &:=&  \mbox{the product of the $-\beta_i$ over the slow $b_i$ times $B_0(0,0)$ }\\
\zeta &:=&  \mbox{the product of the $-\zeta_i$ over the slow $\zeta_i$} \\
\a_i &:=& a_i/t^\gamma \\
\b_i &:=& b_i/t^\gamma \\
\s_i &:=& s_i/t^\gamma \\
\tAA &:=& \alpha\prod (\w - \a_i) \qquad\mbox{(the product is over $i$ such that $a_i = O(t^\gamma)$)}\\
\tBB &:=& i\lambda B_0\beta\prod (\w - \b_i) \\
\tSS &:=& \zeta\prod (\w - \s_i)
\end{eqnarray*}
Thus $\a_i$ is zero for the fast roots $a_i$, but it could also be zero for some roots $a_i$ that
go to zero as $t^\gamma$, and $\a_i$ is not defined for the slow roots.  When $t=0$, $\a_i$ 
becomes $\alpha_i$.  Over $\AA$, $\BB$, or $\SS$, a tilde 
indicates $t$-dependence. When $t = 0$ we drop the tilde:
\begin{eqnarray*}
\AA &=& \alpha \,  \prod(\w-\alpha_i) \\
\BB &=& i \lambda \beta  \, \prod(\w-\beta_i) \\
\SS &=& \zeta\prod(\w-\zeta_i) 
\end{eqnarray*}

Now we consider the relation between $A(t,z)$ and $\tAA(t,w)$.  As discussed above, the roots that go to zero
as $t^\gamma$ or faster contribute factors of $t^\gamma w$ or $t^\gamma(w-\alpha_i)$ to $\tAA$, and the slow 
roots contribute factors $-\alpha_i t^{\gamma_i}$, with $\gamma_i < \gamma$.  The product of the factors
contributed by all the slow roots is thus a constant $\alpha$ (defined above) times a power $t^d$ of $t$, 
where $d$ is the sum of the $\gamma_i$ over all the slow roots $a_i$.    Since the constant $\alpha$ is 
part of the definition of $\AA$, it does not show explicitly in the relation between $A$ and $\AA$:  
\begin{eqnarray}
A = t^d t^{Q\gamma} \AA (1 + O(t)) \label{eq:ddef}
\end{eqnarray}
  Similarly, we have numbers $e$ and $f$ such that
\begin{eqnarray}
B &=&  t^e t^{R\gamma}  \BB (1 + O(t)) \label{eq:edef} \\
S &=&  t^f t^{S\gamma}  \SS (1 + O(t))\label{eq:fdef}
\end{eqnarray}
No ambiguity will result from using $f$ for this power of $t$ since in this meaning, it will always occur in an exponent,
and the Weierstrass function $f$ never occurs in an exponent.  Similarly the letter ``$S\,$'' has 
to do double duty, in exponents as the number of roots $s_i$ that go to zero as $t^\gamma$ or faster,
and not in exponents as a function of $z$, as can be seen in (\ref{eq:fdef}). (We are short of letters in the alphabet, and the function $S$ will eventually disappear.)  We define 
\begin{equation} \label{eq:Kdef}
 \K := (R-Q)\gamma + (e-d).
\end{equation}
The importance of $\K$ stems from the following equation, which results from dividing (\ref{eq:edef})
by (\ref{eq:ddef}):
\begin{equation}
\frac B A = t^{\K} \frac {\BB}{\AA}  + O(t^{\K + 1}).\label{eq:Kprop}
\end{equation}
This equation holds on any compact subset of the $w$-plane away from the zeroes of $\AA$.
When $\gamma = \gamma_1$, so we are working with the   first-ring roots, then $\K = k\gamma_1$, since no roots 
go to zero slower than $\gamma_1$, so $e=d=0$, and $R-Q = k$. 

Our next aim is transform our formulas for the Weierstrass representation from their 
expressions in $z$ and $t$ to express them in terms of $w$ and $t$. Starting with the 
expressions in terms of $z$ and $t$, we have
\begin{eqnarray*}
u  &=& \Re \vector{ \onehalf \int f - fg^2 \, dz }
              { \frac i 2 \int f + fg^2\, dz}
              { \int fg \, dz} \\
   &=&\Re \vector { \onehalf \int A^2 S + B^2 S \, dz}
               { \frac i 2  \int A^2 S - B^2 S \, dz}
               { i \int ABS \, dz} 
\end{eqnarray*}
Now we express $A$, $B$, and $S$ as appropriate powers of $t$ times $\tAA$, $\tBB$, and $\tSS$, as given above;  and
we note that since $z = t^\gamma w$, we have $dz = t^\gamma dw$, so when we change the variable of integration to $w$,
we pick up an extra power of $t^\gamma$.  Replacing $\tAA$ by $(1+O(t))\AA$, and similarly for $\tBB$ and $\tSS$, we have by (\ref{eq:ddef}), (\ref{eq:edef}), and (\ref{eq:fdef})
$$
u = (1 + O(t))\, \Re \vector {\onehalf \int t^{2d + f} t^{(2Q+2S+1)\gamma} \AA^2 \SS + t^{2e+f} t^{2R+2S+1)\gamma} \BB^2 \SS \,dw }
                   {\frac i 2 \int t^{2d + f} t^{(2Q+2S+1)\gamma} \AA^2\SS - t^{2e+f} t^{(2R+2S+1)\gamma} \BB^2 \SS \,dw }
                   { \int t^{d +e + lf} t^{(Q+R+S + 1)\gamma} \AA\BB\SS \, dw}  
$$
Since the expressions for these exponents are a bit unwieldy, we shorten them by introducing
\begin{eqnarray} \label{eq:Mdef}
 \mm &:=& 2d + f + (2Q+S+1)\gamma 
\end{eqnarray}
With this definition we have
\begin{eqnarray}
A^2S &=& t^{\mm-1} \tAA^2\tSS  \ = \ \AA^2\ (1+O(t)) \label{eq:mm} \\
B^2S &=& t^{\mm-1+ 2\K} \tBB^2 \tSS  \ = \ \BB^2\ (1+O(t)) \label{eq:mm2k}
\end{eqnarray}
Note that $dz = t^\gamma dw$, so when we integrate we get an extra factor of $t^\gamma$.
Using $\mm$ and $\K$,  the previous equation becomes
\begin{equation}
u = (1 + O(t))\, \Re \vector {\onehalf \int t^{\mm} \AA^2 \SS +  t^{\mm + 2\K} \BB^2 \SS \,dw }
                   {\frac i 2 \int t^{\mm} \AA^2\SS -  t^{\mm + 2\K} \BB^2 \SS \,dw }
                   { \int t^{\mm+\K} \AA\BB\SS \, dw} \label{eq:WeierstrassAA}
\end{equation}
To get an intuitive grasp of these formulas, consider the special case in which we are 
working with the first-ring roots. Then $2Q+S+1 = 2m+1$, and if there are also no slow 
roots then $e=d=f=0$, so $\mm = (2m+1)\gamma$.  If there are some fast roots they
contribute extra powers of $t$, namely $t^{e+d+f}$.   In general, we will find that $\mm$
plays the role that $(2m+1)\gamma$ would play in the simplest case.

We next review the formula for the unit normal.  
The function $g$ in the Weierstrass representation, which is the stereographic projection 
of the unit normal $N$, is given by $ g(z) = i B/A$, and the normal itself is given by 
\begin{equation}
N = \frac 1 {\vert g\vert^2 + 1} \vector{ 2\Re\,(g)}{ 2\,\Im(g)}{ \vert g \vert^2 -1}. 
\end{equation}
We write $\bar A$ for the complex conjugate of $A$.  
Substituting $g = i  B/A$, we have
$$ N = \frac {1}{1 + \vert B/A\vert^2} \vector{ -2\Im\,(B/A)} { 2\Re(B/A)} {-1 + \vert B/A\vert^2}$$
Expressing this in terms of $\tBB/\tAA$, we have $B/A = t^\K \tBB/\tAA$, and hence
\begin{equation} \label{eq:Nofwexact}
 N = \frac {1}{1 + t^{2\K}\vert \tBB/\tAA\vert^2} \vector{ -2t^{\K}\Im\,(\tBB/\tAA)}
                                          {2 t^{\K} \Re\,(\tBB/\tAA)} 
                                          {-1 + t^{2\K}\vert \tBB/\tAA\vert^2}
\end{equation}
We remind the reader that $\None$ is the first component of $N$, and $\Ntwo$ the second component. We have
\begin{equation} \label{eq:N1OverN2}
 \frac \None \Ntwo =  - \frac {\Im\,(\tBB/\tAA)} {\Re\,(\tBB/\tAA)} 
\end{equation}
We record in a lemma the geometric significance of $\K$:
\begin{lemma} \label{lemma:Kpositive} $\K > 0$ if and only if the unit normal $N$ converges to $(0,0,-1)$ on compact subsets of the $w$-plane  away from the zeroes of $\AA$, as $t$ goes to zero.
\end{lemma}

\noindent{\em Proof}.   The stereographic projection of the unit normal $N$ is $B/A$,
which by (\ref{eq:Kprop}) is given by $t^\K \tBB / \tAA$.   This goes to 
zero on compact subsets of the $w$ plane if and only if $\K > 0$;  but 0 is the stereographic projection of $(0,0,-1)$.
Alternately the lemma can be proved by direct examination of (\ref{eq:Nofwexact}).
 
\begin{lemma} \label{lemma:realAA}
Suppose that on compact subsets of the $\w$-plane away from the zeroes of $\AA$, $N$ converges to $(0,0,-1)$ as $t$ goes to zero.  Then
\smallskip

(i) $\AA$, $\SS$, and $\BB$ are real when $w$ is real;
\smallskip

(ii) the $\alpha_i$  and $\zeta_i$ are real and the $\beta_i$ are
 real and/or occur in complex-conjugate pairs.
\end{lemma}

\noindent{\em Remark}. The hypothesis is certainly true for the first ring of $a_i$, the slowest
roots.  We will later show inductively
that it is true for all the rings of $a_i$, i.e. when each $\gamma_n$ is used to define $\w$.   
\smallskip

\noindent{\em Proof}.  
First we prove that $\alpha_i$ does not have positive imaginary part. 
Suppose, for proof by contradiction, that it does.  Then consider a circle of
radius $\rho$ about $\alpha_i$, where $\rho$ small enough that the entire 
circle lies in the upper half plane and each other $\alpha_i$ is either 
inside or outside the circle.  Then by hypothesis, on the boundary of the circle, 
$N$ goes to $(0,0,-1)$ as $t$ goes to zero.
Hence the Gauss map covers the upper hemisphere inside the circle, so 
the least eigenvalue of $u^t$ is less than 2, for $t>0$ sufficiently small; but that is a contradiction.  That contradiction proves that 
$\alpha_i$ does not have positive imaginary part.

Let $Y(\w) := {}^2u(t^\gamma \w)$, where again, ${}^2u$ is the second component of $u$.
Then for real $w$ 
$$Y(\w) = \frac {C_2}{q+1} \tau^{q+1}(1 + O(s))$$
 where $\tau$, $q$, and $C_2$ are as in (\ref{eq:taudef}), and  $q > 0$.
\begin{eqnarray*}
\tau &=& X(w) = \int_0^z A^2S + B^2S\, dz  \\
 &=& t^\mm \int_0^\w \AA^2\SS  + t^{2\K} \BB^2 \, d\w \ ( 1 + O(t)) 
\end{eqnarray*}
by (\ref{eq:mm})  and (\ref{eq:Kdef}).  The hypothesis that $N$ converges to $(0,0,-1))$ 
away from the zeroes of $\AA$ implies $\K > 0$, by Lemma \ref{lemma:Kpositive}. 
Hence the term in $t^{2\K}$ can be absorbed into the error term, and we have for real $w$
$$ \tau = t^{\mm}  \int_0^\w \AA^2\SS \, d\w \ ( 1 + O(t)) $$
and hence
$$ Y(\w) = \frac {C_2}{q+1}  t^{(q+1)\mm}  \int_0^\w \AA^2\SS \, d\w \,(1 + O(t))$$
 On the other hand from the Weierstrass representation (\ref{eq:WeierstrassAA}), we have 
$$Y(\w) =  t^{\mm}\, \Im \int_0^\w \tAA^2\tSS\, d\w  + O(t^{2d + 2Q\gamma + f + S\gamma + 2 }) $$
Since $q > 0$, the power of $t$ shown explicitly in this equation is less than the power of $t$ in the preceding equation. 
It follows that 
$$\Im \int_0^\w \AA^2\SS\, dw = 0$$
since it is independent of $t$ but from the two previous equations must be $O(t)$.
Therefore the polynomial $\int_0^\w \AA^2\SS \,d\w$ is real on the real axis.  It follows that
its coefficients are real.  Hence the coefficients of its derivative $\AA^2\SS$ are real.  Hence the
roots of $\AA^2\SS$,  which are the $\alpha_i$ and $\zeta_i$, come in 
complex-conjugate pairs.  But we proved above that no $\alpha_i$ has positive imaginary part.
Hence they are all real when $t=0$.  Similarly, the $\zeta_i$ do not have positive imaginary part, since 
$u^t$ has no branch points for $t > 0$.  Hence  the $\zeta_i$ are real when $t=0$.

Similarly, with $ Z(\w) = {}^3u(t^{\gamma}) \w $ in place of $Y$, we find
$$\Im \int_0^\w \AA\BB\SS\, dw = 0$$
That is, $\int_0^\w \AA\BB\SS\, d\w = 0$ is real on the real axis. 
Therefore also its derivative  $\AA\BB\SS$ is real on the 
real axis.   But $\AA\BB\SS$ is a polynomial; since it is real on 
the real axis, its coefficients are real.  We have already proved that 
 $\AA\SS$ is real on the real axis; hence $\BB$ is also real on the real 
axis.  Hence its roots, which are the $\hat \b_i$ when $t=0$, are real or occur in complex-conjugate
pairs.  That completes the proof of the lemma.

\section{The eigenfunction}
In this section, we compute the eigenfunction of the 
variational problem associated with the second variation of area.  This is shown in \cite{part1} to be 
$\phi = u_t \cdot N$, where $ u_t = du/dt$ is the 
``tangent vector'' to the one-parameter family, and $N$ is the unit normal.  
Only one assumption will be made for our eigenfunction computation:  we assume that $\K > 0$.
(See (\ref{eq:ddef}), (\ref{eq:edef}), (\ref{eq:fdef}), and (\ref{eq:Kdef}) for the definition of $\K$ and the quantities
used to define $\K$.) 
  This assumption ensures that 
the unit normal, whose stereographic projection is $B/A$, is equal to $t$ to a positive power times $\BB/\AA$, so that 
on compact subsets of the $\w$-plane away from zeroes of $\AA$, the unit normal tends to $(-1,0,0)$.    Specifically we have
\begin{eqnarray*}
 \frac B A &=& t^{(R-Q)\gamma + (e-d)} \frac \tBB \tAA \\
&=& t^{\K\gamma} \frac \tBB \tAA  \qquad \mbox{by (\ref{eq:Kdef}) } 
\end{eqnarray*}
and   $\K > 0$, by our assumptions. 
Differentiating the Weierstrass representation for $u$, we have 
$$2 u_t = \frac{d}{dt} \Re \vector{ \int_0^z (A^2S+B^2S) \,dz} 
                                { i \int_0^z (A^2S-B^2S) \,dz}
                                { 2i \int_0^z ABS \,dz}
$$
Expressing this in terms of $\tAA$ and $\tBB$ we have by (\ref{eq:mm}) and (\ref{eq:mm2k})
\begin{equation}
2 u_t = \frac{d}{dt} \Re \vector{ t^{\mm}\int_0^\w (\tAA^2\tSS+ t^{2\K}\tBB^2\tSS) \,d\w  }
                                { i t^{\mm}\int_0^\w (\tAA^2\tSS-t^{2\K}\tBB^2\tSS) \,d\w}
                                { 2i t^{\mm+ \K} \int_0^\w \tAA\tBB\tSS \,d\w }
                                \label{eq:1409}
\end{equation}
The rule for differentiating with respect to $t$ in this situation is found by applying the chain rule as follows,
where $d/dt$ means the partial derivative with respect to $t$, holding $z$ fixed, and $\partial/\partial t$ 
means the partial derivative with respect to $t$, holding $\w$ fixed.
\begin{eqnarray*}
dH/dt &=& \frac{\partial H}{\partial t} + \frac{dH}{d\w} \frac{d\w}{dt} \nonumber \\
&=& \frac{\partial H}{\partial t} + \frac{ dt^{-\gamma} z}{dt} \frac{dH}{d\w} \nonumber  \\
&=&  \frac{\partial H}{\partial t} - \gamma t^{-\gamma-1}z \frac{dH}{d\w}\nonumber  \\
&=&  \frac{\partial H}{\partial t} - \gamma t^{-\gamma-1}t^\gamma w \frac{dH}{d\w}\nonumber \\
\end{eqnarray*}
The final result for differentiating with respect to $t$ is then
\begin{equation}
\label{eq:10}
\frac{dH}{dt} =  \frac{\partial H}{\partial t} - t^{-1}\gamma w \frac{dH}{d\w}
\end{equation}
Applying this rule on the right side of (\ref{eq:1409}), we find  
\begin{eqnarray*}
&& 2u_t = \\
&&   \Re \vector
            { t^{\mm-1}\mm\int_0^\w\tAA^2\tSS \,d\w + t^{\mm + 2\K -1}(\mm+2\K)\int_0^w \tBB^2\tSS \,d\w  }
            { it^{\mm-1}\mm\int_0^\w \tAA^2\tSS \,dw  + it^{\mm + 2\K-1}(\mm + 2\K )\int_0^\w \tBB^2\tSS \,d\w   }
            { 2it^{\mm + \K -1} (\mm + \K) \int_0^\w \tAA\tBB\tSS \,d\w  }    \\
 &&- \Re \vector
            { t^{\mm-1}\gamma \tAA^2\tSS\,  \w  + t^{\mm + 2\K -1} \tBB^2\tSS\,  \w }
            { it^{\mm-1}\gamma \tAA^2\tSS\, \w + it^{\mm + 2\K -1} \tBB^2\tSS\,  \w }
            { 2it^{\mm + \K -1}  \gamma \tAA\tBB\tSS \, \w}    \\
&&  + \Re \vector 
       {t^{\mm} \mm \frac \partial {\partial t} \int_0^\w \tAA^2\tSS \,d\w + 
        t^{\mm + 2\K}(\mm + 2\K) \frac \partial {\partial t} \int_0^w \tBB^2\tSS \,\w
      }
      { it^{\mm} \mm \frac \partial {\partial t} \int_0^\w \tAA^2\tSS \,d\w + 
        it^{\mm + 2\K} (\mm + 2\K)\frac \partial {\partial t} \int_0^\w \tBB^2\tSS \,d\w
      }
      { 2it^{\mm + \K }(\mm + \K) \frac \partial {\partial t} \int_0^\w \tAA\tBB\tSS \,d\w
      }  
\end{eqnarray*}

The terms in the third vector, involving $\partial/\partial t$, have one higher power of $t$ and are   $O(t^{\mm})$.  These swamp the $t^{\mm+2\K-1}$ terms.
Eliminating the swamped terms we have
\begin{eqnarray*}
&& 2u_t = \\
&&   \Re \vector
            { t^{\mm-1}\mm\int_0^\w\tAA^2\tSS \,d\w - t^{\mm-1}\gamma \tAA^2\tSS\,  \w  + O(t^{\mm})}
            { it^{\mm-1}\mm\int_0^\w \tAA^2\tSS \,dw -it^{\mm-1}\gamma \tAA^2\tSS\, \w  + O(t^{\mm}) }
            { 2it^{\mm + \K -1} (\mm + \K) \int_0^\w \tAA\tBB\tSS \,d\w-2it^{\mm + \K -1}  \gamma \tAA\tBB\tSS \, \w + O(t^{\mm+\K})}    
 \end{eqnarray*}
 We may also drop the tildes over $\AA$, $\BB$, and $\SS$, since 
the terms with nonzero powers of $t$ can be absorbed in the error terms.  Dropping the tildes
and dividing by the leading power of $t$, we have
\begin{eqnarray*}
\frac {2u_t} { t^{\mm-1}}  
&=&  \Re \vector
                          { \mm \int_0^\w \AA^2 \SS \,d\w - \gamma\AA^2\SS\,  \w  + O(t)}
                          { \mm i\int_0^\w \AA^2\SS \,d\w - i\gamma\tA^2\SS\, \w  + O(t)}
                          {  2it^{\K}(\mm + \K)  \int_0^\w \AA\BB\SS \,d\w -  \gamma\AA\BB\SS \, \w  + O(t^{\K + 1})} 
\end{eqnarray*}
Taking the dot product with $N$, we have the following formula for the eigenfunction, valid in the whole $\w$-plane,
including subsets that contain the zeroes of $\AA$:
\begin{eqnarray} \label{eq:ExactEigenFunction1}
&&\frac 2 { t^{\mm-1}}  \phi  = \\
&&  [ \None]\, \Re\, \bigg(\mm\int_0^\w  \AA^2\SS\,d\w - \gamma\AA^2\SS\, \w + O(t)\bigg)\nonumber \\
&&-   [ \Ntwo]\, \Im\, \bigg(\mm\int_0^\w \AA^2\SS\,d\w - \gamma\AA^2\SS\, \w + O(t) \bigg)  \nonumber \\
 && - t^{\K}(\mm + \K)[\Nthree]2\,\Im\,\bigg(\int_0^\w \AA\BB\SS \,d\w - \gamma \AA\BB\SS \, \w +
O(t) \bigg) \nonumber
\end{eqnarray}

\begin{lemma} [The limit of the eigenfunction]  \label{lemma:eigenfunctionCompact}
The following two formulas for the least eigenfunction of $u^t$ 
in terms of $N$ are valid on each compact subset of the $\w$ plane away from the zeroes of $\AA$,
provided $\K > 0$. 

\begin{eqnarray*} 
\frac {-1 } { t^{(\mm + \K)\gamma-1}}    \phi  &=& 
   \Im\,\bigg\{ \frac \BB \AA \, \bigg(\mm\int_0^\w \AA^2\SS\,d\w - \gamma\AA^2\SS\, \w \bigg) \\
&&  - (\mm + \K)\int_0^\w \AA\BB\SS \,d\w - \gamma\AA \BB \SS \, \w  \bigg\} + O(t) \\
&=&  \Im\,\bigg\{\frac \BB \AA \, \bigg(\mm\int_0^\w \AA^2\SS\,d\w \bigg) 
 - (\mm + \K)\int_0^\w \AA\BB\SS \,d\w  \bigg\}  + O(t)\\
\end{eqnarray*}
\end{lemma}

\noindent{\em Remark}.  The provision that $\K > 0$ is the same as providing that $N$ approach $(0,0,-1)$ on compact subsets 
away from the zeros of $\AA$.
It is certainly true when $\gamma = \gamma_1$, and later we will inductively show it is true on the inner rings of roots as well. 
\medskip

\noindent{\em Proof}. 
Starting with formula (\ref{eq:ExactEigenFunction1}),
we substitute for $\None$, $\Ntwo$, and $\Nthree$ the values obtained in (\ref{eq:Nofwexact}),
 to express the formula  entirely in  terms of $\tBB$ and $\tAA$.
We get (after multiplying by the denominator and changing the signs of both sides)
\begin{eqnarray*}
&& \frac {-(1 + t^{2\K} \vert  \tBB/\tAA\vert^2)} {t^{\mm-1}}    \phi  = \\
&&   \bigg(t^{\K} \Im\,(\tBB/\tAA)\bigg)\, \Re\, \bigg(\mm\int_0^\w \AA^2\SS\,d\w - \gamma\AA^2\SS\, \w + O(t)\bigg)  \\
 &&+  \bigg( t^{\K} \Re\,(\tBB/\tAA)\bigg)\, \Im\, \bigg(\mm \int_0^\w \AA^2\SS\,d\w - \gamma\AA^2\SS\, \w + O(t)\bigg)   \\
&& +  t^{\K\gamma}(\mm + 2\K)\bigg(-1+t^{2\K}\vert \tBB/\tAA\vert^2\bigg)2\,\Im\, \bigg(\int_0^\w \AA\BB\SS \,d\w -  \gamma\AA\BB\SS \, \w +O(t)\bigg) 
\end{eqnarray*}
Now we separate out the error terms:
\begin{eqnarray*}
&& \frac {-(1 + t^{2\K} \vert  \tBB/\tAA\vert^2)} {t^{\mm-1}}    \phi  = \\
&&   \bigg(t^{\K } \Im\,(\tBB/\tAA)\bigg)\, \Re\, \bigg(\mm\int_0^\w \AA^2\SS\,d\w - \gamma\AA^2\SS\, \w  \bigg)  \\
 &&+  \bigg( t^{\K } \Re\,(\tBB/\tAA)\bigg)\, \Im\, \bigg(\mm \int_0^\w \AA^2\SS\,d\w - \gamma\AA^2\SS\, \w  \bigg)   \\
&& +  t^{\K }(\mm + 2\K)\bigg(-1+t^{2\K\gamma}\vert \tBB/\tAA\vert^2\bigg)2\,\Im\, \bigg(\int_0^\w \AA\BB\SS \,d\w -  \gamma\AA\BB\SS \, \w \bigg) \\
&& + t^{\K  } \Im(\tBB/\tAA)\, O(t) + t^{\K  } (\BB/\AA)\, O(t) + t^{\K\ }\, O(t)
 + \vert \BB/\AA\vert^2 \,O(t^{3\K  } )\\
\end{eqnarray*}
Now we combine real and imaginary 
parts according to the rule $$\Re(u)\Im(v) + \Re(v)\Im(u) = \Im(uv),$$ applied to the first two terms on the 
right, and divide both sides by $t^{\K  }$.  That yields  
 \begin{eqnarray*} 
&& \frac {-(1 + t^{2\K } \vert \tBB/\tAA\vert^2)} { t^{\mm + \K -1}}    \phi  \ = \ 
    \Im\,\bigg\{(\tBB/\tAA)\, \bigg(\mm\int_0^\w \AA^2\SS\,d\w -\gamma \AA^2\SS\, \w  \bigg)\bigg\} \\
&& -  \bigg(1-t^{2\K }\vert \tBB/\tAA\vert^2\bigg) 2\, \Im\,\bigg\{(\mm + \K)\int_0^\w \AA\BB\SS \,d\w - \gamma\AA \BB \SS \, \w  \bigg\}  \\
&& + \Re\,(\tBB/\tAA) O(t)+ \Im\,(\tBB/\tAA) O(t) + \vert \tBB/\tAA \vert^2 O(t^{2\K})
\end{eqnarray*}
That formula for the least eigenfunction of $u^t$ 
in terms of $N$ is valid everywhere, even near the zeroes of $\AA$, in the sense 
 that the big-$O$ terms are analytic in $t$ and $w$.
 
Now replace
 $t^\K \tBB/\tAA$ with $O(t)$
on the left, and noting that $\tBB = \BB + O(t)$, and $\tAA = \AA + O(t)$, and  on a given compact subset away from zeros of $\AA$, we have $\tBB/\tAA = \BB/\AA + O(t)$.  That yields
the first formula of the lemma.
The second formula follows by canceling two terms in the first formula.  That completes the proof of the lemma.
\smallskip

We now define the function $\HH$ by 
\begin{eqnarray} \label{eq:HHdef}
\HH(\w) &=& \frac \BB \AA\, \mm\int_0^\w \AA^2\SS\,d\w  - (\mm + \K)\int_0^\w \AA\BB\SS \,d\w  \,  
\end{eqnarray}
Then by the previous lemma we have 
\begin{eqnarray*}
\Im\, \HH &=& \lim_{t \to 0} \frac {-1} { t^{\mm + \K -1}}    \phi  
\end{eqnarray*}
on compact subsets of the $\w$-plane away from the zeroes of $\AA$.
A crucial question is whether $\HH$ is constant or not.  If it is constant, then all the visible terms cancel out, and 
we can thus get no information from them.  If it is not constant, we can get some information by the asymptotic behavior of 
$\HH$ near the origin or near $\alpha_i$.

\section{Gaussian area and the eigenfunction}

We next draw some consequences from the calculations of the eigenfunction and the properties of $\HH$.
Recall that by Lemma~\ref{lemma:realAA},  all the $\alpha_i$ are real.  Hence we may guess that the 
$X$-component of $u_t$ does not change sign near $\alpha_i$.  Now consider the ``east pole''
and the ``west pole'' on the Riemann sphsere, namely $(1,0,0)$ and $(-1,0,0)$.  Can the unit 
normal take on both values near $\alpha_i$?  If it did, then since the eigenfunction $\phi = u_t \cdot N$,  that would imply that the $X$-component of $u_t$ would have to change sign near $\alpha_i$,
which seems unlikely.  Thus, we hope to show that $N$ cannot take on both the east pole and the 
west pole near $\alpha_i$.  But since $N$ is confined to be perpendicular to $\Gamma$ on the
boundary,  and $\alpha_i$ is real, it intuitively seems that $N$ must then be confined to one 
hemisphere.   These ideas are the basis of the following rigorous developments.

\begin{lemma} \label{lemma:hemisphericcovering}
Suppose that on compact subsets of the upper half of the $\w$-plane away from the zeroes of $\AA$, 
$N$ converges to $(0,0,-1,)$ as $t$ goes to zero. Then 
exactly one of the points $(1,0,0)$ and $(-1,0,0)$ (the ``east pole'' and ``west pole'') can be taken by $N$  at points 
 in the upper half plane in the  vicinity of one $\alpha_i$, and that point is taken on at most $J$ times, where $J$ is the number of 
 $\a_j$ that converge to $\alpha_i$.  
\end{lemma}

\noindent{\em Proof}.  Fix $i$, and fix a disk $U$ centered at $\alpha_i$ and excluding all other $\alpha_i$ and zero.
Consider a point $\xi$ in $U$ at which the normal takes on $(1,0,0)$ or $(-1,0,0)$.  Let $U^+$ be the intersection 
of $U$ with the upper half plane and let $U^-$ be the intersection of $U$ with the lower half plane. 
Since $N = t^\K \tBB/\tAA$, the hypothesis about the convergence of $N$ is equivalent to 
$\K > 0$.  Hence $N$ is $O(t)$ also on the boundary of $U$.  
There are $J$  values of $j$ for which $\a_j$ lies in $U$ (for sufficiently small $t$), and at these points $N$ takes
on the north pole.  Hence in the whole neighborhood $U$,  the map $N$ is a $J$-sheeted covering of the portion of the 
Riemann sphere lying north of a small circle of latitude within $O(t)$ of the south pole.  Hence each of the east pole
and the west pole can be taken on at most $J$ times over the whole neighborhood $U$.

We will show that if one of the east or west poles is taken on in $U^+$, the other cannot be.   At $\xi$ we have $\None = \pm 1$,
and $\Ntwo = \Nthree = 0$.  The formula (\ref{eq:ExactEigenFunction1}) calculated above for the eigenfunction $\phi$ in terms of the components of $N$ simplifies greatly when we substitute $\Ntwo = \Nthree = 0$.  According to that formula, at $\xi$ we have
\begin{eqnarray*}
\frac 2 {t^{\mm -1}  } \phi &=& \None\,\Re\,\bigg( \mm \int_0^\xi \AA^2\,\SS\,d\w - \AA^2(\xi)\,\SS(\xi)\,\xi + O(t)\bigg)
\end{eqnarray*}
Note that 
$$ \int_0^\xi \AA^2\SS\, dw = \int_0^{\alpha_i} \AA^2\SS\, dw  + \int_{\alpha_i}^\xi \AA^2\SS\, dw $$
The first term is a constant and the second is $O(\xi-\alpha_i)$.   Hence the previous equation can be written as
\begin{eqnarray*}
\frac 2 {t^{\mm -1}  } \phi &=&  \None \, c(1+O(\xi-\alpha_i)) - \tAA^2(\xi)\,\tSS(\xi)\,\xi + O(t)
\end{eqnarray*}
where $c$ is the constant $\mm \int_0^{\alpha_i} \AA^2\,\SS\,d\xi$. The constant $c$ is nonzero since (i) $\alpha_i$ is real, 
by Lemma \ref{lemma:realAA}, and
and (ii) the integrand $\AA^2\,\SS$ is nonnegative on the real axis, 
because $u$ takes the boundary monotonically and $u_x = t^{\mm} \AA^2\,\SS(1 + O(t))$, and (iii) 
$\alpha_i \neq 0$, since   $a_i$ goes to zero as $t^\gamma$.

  As $t$ goes to zero, $\xi$ will converge to $\alpha_i$, since by assumption 
outside every disk centered at $\alpha_i$, $N$ converges to $(0,0,-1)$.
We do not require that $\xi$ depend smoothly on $t$; a sequence of values of $t$ going to zero will suffice, or even
a single choice of $t$ depending appropriately on the constants hidden in the $O(t)$ and $O(w)$ terms.
Therefore the $O(\xi-\alpha_i)$ term goes to zero too; that is, $\xi -\alpha_i(t) = O(t)$. 
We have 
\begin{eqnarray}
\frac 2 {t^{\mm -1}  } \phi(\xi) &=&  \None(\xi) \, c + O(t)   \mbox{\qquad where $c = \mm \int_0^{\alpha_i} \AA^2 \SS$} \label{eq:phieastwest}
\end{eqnarray}

For small enough $t$, we see that the sign of $\phi(w)$ is equal to the sign of $\None$ if $\alpha > 0$, or is opposite to the sign 
of $\None$ if $\alpha < 0$.  In particular, for all such points in the vicinity of a fixed $\alpha_i$, the 
sign of $\None$ is the same, since the sign of $\phi$ is the same in the entire upper half plane, because $\phi$ is 
the eigenfunction of the least eigenvalue. 
That proves that not both the east and the west pole can be taken on over $U^+$ near $\alpha_i$.
But the west (or east) pole can be covered at most $J$ times, since $N$ is a $J$-sheeted covering over $U$.
That completes the proof.
\smallskip

{\em Discussion}. The significance of this lemma is that the Gauss map, considered as a covering map, covers at most one hemisphere
(plus $O(t)$ in a narrow band near the great circle dividing the two hemispheres, whose plane is perpendicular to the $X$-axis), since if it could be extended across the narrow strip near the $YZ$ plane where the image of the boundary
must lie, it would reach the other of the two points in question. Of course that hemisphere could be covered multiple times;
or possibly not at all, as it is conceivable that $\tilde a_i$ has a slightly negative imaginary part for $t>0$ and that $N$ then 
covers the entire sphere in the lower half plane near $a_i$, even though $\tilde a_i$ tends to the real number $\alpha_i$.
At this point we have not ruled 
out the possibility that one or more sheets of the covering might lie entirely over the lower half $U_-$ of $U$.
The argument does not work for $U_-$, because although $\phi$ continues analytically into $U_-$,   
 we do not know that it has the same sign throughout $U_-$.
 
However,  the proof does permit us to extract a bit more information than is stated in the above lemma.  For each $\alpha_i$,
we have shown that there is one ``permitted hemisphere'',  either the west hemisphere or the east hemisphere, according as 
the east pole or the west pole may possibly be a value of $N$ near that $\alpha_i$.  We can say more about which hemisphere
is permitted:

\begin{lemma} \label{lemma:permittedhemisphere}
All positive $\alpha_i$ have the same permitted hemisphere, and all negative $\alpha_i$ have the 
opposite permitted hemisphere.
\end {lemma}

\noindent{\em Proof}.  Let $\xi$ be a point near $\alpha_i$ where the east or west pole is taken on.  
The formula for $\phi(\xi)$ given in (\ref{eq:phieastwest}) shows that  $\phi(\xi)$ has the same sign as $\None(\xi)$
when $\alpha_i > 0$, since then $\int_0^{\alpha_i} \AA^2\SS\, dw > 0$,  and the opposite sign with $\alpha_i < 0$,
since then the integral is negative.  But $\phi$, since it is the eigenfunction for the least eigenvalue, has the 
same sign throughout the upper half plane.  Hence the permitted sign of $\None$ is the same for all positive 
$\alpha_i$, and the opposite sign for all negative $\alpha_i$.   That completes the proof.

\begin{definition}
We say that $\alpha_i$ {\em contributes Gaussian area} if, in the vicinity of $\alpha_i$, the 
normal restricted to a small upper half-disk covers an area at least $2\pi - O(t)$ on the sphere.
To make this definition precise, let $D_\rho$ be a disk of radius $\rho$ about $\alpha_i$ in the $w$-plane,
excluding all other $\alpha_j$ and zero, 
and let $D_\rho^+$ be its intersection with the upper half plane.  Then $\alpha_i$ contributes Gaussian
area if, for every $\rho > 0$,  $N(D_\rho^+)$ has area at least $2\pi - O(t)$ for sufficiently small $t$.
In case $N(D_\rho^+)$  has area $J\pi - O(t)$ then we say $\alpha_i$ contributes Gaussian area $J\pi$.
\end{definition}

Since the normal $N$ is a projection mapping (away from its ramification points), and since the 
image of the boundary is confined to within $O(t)$ of the $YZ$-plane,  as soon as $N$ takes on a value
more than $O(t)$ from the $YZ$-plane (in $D_t^+$)  then $\alpha_i$ contributes Gaussian area, and in that
case it must contribute Gaussian area $2J\pi$ for some $J$, i.e. area equal to an integer number of hemispheres.
  In particular, if $N$ takes on the ``east pole'' $(1,0,0)$ or the ``west pole'' $(-1,0,0)$ in $D^+_t$
for every $t > 0$, then $\alpha_i$ contributes Gaussian area.

{\em Discussion}.  How much Gaussian area is contributed by each $a_i$?  We hope to show just 
one hemisphere, rather than one sphere.  We know that no $a_i$ can contribute a whole sphere
in the most obvious way, as then the Gaussian image would contain more than a hemisphere and the 
eigenvalue would be less than 2.  Moreover in the previous lemma we showed that the image of the 
boundary sticks up at least to the equator, hinting that it divides the sphere into hemispheres.
But this still does not rule out the possibility that somehow a whole sphere of Gaussian area
might be contributed.  In the end, we are going to count the Gaussian area, and it will be 
vital that each $a_i$ contributes only a hemisphere.  The next lemma addresses this crucial issue.
The key to the proof is the fact that near $a_i$, and in the upper half-plane, $N$ can take on the 
east or the west pole, but not both.  That means it is confined to a hemisphere.  Here are the details:
\smallskip

\begin{lemma} \label{lemma:GaussAreaPrincipal}
Assume that $N$ goes to $(0,0,-1)$ on compact subsets of the $w$-plane away from the zeroes of $\AA$.
Then every $\alpha_i$ contributes Gaussian area at most equal to one hemisphere for each $\a_j$
that converges to $\alpha_i$. 
\end{lemma}

\noindent{\em Proof.}  According to the Gauss-Bonnet theorem (discussed in Section \ref{section:Gauss-Bonnet}), we 
have when $t=0$
$$  \int KW \, dx\, dy + \int_\Gamma \kappa_g = 2\pi + 2M \pi$$
where $M$ is the sum of the orders of the boundary branch points and twice the sum of the 
orders of the interior branch points.   When $t>0$, there are no branch points in the closed upper half plane,
so we have 
$$  \int KW \, dx\, dy + \int_\Gamma \kappa_g = 2\pi  $$
We let $t$ go to zero in this formula.  
Since $\kappa_g$ is bounded by the curvature of $\Gamma$,  the integral containing $\kappa_g$ converges as $t$ goes
to zero to the corresponding term in the first formula.   Hence the first term, which represents the Gaussian area, 
must jump upwards by $2M\pi$.  That is, each interior branch point contributes one sphere  (times the order)
and each boundary branch point contributes a hemisphere (times the order).   

Near each $a_i$,  disregarding the boundary,  one whole sphere of Gaussian area (times the order) will exist in a small
neighborhood of $a_i$ for  small positive $t$, within $O(t)$.   Since the Gaussian image of the boundary is confined
to lie near the $YZ$ plane,  and since (in the case of several nearby $a_i$'s) the Gauss map may be a multi-sheeted 
covering of the sphere,  there are ({\em a priori}) several possibilities.  For example, there may be no Gaussian area contributed
(in the upper half plane), if the whole ``bubble'' is taken on in the lower half plane.  Or, the ``bubble'' may 
be divided by the image of the boundary, so that one hemisphere is taken on in the upper half plane, and one in the 
lower half plane.  Or, if there is more than one sheet of the covering in a small neighborhood of $\alpha_i$,  which
is possible if several $\a_j$ converge to the same $\alpha_i$,  then some sheets may be divided by the boundary
and others not.   {\em A priori} it is also possible that the image of the boundary may rise part way up the sphere and 
then reverse direction, so as to only partly divide the sphere.

However, we can now rule out some of these possibilities.  Let $p$ be one of the $\alpha_i$ or zero.
Fix a neighborhood $U$ of $p$ in the $w$-plane
small enough to exclude all other $\alpha_i$, or all $\alpha_i$ if $p$ is zero, 
and let $U^+$ be the intersection of $U$ with the upper half-plane.
We know that each $\alpha_i$ is real, so $U^+$ is half of $U$.  Let $J$ be the number of $\a_j$ that converge to $p$.
Then the total extra Gaussian area over $U$ is $4\pi J$, and over $U$, the Gauss map is a $J$-sheeted 
branched covering of the part of the Riemann sphere north of any fixed circle of latitude (say in the southern 
hemisphere).   That is, fixing such a circle of latitude, for $t$ sufficiently small $N$ will be a branched covering
of the upper part of the Riemann sphere, north of the fixed circle of latitude. 
 
  If this $\alpha_i$ contributes any Gaussian area at all (in the upper 
half plane),  then $N$ takes on either the east or west pole somewhere in $U^+$.  Let us suppose that $N$ takes on 
the west pole $W$ (the other case is symmetric).   It may do so on more than one sheet.
On each sheet (that is, starting from each different pre-image of $W$,  we can continue a path $\pi$ such 
that $N(\pi(\tau))$ lies on the equator (even if the path must path through ramification points of $N$).  But the image $N(\pi(\tau))$ cannot extend into the eastern hemisphere by more than $O(t)$,
or else it could be continued all the way to the east pole, contradicting Lemma~\ref{lemma:permittedhemisphere}. 
  Over the whole neighborhood $U$ there are  $J$
sheets covering the sphere $J$ times, for a total area of $4\pi J$.  But now we have proved that each 
of these $J$ sheets has at most the west hemisphere covered by $N(U^+)$.

 Hence the total Gaussian area contributed by $\alpha_i$ 
is at most $2\pi J$; that is, one hemisphere per $a_i$ at the most.  
That completes the proof of the lemma.


\section{What is true if $H$ is not constant}

In this section, we suppose that $\gamma$ is one of the $\gamma_n$, that $z = t^\gamma \w$, 
and that $N$ goes to $(0,0,-1)$ on compact subsets away from the zeroes of $\AA$.  This condition 
we call the ``$N$-condition''.   Since 
$$ \frac B A = t^{\K} \frac \BB \AA$$
the $N$-condition can be succinctly expressed as $\K > 0$.
The $\alpha_i$ are the coefficients of the $a_i$ in this ring of roots, 
i.e.,  $a_i = \alpha_i t^\gamma ( 1 + O(t))$,  and we define $\alpha_i = 0$ for 
roots that go to zero faster than $t^\gamma$.  Similarly for the $\beta_i$ and $\zeta_i$.
We refer to the roots that go to zero slower than $t^\gamma$ as ``slow roots'', and 
the roots that go to zero faster as ``fast roots''.

\begin{lemma} \label{lemma:H}
 If\, $\HH$ is not constant, then there are no fast roots; that is, no $a_i$ or $b_i$ or $s_i$ goes to zero faster than $t^\gamma$.  
\end{lemma}

\noindent{\em Proof}.  Recall that $Q$, $S$, and $R$, are the numbers of $a_i$, $s_i$, and $b_i$
going to zero as fast or faster than $t^\gamma$.  But now we need notation for the numbers of 
fast roots. 
\def\Qf{Q_f}
\def\Sf{S_f}
\def\Rf{R_f}
 Let $\Qf$, $\Sf$, and $\Rf$ be respectively the number of fast $a_i$, fast $s_i$, and fast $b_i$.
 The total number of   roots going to zero faster than $t^\gamma$, is $$ P := \Qf + \Rf + \Sf.$$
Let $\alphat$ be the product of the nonzero $\alpha_i$ (all of which are real) over the $a_i$
that go to zero as $t^\gamma$.  Let
$\betat$ be the product of the nonzero $\beta_i$.
If the $\beta_i$ are not real, they come in complex conjugate pairs, by Lemma \ref{lemma:realAA}, 
so the product $\betat$ is real. Let $\zetat$ be the product of the nonzero $\zeta_i$, which is 
real since all the $\zeta_i$ are real.  (Do not confuse $\alphat$  with 
$\alpha$,  which is the product of the $\alpha_i$ over the slow roots.)
\
 We will analyze the asymptotic form of $\HH$ near the origin (without yet assuming $\HH$ is not constant.) 
 We will show that, if
$\HH$ is not constant, its
leading term is a nonzero real constant times $ \beta \w^{P+1}$.

 We  analyze the first formula of Lemma \ref{lemma:eigenfunctionCompact}, rather than the 
second one.  Letting $t$ go to zero in that formula and using the definition of $\HH$ we have
\begin{eqnarray} 
&&  \HH =    \frac \BB \AA \, \bigg(\mm\int_0^\w \AA^2\SS\,d\w - \gamma \AA^2\SS\, \w \bigg)  \label{eq:25} \\
&&  
 - (\mm + \K)\int_0^\w \AA\BB\SS \,d\w - \gamma\AA \BB \SS \, \w   \nonumber
\end{eqnarray}
 Consider the second term in that equation, namely
\begin{equation} \label{eq:secondterm23}
(\mm + \K)\int_0^\w \AA\BB\SS\, d\w - \gamma \AA\BB\SS \,\w.
\end{equation}
Asymptotically for small $w$, $\AA\BB\SS$ is a 
constant times $w^P$, since each fast root contributes a factor of $w$ and there are $P$ fast roots in all.
  Hence, asymptotically for small $\w$, the second term (\ref{eq:secondterm23})    has the form
\begin{eqnarray}
(\mm + \K)\int_0^\w \AA\BB\SS\, d\w - \gamma \AA\BB\SS \,\w &=& cw^{P+1} + O(w^{P+2}) \label{eq:secondterm}
\end{eqnarray}
It does not matter whether the coefficient $c$ is zero or not (although a closer analysis 
shows that it is not).

Now consider the behavior of the first term in (\ref{eq:25}) near the origin, namely
\begin{equation}
\label{eq:27}
\frac{\BB}{\AA} \bigg(\mm\int_0^w \AA^2\SS \, d\w - \gamma\AA^2\SS\, \w\bigg).
\end{equation}
We have 
\begin{eqnarray*}
\AA  &=& \alpha \alphat w^{\Qf} + O(w^{\Qf+1}) \\
 \SS &=& \zeta\zetat w^{\Sf} + O(w^{\Sf+1}) \\
 \AA^2 \SS &=& \alpha^2\alphat^2\zeta\zetat w^{2\Qf + \Sf} + O(w^{2\Qf+\Sf+1}) 
 \end{eqnarray*}
 To simplify this equation we introduce a new constant:
 $$ c_1 = \alpha^2 \alphat^2\zeta\zetat.$$
 Then we have
 
 \begin{eqnarray*}
 \AA^2 \SS &=& c_1 w^{2\Qf + \Sf} + O(w^{2\Qf+\Sf+1}) \\
 \mm \int_0^w \AA^2 \SS \, d\w &=& \frac {\mm c_1} {2\Qf+\Sf+1} w^{2\Qf+\Sf+1}  + O(w^{2\Qf+\Sf+2})\\
 \gamma\AA^2 \SS \w &=& \gamma c_1 \w^{2\Qf+\Sf+1} + O(w^{2\Qf+\Sf+2}) 
\end{eqnarray*}
Subtracting the last two equations, 
\begin{eqnarray*}
  \mm \int_0^w \AA^2 \SS \, d\w -\gamma\AA^2 \SS \w &=& c_1\bigg( \frac \mm {2Q+S+1} - \gamma\bigg) w^{2\Qf+\Sf+1}  + O(w^{2\Qf+\Sf+2})
\end{eqnarray*}
We have 
\begin{eqnarray*}
 \frac {\BB}{\AA} &=& \frac { \BB} {\alpha\alphat  w^{\Qf} + O(w^{\Qf+1})} \\
 &=& \frac {\beta\betat}{\alpha\alphat} w^{\Rf-\Qf} (1 + O(w))\\
 &=& \frac {\beta\betat}{\alpha\alphat} w^{\Rf-\Qf} + O(w^{\Rf-\Qf+1})
 \end{eqnarray*}
so the first term of (\ref{eq:25}), that is (\ref{eq:27}), comes to 
\begin{eqnarray}
&&c_0\bigg( \frac \mm {2\Qf+\Sf+1} - \gamma\bigg) w^{\Rf+\Qf+\Sf+1}  + O(w^{\Rf+\Qf+\Sf+2}) \nonumber\\
&=& c_0\bigg( \frac \mm {2\Qf+\Sf+1} - \gamma\bigg) w^{P+1} + O(w^{P+2}) \label{eq:firstterm}
\end{eqnarray} 
Subtracting the contributions of the two terms as given in (\ref{eq:firstterm}) and (\ref{eq:secondterm}), we have\begin{eqnarray*}
  \HH  &=&  w^{P+1}(c_0 \bigg( \frac \mm {2\Qf + \Sf+1} - \gamma\bigg)  - c)   +O(w^{P+2})  \\
  &=& O(w^{P+1})
\end{eqnarray*}
It does not matter whether the leading term is zero or not.
$\HH$ is a rational function of $w$, whose poles are at the $\alpha_i$, by definition all different
from zero.  If $\HH$ is not identically zero,  then it has a leading term with power at least
$P+1$;  but since $\Im\, \HH$ takes only one sign in the upper half plane, we must then have 
$P=0$.   But $P=0$ means there are no fast roots.  That completes the proof of the lemma.
\medskip

{\em Remark}.  A similar analysis near each $\alpha_i$ yields the following lemma,
which we do not prove here because it turns out not to be needed.  Nevertheless, 
stating it may assist the reader to understand the situation.
\medskip

\begin{lemma} \label{lemma:H2}   
If $H$ is not constant, and the $N$-condition holds, then at each $\alpha_i$,  either $\BB/\AA$ has a simple zero or a simple pole, or is 
analytic non-vanishing at $\alpha_i$.
\end{lemma}

\noindent{\em Remark}. This means that the number of $\b_j$ converging to $\alpha_i$  and the number of $\a_i$ converging to $\alpha_j$
differ by at most 1.
\smallskip

\noindent{\em Proof}.  Omitted, because the lemma is not needed.

 \section{What is true if $\HH$ is constant}
In this section we draw out some consequences of the assumption that $\HH$ is constant.
We start with a formula for the derivative of $\HH$:

\begin{lemma} \label{lemma:Hdif}
$$\frac d {dw} \HH = \frac d {dw} \bigg( \frac \BB \AA \bigg) \int_0^\w \AA^2\SS\,d\w
-\K \AA\BB\SS $$
\end{lemma}

\noindent{\em Proof}. By the definition  in (\ref{eq:HHdef}) we have 
\begin{eqnarray*}  
\HH(\w) &=& \frac \BB \AA\, \mm\int_0^\w \AA^2\SS\,d\w  - (\mm + \K)\int_0^\w \AA\BB\SS \,d\w  \,  
\end{eqnarray*}
Differentiating we have
\begin{eqnarray*} 
\frac d {dw} \HH  &=& \frac d {dw} \bigg( \frac \BB \AA \bigg) \int_0^\w \AA^2\SS\,d\w
                       + \bigg(\frac \BB \AA\bigg) \AA^2 \SS - (\mm +\K) \AA \BB \SS
\end{eqnarray*}
Simplifying, we have the formula of the lemma.  That completes the proof.

\begin{lemma} \label{lemma:Hconstant}
If $H$ is constant, then 
$$ \frac \BB \AA = C \bigg(\int\,\AA^2\SS\,d\w\bigg)^{ \K/\mm}$$
for some constant $C$,  
and $\BB/\AA$ is a polynomial.  
\end{lemma}

\noindent{\em Proof}.  
 Let $\Q$ be the rational function $\BB/\AA$. 
By Lemma~\ref{lemma:Hdif}, we have
\begin{eqnarray*}
\HH^\prime = 0   
&=& \mm \Q^\prime \int\, \AA^2\SS\,d\w - \K\AA^2\SS\Q   
\end{eqnarray*}
The equation has one solution 
$$\Q = \bigg(\int\,\AA^2\SS\,d\w\bigg)^{ \K/\mm}$$
Since it is a linear equation, every solution is a multiple of that one, so the formula in the lemma 
is proved.  
Whenever $C$ occurs in the context of this integral, it will mean this constant.

Multiplying by $\AA$ we have 
$$\BB  =  C \AA \bigg(\int\,\AA^2\SS\,d\w\bigg)^{ \K/\mm}.$$
This equation shows that $\BB$ vanishes wherever $\AA$ does, and to at least
as great an order.  Hence the rational function $\BB/\AA$ is actually a polynomial.
That completes the proof of the lemma.
\medskip

The converse is also true:
\begin{lemma} \label{lemma:impliesHconstant}
If for some constant $C$, we have
$$ \frac \BB \AA = C \bigg(\int\,\AA^2\SS\,d\w\bigg)^{ \K/\mm}$$
then $\HH$ is constant.
\end{lemma}

\noindent{\em Proof}.  By Lemma~\ref{lemma:Hdif}, 
\begin{eqnarray*}
\frac d {dw} \HH   &=& \mm \frac d {dw} \bigg( \frac \BB \AA \bigg) \int\, \AA^2\SS\, d\w   - \K\AA\BB\SS \\
\end{eqnarray*}
Substituting the expression assumed in the lemma for $\BB/\AA$ we have
\begin{eqnarray*}
\frac d {dw}\HH   &=& \mm C \frac d {dw}\bigg( \int\,\AA^2\SS\,dw \bigg)^{\K/\mm}\int\, \AA^2\SS\, d\w   - \K\AA\BB\SS \\
&=& \K C\bigg( \int\,\AA^2\SS\,dw \bigg)^{\K/\mm-1}\AA^2 \SS \int\, \AA^2\SS\, d\w   - \K\AA\BB\SS \\
&=&\K C\bigg( \int\,\AA^2\SS\,dw \bigg)^{\K/\mm} \AA^2\SS   - \K\AA\BB\SS \\
&=& \K \bigg(\frac \BB \AA \bigg) \AA^2 \SS - \K\AA\BB\SS \\
&=& 0
\end{eqnarray*}
Since $\HH$ has derivative 0,  $\HH$ is constant.  That completes the proof of the lemma.

\begin{lemma} \label{lemma:H-constant} 
Let $R$, $Q$, and $S$ again be the numbers of $b_i$, $a_i$, and $s_i$ going to zero as $t^\gamma$ or faster. (In other words, the degrees of $\BB$, $\AA$, and $\SS$ respectively.) Let 
$R_f$, $Q_f$, and $S_f$ be the numbers of $b_i$, $a_i$, and $s_i$ going to zero strictly faster than $t^\gamma$ (the ``fast roots'').
If \,$\HH$ is constant then
\smallskip

(i)  $ R-Q = (2Q+S+1) \K/\mm$.
\smallskip

(ii) $R_f - Q_f = (2Q_f+S_f+1) \K/\mm$.
\smallskip

(iii) if $\K > 0$ then $R_f > Q_f$.
\end{lemma}

\noindent{\em Proof}. 
 By the previous lemma we have 
$$ \frac \BB \AA = C\bigg(\int\,\AA^2\SS\,d\w\bigg)^{ \K/\mm}$$
We will derive (i) and (ii) by looking at this formula for large $w$ and small $w$, 
respectively.   Assertion (iii) follows immediately from (ii), so it suffices to 
prove (i) and (ii).
\smallskip

 Ad (i).  We look at the asymptotic expansion for large $w$.
 We note that as $\w$ goes to infinity, $\Q = \BB/\AA$ must be asymptotic to $  (\beta/\alpha)\w^{R-Q}$, where
$R$ and $Q$ are degrees of $\BB$ and $\AA$ respectively.   The right hand side gives us
\begin{eqnarray*}
 \bigg(\int\,\AA^2\SS\,d\w\bigg)^{ \K/\mm} &=& 
 \bigg(  \frac {\alpha^2\zeta}{2Q+S+1}( w^{2Q+S+1} + O(w^{2Q+S}))\bigg)^{\K/\mm} 
\end{eqnarray*}
\begin{eqnarray*}
 &=& \bigg(\frac {\alpha^2\zeta}{2Q+S+1}\bigg)^{\K/\mm} w^{(2Q+S+1)\K/\mm} + O(w^{(2Q+S+1)\K/\mm -1})
 \end{eqnarray*}
The exponents of the leading power of $w$ must be equal; that yields the equation in (i).  Comparing the coefficients 
also enables us to determine the constant $C$; it turns out to be 
$$ C = \frac { \beta}{\alpha} \bigg(\frac{2Q+S+1}{\alpha^2\zeta}\bigg)^{\K/\mm} $$

Ad (ii).  Now we look at the behavior for small $w$.
Asymptotically near the origin we have 
$$ \int\,\AA^2\SS\,d\w = \frac {\alpha^2 \zeta} {2Q_f+S_f+1} \w^{2Q_f + S_f + 1} (1 + O(\w))$$
and hence 
$$\bigg(\int\,\AA^2\SS\,d\w\bigg)^{ \K/\mm} =  \bigg( \frac {\alpha^2\zeta} {2Q_f+S_f+1}\bigg)^{\K/\mm} \w^{(2Q_f + S_f + 1)\K/\mm} (1 + O(\w))$$
and hence 
\begin{eqnarray*}
 \BB  &=& C \bigg( \frac {\alpha^2\zeta} {2Q_f+S_f+1}\bigg)^{\K/\mm} \AA\, \w^{(2Q_f + S_f + 1)\K/\mm} (1 + O(\w)) \\
&=&   c  \w^{(2Q_f + S_f + 1)\K/\mm + Q_f} (1 + O(\w))  \qquad\mbox{for some constant $c \neq 0$} \\
\beta \w^{R_f} &=&  c \w^{(2Q_f + S_f + 1)\K/\mm + Q_f} (1 + O(\w)) \\
\end{eqnarray*}
Equating the exponents of the leading terms on left and right we have 
$$ R_f = (2Q_f+S_f+1)\K/\mm + Q_f$$
as required in (ii). 
 That completes the proof of the lemma.
\smallskip


\section{Descending through the rings of roots}
Recall that the roots $a_i$, $b_i$, and $s_i$ fall into ``rings'',  where the $n$-th ring of roots
goes to zero as $t^{\gamma_n}$, and $\gamma_1 > \gamma_2 \ldots$.  For each ring we can 
consider the corresponding $\w$-plane, where $z = t^{\gamma_i} \w$.   In each ring there 
is a different $\AA$, $\BB$, $\SS$, and $\HH$. Up until now, we have considered only one 
ring of roots at a time.  Now, we must consider more than one ring at a time.  

The
  results of this section are (i) the branch points $s_i$ do not occur, and (ii) 
$\HH$ is constant on all the rings but the innermost, and (iii) $\HH$ is not constant on 
the innermost ring, and (iv) there is at least one $\alpha_i$ on that innermost ring.
 
  As above, we abbreviate by the phrase ``$N$-condition''
the proposition that $N$ goes to $(0,0,-1)$ on compact subsets of the $w$-plane apart from 
the zeros of $\AA$.  The $N$-condition depends on $n$, since the power of $t$ used to 
define the $w$-plane for the $n$-th ring of roots is $\gamma_n$.  If the $n$-th ring satisfies
the $N$-condition, and $\HH$ for ring $n$ is constant, then by 
Lemma \ref{lemma:H-constant}, there are more fast $b_i$'s than there are fast $a_i$'s,  so the $n+1$-st ring
also will satisfy the $N$-condition.   This is the fundamental observation behind the 
following lemma:

\begin{lemma}  \label{lemma:descendingrings} 
 Suppose that the $N$-condition holds on all rings from the first up to and including the $n$-th ring,
and that $\HH$ is constant on all those rings.  Then
the $N$-condition  holds on ring $n+1$. 
\end{lemma}

\noindent{\em Proof}.  The essential point is that on each of the first $n$ rings, there are 
more $b_i$ roots than $a_i$ roots, as we shall see.

 The $N$-condition is equivalent to  $\tBB /\tAA = O(t)$ on compact subsets  of the $w$-plane
  away from the $\alpha_i$, and also (by Lemma~\ref{lemma:Kpositive}) to $\K > 0$.
Fix one of the rings,  say the one given by $\gamma_j$,  and let $\gamma = \gamma_j$.
Recall that when we set $z = t^\gamma w$,  and convert $z-b_i$ to the $w$-plane,
we get $t^\gamma(w-\beta_i)$ if $b_i$ is on the $j$-th ring,  or if $b_i$ goes to
zero as $t^{\gamma_i}$ with $\gamma_i < \gamma$,   we get $t^\gamma_i$ times a constant,
or if $b_i$ goes to zero faster, we just get $t^\gamma w$.   The power of contributed
by the slow roots is $t^{\gamma_i}$ (whatever $\gamma$ is),  so if there are more $b_i$ than $a_i$
on every ring slower than $\gamma$,  then when $B/A$ is expressed in the $w$-plane,
we get a positive power of $t$ from the slow roots.   Each root going to zero as $t^\gamma$
or faster contributes a power of $t^\gamma$, so if the degree of $\BB$ is more than 
the degree of $\AA$, we will get a positive power of $t$ from those roots too.

Now suppose that for all rings up to and including the $j$-th one, there are more $b_i$ roots than 
$a_i$ roots on that ring.  
To make the argument precise requires numerous subscripts, but it seems unavoidable:

Let $R_n$, $Q_n$, and $S_n$ be the number of roots $b_i$, $a_i$ and $s_i$
 that go to zero as or faster than $\gamma_n$; 
in other words the degrees of $\BB$, $\AA$, and $\SS$ respectively.
Let $e_n$ be the sum of $n_i \gamma_i$ for $i < n$, where $n_i$ is the number of $b_i$ that go to zero as $t^\gamma_i$;
and let $d_n$ be the sum of $n_i \gamma_i$ for $i < n$, where $n_i$ is the number of $a_i$ that go to zero as $t^\gamma_i$.
Then in the $n$-th ring we have
\begin{eqnarray*} 
\frac B A &=&  t^{\K} \frac \BB \AA  + O(t^{\K + 1})\\
&=& t^{e_n-d_n + \gamma(R_n-Q_n)} \frac \BB \AA + O(t^{\K + 1})
\end{eqnarray*}
Since the $N$-condition holds for ring $n$, we have
$$ e_n - d_n + \gamma(R_n-Q_n) > 0.$$
When we pass to ring $n+1$, the roots that go to zero slower than $t^\gamma_{n+1}$ include those that 
contributed $t^{e_n - d_n}$ at ring $n$, but also those that go to zero as $t^{\gamma_n}$.  Those now
contribute to $t^{e_{n+1}-d_{n+1}}$ the powers of $t$ that they were contributing to $B/A$, making
no net change in the power of $t$ in $B/A$.  The remaining roots, the roots that go to zero
faster than $t^{\gamma_n}$,  now contribute a factor $t^{\gamma_{n+1}}$ instead of the former
$t^{\gamma_n}$.   

 Therefore the $N$-condition will hold at ring $n+1$ if there are more of the
$b_i$ going to zero faster than $t^{\gamma_n} $ than there are $a_i$ going to zero faster than
$t^{\gamma_n}$.  But that has been proved in Lemma~\ref{lemma:H-constant}, part (ii),
 which is applicable since 
$\HH$ is constant on ring $n$, and the $N$-condition is equivalent to $\K > 0$.
That completes the proof of the lemma.  (Actually $\K \ge 0$ in Lemma~\ref{lemma:H-constant}, part (ii),
 would have been enough.)

\begin{lemma} \label{lemma:Hnotconstant} On the ``innermost'' ring of roots, i.e. the roots that go to zero as the 
highest power of $t$,  $\HH$ is not constant, and the $N$-condition holds (i.e. $\K > 0$) on that ring, and 
indeed on all the rings.   On all but the innermost ring $\HH$ is constant.  On all the rings we have 
$$ \frac {\K} \mm = \frac k {2m+1} $$ 
\end{lemma}

\noindent{\em Proof.} There are only finitely many rings, because there are only finitely many of the $a_i$, $b_i$, 
and $s_i$.   Consider the first one on which $\HH$ is not constant, if any.
By Lemma \ref{lemma:descendingrings},
each ring of roots down to and including that ring satisfies the $N$-condition, i.e. $\K > 0$.
 On a ring where $\HH$ is constant and the $N$-condition is satisfied, there are some fast roots, by Lemma \ref{lemma:H-constant}, part (iii).
 Assume, for proof by contradiction, that $\HH$ is constant on all the rings.  Then the $N$-condition holds on the 
 innermost ring.  But by Lemma \ref{lemma:H-constant}, part (iii) there are some fast roots on that ring, which is impossible
 for the innermost ring, as the fast roots would lie on another ring of roots.  That contradiction shows that $\HH$
 cannot be constant on all the rings, so there is a ring on which $\HH$ is not constant.  By Lemma \ref{lemma:descendingrings},
 the $N$-condition holds
 on the slowest (outermost) ring for which $\HH$ is not constant.   Then by Lemma \ref{lemma:H}, there are no 
 fast roots on that ring.  Hence, it is the innermost ring.  It now remains only to prove
 the last equation of the lemma.
 
 On the outermost ring we have $\K = k$ and $\mm = 2m+1$, so $\K/\mm = k/(2m+1)$ on the outermost ring.  Now look at 
 the formulas in Lemma \ref{lemma:H-constant}.  These involve the numbers $Q$, $R$, and $S$ (numbers of roots that 
 go to zero as $t^\gamma$ or faster) and $Q_f$, $R_f$, and $S_f$ (numbers of roots that go to zero faster than $t^\gamma$).
  When we 
 pass from a ring on which $\HH$ is constant to the next ring, the $Q_f$, $R_f$, and $S_f$ of the first ring become 
 the $Q$, $R$, and $S$ on the next ring in.  According to Lemma \ref{lemma:H-constant},  we have 
 $\K/\mm$ on the first ring equal to $(R_f -Q_f)/(2Q_f+S_f+1)$.   Now passing to the second ring, the values on the right 
 become $(R-Q)/(2Q+S+1)$, which, according to the other formula in Lemma \ref{lemma:H-constant}, is equal to 
 $\K/\mm$ on the next ring.   Hence, the value of $\K/\mm$ never changes as we pass from one ring to the next,
 all the way down to the innermost ring.  Hence, it retains the value $k/(2m+1)$ that it has on the outermost ring.
  That completes the proof of the lemma. 

\begin{lemma} \label{lemma:NoS} The $s_i$ do not occur.
\end{lemma}

\noindent{\em Proof}.   According to the Gauss-Bonnet-Sasaki-Nitsche formula, a total of $2\pi m$ in extra 
Gaussian area must be contributed by the $a_i$, since the branch point has order $2m$,
and each boundary branch point contributes an amount equal to a hemisphere of Gaussian
area, times half the order.   If any $s_i$ exist, there are 
fewer than $2m$ of the $a_i$. 
We shall show that each $a_i$ can contribute at most a hemisphere, so the total will be insufficient unless there are actually $m$ of the $a_i$, in which case, no $s_i$ exist.   By Lemma \ref{lemma:descendingrings},
on the innermost ring $\HH$ is not constant, and the $N$-condition holds on each ring.   
 By Lemma \ref{lemma:hemisphericcovering}, all the $a_i$ that belong to that ring or to slower rings
 contribute only one hemisphere of Gaussian area per $a_i$.  Thus every $a_i$ contributes only one hemisphere.
That completes the proof of the lemma.
\smallskip

\noindent{\em Remark.}  Tony Tromba asked whether the proof (of our main theorem about finiteness)
 requires an appeal to external references to know that 
the branch point is not a false branch point.  It does not, and the key reason why not is the application of the Gauss-Bonnet 
formula.  The surfaces $u^t$ for $t > 0$ must contribute Gaussian area.  Therefore $\HH$ is not constant, as we have shown.
That means that the surfaces $u^t$ are not simply ``sliding'' over the same geometrical surface as $t$ becomes positive; there is 
a normal component, with {\em some} power of $t$, to the motion. 

\begin{corollary} \label{lemma:HwithoutS}
 We have
$$ \HH =  \frac \BB \AA \, \bigg(\mm\int_0^\w \AA^2 \,d\w \bigg) 
 - (\mm + \K)\int_0^\w \AA\BB  \,d\w   
 $$
and $\Im\,\HH$ cannot take two different signs in the upper half plane.
\end{corollary}

\noindent{\em Proof}.   By the previous lemma we can now put $\SS = 1$ in the formula
\ref{eq:HHdef}.
  Since $\Im\,\HH$ is (according to  Lemma~\ref{lemma:eigenfunctionCompact}) the limit of a power of $t$ times
the first eigenfunction,  and the first eigenfunction has just one sign in the upper half plane,
either $\Im\, \HH$ takes just one sign or possibly is identically zero.  That completes the proof.

\begin{lemma} \label{lemma:crucialidentity}
Suppose that on the outermost ring of roots, $\HH$ is constant.  Then on the innermost ring,  we have 
$ R-Q = (2Q+1) \K/\mm$, where $Q$ and $R$ are respectively the degrees of $\AA$ and $\BB$;  that is, the number
of $a_i$ and $b_i$ belonging to the innermost ring. 
\end{lemma}

\noindent{\em Proof}.  If on the outermost ring $\HH$ is constant, then the innermost ring (where $\HH$ is not constant)
is not the outermost ring.   Consider the next-to-innermost ring.   On that ring, the ``fast roots'' are the 
roots that belong to the innermost ring.  We apply Lemma \ref{lemma:H-constant} on the next-to-innermost ring,
where $\HH$ is constant.   Formula (ii) of that lemma tells us that 
$$ R_f - Q_f = (2Q_f + S_f + 1) \frac \K \mm $$
where the numbers $Q_f$, $R_f$, and $S_f$ are the numbers of $a_i$, $b_i$, and $s_i$ that go to zero faster
than the roots of the next-to-innermost ring.  But those are exactly the roots of the innermost ring.  By Lemma \ref{lemma:NoS},
we have $S_f = 0$.  Hence 
$$ R - Q = (2Q + 1) \frac \K  \mm$$
That completes the proof of the lemma.
 
\begin{lemma} \label{lemma:someA} [This claim is withdrawn.  The step $e \le d$ near the end of the proof is 
wrong.]
On the innermost ring of roots,  there exists some $a_i$; that is, 
some $a_i$ goes to zero as $t^\gamma$, where $\gamma$ is the largest number such that 
some $a_i$, $b_i$, or $s_i$ goes to zero as $t^\gamma$.
\end{lemma}

\noindent{\em Proof}.   On the innermost ring, $\HH$ is not constant, and $K>0$, as shown in Lemma \ref{lemma:Hnotconstant}.
Suppose, for proof by contradiction, that there is no  $a_i$ on this ring.  Then $\AA$ is constant, $\AA = \alpha$.
 By Corollary~\ref{lemma:HwithoutS}, we have 
\begin{eqnarray*} 
\HH(\w) &=& \frac \BB \AA\, \mm\int_0^\w \AA^2\,d\w  - (\mm + \K)\int_0^\w \AA\BB\,d\w  \,  
\end{eqnarray*}
Substituting $\AA= \alpha$ we have 
\begin{eqnarray*} 
\HH(\w) &=&  \BB  \alpha \mm w  -  \alpha (\mm + \K)\int_0^\w \BB\,d\w  
\end{eqnarray*}
Looking at the asymptotic behavior for large $w$, we have $\BB = \beta w^R + O(w^{R-1})$, 
where $R$ is the degree of $\BB$, and $\beta$ is a nonzero constant.  Then
\begin{eqnarray*} 
\HH(\w) &=&  w^R  \alpha \beta\mm w  -  \alpha \beta(\mm + \K)  \frac {w^{R+1}}{R+1}  + O(w^R)\\
&=&\alpha\beta  w^{R+1}\bigg( \mm -  \frac {\mm + \K}{R+1} \bigg)  + O(w^R)
\end{eqnarray*}
We have $R > 0$ since there must be {\em some} root on the innermost ring, and we have supposed there are no $a_i$
and proved in Lemma \ref{lemma:NoS} that there are no $s_i$.  Hence, unless the coefficient is zero, 
$\HH$ will have a leading term in a power at least 2, and hence the eigenfunction, which is a multiple of $\Im\, \HH$,
will not have one sign in the upper half plane.  But it must have one sign.  Hence the coefficient is zero:
$$\mm = \frac {\mm + \K}{R+1} $$
Clearing the denominator we have 
\begin{eqnarray*}
\mm(R+1) &=& \mm + \K  
\end{eqnarray*}
Subtracting $\mm$ from both sides we have 
\begin{eqnarray}
\mm R &=& \K \label{eq:2197}
\end{eqnarray}
Now we recall the definitions of $\mm$ and $\K$ from (\ref{eq:Mdef}) and $\ref{eq:Kdef}$.  Since there   
is no $a_i$ on this ring, we have $Q = 0$ in those definitions; also since there no $s_i$ at all, we have 
 $S = 0$ and $f = 0$.   The equations for $\K$ and $\mm$  become
\begin{eqnarray*}
\K &=& R\gamma + (e-d) \\
\mm &=& 2d +\gamma
\end{eqnarray*}
Putting those expressions into \ref{eq:2197} we get
\begin{eqnarray*}
(2d + \gamma)R &=&  R\gamma + e-d  
\end{eqnarray*}
Subtracting $R\gamma$ from both sides we have 
\begin{eqnarray}
2d R &=& e-d  \nonumber \\
R &=& \frac {e-d}{2d} \label{eq:31}
\end{eqnarray}
Now recall what $e$ and $d$ are.  Let $\gamma$ be the exponent associated with the innermost ring, 
i.e., the largest number such that some root goes to zero as $t^\gamma$.  
 The number $d$ is, according to its definition in (\ref{eq:ddef}),
the product of one factor $t^{\gamma_j}$ for each $a_i$ that 
goes to zero on a slower ring, i.e., as $t^{\gamma_j}$ for $\gamma_j < \gamma$.
 Since we are working on the innermost ring and assuming that there are no $a_i$
on that ring,  that is one factor of $t^{\gamma_j}$ for {\em every} $a_i$.   Similarly, according to (\ref{eq:edef}),
 $e$ is a product of 
one factor $t^{\gamma_j}$ for each root $b_i$ on any but the innermost ring.  Now, on all the rings but the 
innermost, $\HH$ is constant, so on those rings, each   $a_i$ is matched by a $b_i$ converging to 
the same $\alpha_i$.  That means that for each factor $t^{\gamma_j}$ in the product defining $d$, the same 
factor occurs in the product defining $e$, and with at least as great a multiplicity,
since when $\HH$ is constant, $\BB$ is a polynomial multiple of $\AA$.
 Since all factors are powers of $t$, the more factors,
the smaller the number.   Hence $e \le d$.  Therefore $e-d \le 0$ and since $d > 0$, (\ref{eq:31}) gives us
$R \le 0$.  But as remarked above, $R >0$, so this is a contradiction.   That completes the proof of the lemma.

 \section{Finiteness}

Now we know that on the innermost  ring of roots, $\HH$ is not constant.
The rest of the proof will focus on that ring of roots only.  Henceforth,
$\gamma$ will be the number such that the roots on that ring go to zero as $t^\gamma$.
We have already established the following facts:
\begin{itemize}
\item The possible branch points $s_i$ in the lower half plane do not occur, so $\SS = 1$.
\item $\AA$ and $\BB$ do not vanish at the origin (i.e.,
there are no fast roots). 
\item The degree $R$ of $\BB$ and the degree $Q$ of $\AA$ are related by 
$$ R-Q = \frac \K \mm (2Q+1).$$
\item $Q > 0$, i.e. there is at least one $\alpha_i$.
\end{itemize}
 We define
\begin{equation}
\sigma := \int_0^w \AA^2 \, dw  \label{eq:sigmadef}
\end{equation}
Then $\sigma$ is a polynomial of degree $2Q+1$.
  As an example of the use of $\sigma$, 
the formula for $\HH$ in Corollary~\ref{lemma:HwithoutS} becomes
\begin{equation} \label{eq:HHdef2}
\HH(\w) = \frac \BB \AA\, \mm\sigma - (\mm + \K)\int_0^\w \AA\BB \,d\w    
\end{equation}
To assist the intuition: $\sigma$ is, except for a factor of a power of $t$, 
the first term in the arc length along the boundary $\Gamma$ from the origin.

A second example of the use of $\sigma$:  we can express the condition shown in
in Lemma~\ref{lemma:Hconstant} and Lemma~\ref{lemma:impliesHconstant} to be 
equivalent to ``$\HH$ is constant''  as 
\begin{equation}
 \frac \BB \AA = c \sigma^{\K/\mm} \label{eq:2528}
\end{equation}
We have dropped the factor $\SS$ occurring in the cited lemmas, because it is 
now known to be 1.

\subsection{A differential equation for $\BB/\AA$ and its solution}
Starting with the formula for $\HH$,  we will differentiate,  change variables,
and integrate.  The result is a formula that connects the (stereographic projection of the) 
normal $\BB/\AA$ with the derivative of the 
eigenfunction.  The formula involves a complex path integral.  The basic idea 
of this final section of the paper is to analyze the behavior of the terms in 
this formula for $\BB/\AA$, especially the complex integral, at zero,  at the complex
zeroes of $\sigma$, and finally at infinity.   Using the facts that $\HH$ must have 
just one sign in the upper half plane, and is not constant,  we are able to reach 
a contradiction.   Here is the key formula on which these results are based:

\begin{lemma} \label{lemma:exactB/A}
Where both sides are defined we have for some constant $c$
\begin{eqnarray*}
\frac {\BB}{\AA} &=& c\sigma^{\K/\mm} \exp \bigg(\int \frac {\HH_w}{\mm \sigma \BB /\AA} \, d\w\bigg)
\end{eqnarray*}
 The path of integration begins somewhere on the positive real axis,  at a point
 greater than all the $\alpha_i$ and $\beta_i$. 
For real $w$ we assume the path
of integration is on the real axis except for small semicircles in the upper half plane 
around the  poles of $\BB/\AA$.  
\end{lemma}
\medskip

\noindent{\em Remarks}.  If $\HH$ is constant, then $\HH_w = 0$, so the exponential factor is 1,
and the formula simplifies to the formula given in (\ref{eq:2528})
 for $\BB/\AA$ when $\HH$ is constant.  In some
sense, this new formula measures the deviation from what is true when $\HH$ is constant.  
The lemma mentions the zeroes of $\sigma$; of course 0 is the only real zero of $\sigma$, but 
it may have complex zeroes. Even for real $w$,  we must use a 
complex path of integration, that deviates from the real axis to make small detours in the 
upper half-plane around the  poles of $\BB/\AA$.  Then the constant $c$ depends only
on the starting point of the path of integration.
 
\medskip

\noindent{\em Proof}.
By Lemma~\ref{lemma:Hdif} (with $\SS = 1$)
we have
\begin{eqnarray}
\HH_w &=&\mm \int_0^w \AA^2 \, d\w \ \frac d { d\w} \bigg( \frac \BB \AA \bigg)  - \K\, \AA\BB 
\label{eq:Hw}
\end{eqnarray}
Recall that $\sigma = \int_0^w \AA^2 \, dw$, so $\sigma_w = \AA^2$.   Introducing $\sigma$, the 
previous equation becomes
\begin{eqnarray*}
\HH_w &=& \mm \sigma \frac d { d\w} \bigg( \frac \BB \AA \bigg) -\K \sigma_w \frac {\BB}{\AA}
\end{eqnarray*}
This can be considered as a differential equation satisfied by $f = (\BB/\AA)$:
\begin{eqnarray*}  
\HH_w = \mm \sigma f_w - \K\sigma_w f 
\end{eqnarray*}
with $w$ as the independent variable.  We now want to regard $\sigma$ as the independent
variable. By definition $\sigma = \int_0^w \AA^2\, dw$.  Hence
 $\sigma_w = \AA^2$.  Then $\sigma$ is a one-one function of $w$ at least locally, 
 away from the $\alpha_i$ (which are zeroes of $\AA$),  so $f(w)$ can be regarded
 (at least locally) as a function of $\sigma$.  We write $f_\sigma$ for the derivative
 of this function.
 By the chain rule $f_w = f_\sigma \sigma_w = \AA^2 f_\sigma$.
Then the previous equation becomes
\begin{eqnarray*}
\mm  \sigma f_\sigma \sigma_w&=&  \K \sigma_w f + \HH_w \\
\mm \AA^2 \sigma f_\sigma &=&  \K \AA^2 f + \HH_w  \mbox{\qquad since $\sigma_w = \AA^2$} 
\end{eqnarray*}
Dividing both sides by $\mm \AA^2 \sigma f$ we have
\begin{eqnarray*}
\frac{f_\sigma} f &=& \frac \K { \mm \sigma} + \frac{\HH_w} { \mm\AA^2 \sigma f} 
\end{eqnarray*}
Replacing $f$ by $\BB/\AA$ in the last term we have
\begin{eqnarray*}
\frac{f_\sigma} f &=& \frac \K { \mm \sigma} + \frac{\HH_w} { \mm \sigma \AA \BB } 
\end{eqnarray*}
We are going to integrate; this will be a complex path integral and we want to 
specify the path.  Fix a starting point $q$ on the positive real $w$-axis away from 0 and the $\alpha_i$.  (We take $q > 0$ so that $\sigma(q)$ will be positive.)
Let $P$ be any path from $q$ to $w$, avoiding 0, any complex zeroes of $\sigma$, and
the zeroes and poles of $\BB/\AA$; to avoid those zeroes and poles it suffices to 
avoid the $\alpha_i$ and $\beta_i$.    By $P\sigma$ we mean the path in the $\sigma$-plane defined by 
composing $P$ with $\sigma$.  
We integrate with respect to $\sigma$ along the path $P\sigma$
(introducing $\log c$ as a constant of integration):
\begin{eqnarray*}
\log f &=& \frac \K \mm \log \sigma + \int \frac {\HH_w}{\mm \sigma \AA \BB} \, d\sigma + \log c 
\end{eqnarray*}
where the constant $c$ is explicitly given by 
$$ \log c = \log f(q) -\log \sigma(q)  = \log \frac {\BB}{\AA \sigma} \mbox{\qquad evaluated at $w=q$}
$$
Since $\BB/\AA$ is finite and non-zero at $q$, and $\sigma$ is non-zero at $q$, the 
constant $\log c$ is well-defined;  in particular $c \neq 0$.  (If $f(q) < 0$ then both 
$\log c$ and $\log f$ have an imaginary part $i\pi$.)  We have
\begin{eqnarray*}
\log f &=& \log \bigg( c \sigma^{k/\mm}\bigg) + \int \frac {\HH_w}{\mm \sigma \AA \BB} \, d\sigma   \\
&=& \log \bigg( c \sigma^{\K/\mm}\bigg) + \int \frac {\HH_w}{\mm \sigma \AA \BB} \, \AA^2 \, d\w 
     \mbox{\qquad since $d\sigma = \AA^2\, d\w$} \\
 &=& \log \bigg( c \sigma^{\K/\mm}\bigg) + \int \frac {\HH_w}{\mm \sigma  \BB /\AA} \,  d\w     
\end{eqnarray*}
Here the path of integration must avoid the singularities 
of the integrand, which occur at $\sigma = 0$ (the only real zero of $\sigma$ is $w=0$, but there may be
complex zeroes), at zeroes of $\BB/\AA$, and possibly at singularities of $\HH_w$.   
The branch of the logarithm is chosen to be continuous along the path of integration.

    Applying the function $\exp(x) = e^x$ to both 
sides of the equation, we have
\begin{eqnarray*}
f &=& c\sigma^{\K/\mm} \exp \bigg(\int \frac {\HH_w}{\mm \sigma \BB/\AA} \, dw \bigg)
\end{eqnarray*}
Recalling that $f = \BB/\AA$ we have proved the lemma. 
\medskip

\subsection{$\K/\mm$ is an integer}
This is an important step in the proof.  Among other things,
it implies that $\sigma^{\K/\mm}$ is 
a polynomial.  We prove it by analyzing the behavior of the path integral 
in Lemma~\ref{lemma:exactB/A} at the origin.

\begin{lemma} \label{lemma:K/minteger}[This claim is withdrawn.  The last line of the proof is incorrect.]
$\K$ is a multiple of $\mm = 2m+1$; in other words $\K/\mm$ is an integer.
\end{lemma}

\noindent{\em Proof}.  First we calculate $\HH_w$ at the origin. 
We have by (\ref{eq:Hw})
\begin{eqnarray*}
\HH_w &=&\mm \sigma \ \frac d { d\w} \bigg( \frac \BB \AA \bigg)  - \K\, \AA\BB 
\end{eqnarray*}
At the origin, $\sigma$ is zero.  Let $\alpha$ and $\beta$ be the constants such that
  $\AA = \alpha \prod(w-\alpha_i)$ and $\BB = \beta \prod(w-\beta_i)$.  Let 
  $\eta$ be the product of all the $-\alpha_i$, and $\delta$ the product of all the 
  $-\beta_i$.  Then 
the value of $\AA\BB$ at the origin is $\alpha\beta\eta\delta$.  Hence
the value of $\HH_w$ at the origin is $-\K\alpha\beta \eta\delta \neq 0$.

By Lemma~\ref{lemma:exactB/A},  
\begin{eqnarray*}
\frac \BB \AA &=& c\sigma^{\K/\mm} \exp \bigg(\int \frac {\HH_w}{\mm \sigma \BB/\AA} \, dw \bigg)
\end{eqnarray*}
As $w$ moves from small positive values to small negative ones, the left side
remains real, and since neither $\BB$ nor $\AA$ is  zero at the origin, the left side does not change 
sign.    Since $\HH_w \neq 0$ at the origin, the integrand has a simple pole, due to the 
simple zero of $\sigma = \int_0^w \AA^2\, dw$, as we now show in detail. 
 Near the origin,
$\sigma = \int_0^w \AA^2\,dw = \alpha^2 \eta^2 w + O(w^2)$. Hence the 
denominator is 
\begin{eqnarray*}
 \mm\sigma\BB/\AA &=& \mm (\alpha^2\eta^2 w) \frac {\beta\delta}{\alpha\eta} + O(w^2) \\
 &=& \mm \alpha \eta \beta \delta w + O(w^2)
\end{eqnarray*}
Since the value of 
the numerator at the origin is $-\K\alpha\beta \eta\delta$,  
the residue is $-\K /\mm$.  In order that the right side remain real when 
$w$ passes from small negative to small positive values,  this residue must be 
an integer.   That completes the proof.
\medskip

\noindent{\em Remark}.  The factor $\sigma^{\K/\mm}$ changes sign if and only if $k$ is odd,
because $\K/\mm = k/(2m+1)$. 
Hence the exponential factor changes sign if and only if $k$ is odd.
For the exponential factor to change sign, the complex path integral in 
the exponent must pick up an odd multiple of $i\pi$ on a small 
semicircle in the upper half plane around the origin.  That also happens if and only 
if $\K/\mm$ is odd,  so in either case, the right side does not change signs.

\subsection{The non-real zeroes of $\sigma$ are also zeroes of $\BB$}
In this section, we analyze the behavior of the formula in Lemma~\ref{lemma:exactB/A}
at the non-real zeroes of $\sigma$.  The conclusion is that these zeroes are 
also zeroes of $\BB$, and moreover, they have the same multiplicity as zeroes of $\BB$
as they do as zeroes of $\sigma^{\K/\mm}$ (which is a polynomial because $\K/\mm$ is an 
integer).
\smallskip
 
We start with the following fact from complex analysis:

\begin{lemma} \label{lemma:complexpole}
Let $f(w)$ be a rational function with a pole at $q$.
Let $\gamma:[0,1] \rightarrow \C$
be a path in the complex plane avoiding the poles of $w$, except that 
at the last point of $\gamma$,  $\gamma$ reaches $q$, i.e. $\gamma(1) = q$.
Let $p = \gamma(0)$.  Then  
 $$\lim_{\tau \to 1} \int_p^{\gamma(\tau)} f(w)\, dw$$
does not have a value independent of $\gamma$; indeed, for appropriate choices 
of $\gamma$ it can be made $-\infty$, $\infty$, or 0.
\end{lemma}

\noindent{\em Proof}. By considering the new rational function $f(w-q)$,
 we can assume without loss of generality that the pole $q$ is 
at $w=0$ and the function $f$ has the asymptotic form
$$ f(w) = Cw^{-j} + O(w^{-j+1}) = Cw^{-j} (1 + O(w))$$
for some constant $C \neq 0$.  Also without loss of generality, we can assume that 
$C$ is real.  We seek a change of variable $w \mapsto \zeta$ such that as a 
function of $\zeta$ we have
$$ f(w) = C\zeta^{-j} \mbox{\qquad exactly} $$
Once we have such a variable $\zeta$,  the desired paths of integration in the $\zeta$ plane
are simply lines given by $ g(\tau) =(1-\tau)e^{i\theta}$ for various $\theta$.  Along 
such a path we have
\begin{eqnarray*}
\int f(w)\, dw &=&\int_0^{g(1-\tau)}\zeta^{-j} d\zeta \\
&=&  \int_0^{1-\tau} g(t)^{-j} g^\prime(t)\, dt \\
&=& \int_0^{1-\tau} (1-t)^{-j} e^{-ij\theta} g^\prime(t)\, dt\\
&=& e^{-ij\theta} \int_0^{1-\tau} (1-t)^{-j}  g^\prime(t)\, dt\\
&=& e^{-ij\theta} \int_0^{1-\tau} (1-t)^{-j}  (-e^{i\theta})\, dt\\
&=& -e^{-(j+1)i\theta} \int_0^{1-\tau} (1-t)^{-j}\, dt \\
&=&  e^{-(j+1)i\theta} \int_1^\tau x^{-j}\, dx 
\end{eqnarray*}
Since the pole at $q$ is at least a double pole, we have $j \neq 1$, and 
\begin{eqnarray*}
\int f(w)\, dw &=&  e^{(j-1)i\theta} \bigg(1 - \frac{ \tau^{-j+1}}{-j+1}\bigg) 
\end{eqnarray*}
and  choosing $\theta$ appropriately we can make the right side approach 
$\pm \infty$ or 0 as desired.  
If $j=1$, we have
\begin{eqnarray*}
\int f(w)\, dw &=&  e^{-2i\theta}\log \tau 
\end{eqnarray*}
which also can be made to approach $\pm \infty$ or 0 by an appropriate choice of $\theta$.

It remains to show how to find $\zeta$ as a function of $w$.  We have 
$$f(w)= Cw^j( 1+ g(w))$$
for some analytic function $g(w)$ with $g(w) = 0$.  Then by the binomial series, 
$(1+g(w))^{1/j}$ is   analytic  in some neighborhood of the origin.  That is, 
there is a function $h$ analytic in some neighborhood of the origin such that 
$$(h(w))^j = 1 + g(w)$$
Now define 
$\zeta = w h(w).$
Then $\zeta^j = w^j (h(w))^j = w^j(1+g(w)) = f(w)$, so we have found the required 
change of variables.  That completes the proof of the lemma.
\medskip

Next we extract information by looking at the formula of Lemma~\ref{lemma:exactB/A} 
at the  non-real zeroes of $\sigma$.   Those zeroes occur in complex-conjugate pairs,
but our argument does not depend on which half-plane contains the zero in question.
\smallskip

\begin{lemma} \label{lemma:nonrealzeroes}
 Every non-real zero of $\sigma$ is also a zero of 
$\BB$, and as a zero of $\BB$ it has multiplicity $\K/\mm$ times greater
than its multiplicity as a zero of $\sigma$.
\end{lemma}

\noindent{\em Proof}. 
By (\ref{eq:Hw}) we have 
\begin{eqnarray*}
\HH_w &=&\mm \sigma \ \frac d { d\w} \bigg( \frac \BB \AA \bigg)  - \K \AA\BB.
\end{eqnarray*}
Suppose that $q$ is a non-real zero of $\sigma$.     Since the zeroes of
$\AA$ are all real,  $\BB/\AA$ is analytic at $q$.  Then at $q$ we have
$$ \HH_w = -\K\AA\BB $$
Therefore 
$$\frac {\HH_w} {\BB/\AA} = -\K\AA^2 \neq 0.$$
Therefore
$$ \frac {\HH_w \,dw}{\mm \sigma \BB/\AA} \mbox{\qquad has a pole at $q$} $$
 By Lemma~\ref{lemma:exactB/A}, we have
\begin{eqnarray*}
\frac \BB \AA &=& c\sigma^{\K/\mm} \exp \bigg(\int  \frac {\HH_w \,dw}{\mm \sigma \BB/\AA}  \bigg).  
\end{eqnarray*}
Now assume, for proof by contradiction, that $\BB$ is not zero at $q$. Then
the left side, $\BB/\AA$,  is not zero at $q$ and is continuous there, since
$q$ is not real but the zeroes of $\AA$ are real. 
On the right side, $\sigma^{k/\mm}$ is zero at $q$.  Therefore the exponential factor
has to tend to infinity as $w$ approaches $q$.  That means that the integral 
also tends to infinity as $w$ approaches $q$.  But by Lemma~\ref{lemma:complexpole},
 the integral does not 
have a limit as $w$ approaches $q$; choosing an appropriate path approaching
$q$ the real part of the (absolute value of the) integrand will remain bounded.   Hence the 
exponential factor remains bounded away from both zero and infinity.
 Hence the right side of the equation goes to zero. 
That is a contradiction, since the left side approaches a nonzero value.
Hence in fact $\BB$ is zero at $q$.

It remains to prove that $q$ has the same multiplicity as a zero of $\B$ as 
it does as a zero of $\sigma^{\K/\mm}$.  Suppose, for proof by contradiction,
that it has a different multiplicity.  Because the exponential factor has 
a non-zero finite limit (for a suitable approach to $q$),  along that 
approach $\BB/\AA$ and $\sigma^{\K/\mm}$ must approach zero at the same rate.
Since $\AA$ is not zero at $q$, it follows that $q$ has the same multiplicity
as a zero of $\BB$ and as a zero of $\sigma^{\K/\mm}$.  That completes the proof
of the lemma.

\subsection{Laurent expansions at infinity}
Recall that $R$ is the degree of $\BB$, $Q$ is the degree of $\AA$, 
and 
$$ R-Q = \frac \K \mm (2Q+1).$$
Looking at the formula of Lemma~\ref{lemma:exactB/A}, we see that the 
left side $\BB/\AA$ and the factor $\sigma^{\/\mm}$ on the right  both 
are asymptotic to $w^{R-Q}$ for large $k$.  The ratio of the two therefore
goes to 1 as $w$ goes to infinity.  But what is the next power in its 
Laurent expansion?  We need to show that it can't be too very far out.
Our answer is that it is no farther out than the term in $w^{-(R-Q)}$. 
It is important that $\sigma^{\K/\mm}$ is a polynomial--otherwise we would 
have no hope of answering this question!  Therefore it is natural to 
begin with a lemma about the Laurent expansions of rational functions.
\smallskip

\begin{lemma} \label{lemma:Laurent} Let $f$ and $g$ be two unequal
rational functions of the same degree $n$
with the same leading coefficient.  Let $n-J$ be the highest power of $w$
such that the coefficients of $w^J$ in $f$ and $g$ are different.  Then
the Laurent expansion of $f/g$ at infinity is given by 
$$ \frac {f(w)}{g(w)} = 1 + \mu w^{-J} + O(w^{-J-1})$$
for some constant $\mu \neq 0$.  In particular $J \le n$.
\end{lemma}

\noindent{\em Proof}.
   The worst case (largest $J$) occurs when the numerator 
and denominator have their corresponding coefficients equal.  Here is 
an example:
\begin{eqnarray*}
 \frac {w^n + 3}{w^n + 5} &=& \frac { 1 + 3w^{-k}}{1 + 5w^{-k}} \\
   &=& (1 + 3w^{-k})(1-5w^{-k})(1 + O(w^{-1}))\\
   &=& (1-2w^{-k})(1+O(w^{-1}))
\end{eqnarray*}
The general proof is only notationally more complicated: Let
$f(w) = a_0w^n + \ldots +a_n$ and $g(w) = b_0w^n + \ldots + b_n$,
with $a_0 = b_0$.
Let $J$ be 
the smallest integer such that $a_J \neq b_J$.
Then with $h(w)$ defined as the sum of the 
terms with lower-indexed coefficients  we have 
\begin{eqnarray*}
\frac {\sum a_i w^{n-i}} {\sum b_i w^{n-i}} 
  &=& \frac {h(w) + \sum_{i=J}^n a_i w^{n-i}} 
          {h(w) + \sum_{i=J}^n b_i w^{n-i}} 
\end{eqnarray*}
Dividing numerator and denominator by $h(w)$ we have 
\begin{eqnarray*}
\frac {\sum a_i w^{n-i} }{\sum b_i w^{n-i}}
&=& \frac {1 + (a_J/a_0) w^{-J} + O(w^{-J-1})}
          { 1 +(b_J/b_0) w^{-J} + O(w^{-J-1})} \\
&=& 1 + \frac {a_J} {a_0} - \frac {b_J}{b_0} w^J + O(w^{-J-1}) \\
&=& 1 + \frac {a_J-b_J}{a_0}w^{-J} + O(w^{-J-1}) \mbox{\qquad since $a_0 = b_0$}
\end{eqnarray*}
Since by hypothesis $a_J \neq b_J$, we have found the first nonzero 
term in the Laurent expansion. That completes the proof of the lemma.

 \begin{lemma} \label{lemma:Laurent2}  For some constant $\mu \neq 0$
and some integer $J$ with $1 \le J \le R-2Q$ we have the Laurent 
expansion at infinity:
$$ \sigma^{-k/\mm} \bigg( \frac \BB \AA \bigg) = 1 + \mu w^{-J} + O(w^{-J-1}).$$
\end{lemma}

\noindent{\em Remark}.  Whether or not $\K/\mm$ is an integer, we always 
have a Laurent expansion.  The point is that since $\K/\mm$ is an 
integer, we have a rational function on the left and can therefore
say where the first nonzero term in the Laurent expansion must occur,
while without knowing that $\K/\mm$ is an integer, the first nonzero 
term could occur very far out in the Laurent series.
\medskip

\noindent{\em Proof}.  Since $\K/\mm$ is an integer, the function
$\sigma^{\K/\mm}$ is a polynomial.  Since the degree of $\sigma$ is 
$2Q+1$, the degree of $\sigma^{\K/\mm}$ is $(\K/\mm)(2Q+1) = R-Q$.   The
 expression in the lemma can be written as 
 $$   \frac {\BB} {\AA\sigma^{\K/\mm}}.$$
 This is a rational function.  It is not constant, since 
 if it were constant, then by Lemma~\ref{lemma:impliesHconstant},  $\HH$
 would be constant, but $\HH$ is not constant.
  On the face of it, the numerator has 
 degree $R$, and the denominator has degree $Q + (\K/\mm)(2Q+1) = Q + (R-Q) = R$, the
 same as the numerator.  But appearances are deceptive:
 by Lemma~\ref{lemma:nonrealzeroes}, the non-real zeroes of the denominator
 are all zeroes of the numerator with the same multiplicity.  Now $\sigma$
 is of degree $2Q+1$ and has only one real zero (namely zero), so there 
 are $2Q(\K/\mm)$ of these zeroes, counting multiplicities.  After these
 common zeroes are cancelled out, the degrees of numerator and denominator 
 will be 
 \begin{eqnarray*}
 n &:=&  R - \frac {2Q\K} \mm  
 \end{eqnarray*}
 Since $\K/\mm$ is an integer, the last term being subtracted off is at least $2Q$. Hence
$$ n \le R-2Q.$$

  Applying Lemma~\ref{lemma:Laurent} to this rational function,
there exists a term $w^{-J}$ in the Laurent expansion with $J \le R-2Q$.
That completes the proof of the lemma.

\subsection{The final contradiction}
Now we are ready to analyze the behavior of the formula of Lemma~\ref{lemma:exactB/A}
at infinity. The numerator of the integrand is $\HH_w$.  
On the face of it, by (\ref{eq:Hw}), $\HH_w$ is
a rational function asympotic to $w^{R-Q}$.  One can easily calculate that the leading
coefficient vanishes, but since $\HH$ has only one sign in the upper half plane, we must 
have $\HH$ asymptotic to $w$ or $w^{-1}$.  That means that the integrand goes to zero 
as a large negative power of $w$,  as those higher powers of $\HH$ are not there to cancel
the denominator.  That makes the integral (to infinity) finite; that is the first 
key point of the proof.  Then we divide the formula of Lemma~\ref{lemma:exactB/A}
by $\sigma^{k/\mm}$,  isolating the exponential factor on the right:
\begin{eqnarray*}
\sigma^{-\K/\mm} \frac \BB \AA &=& c  \exp \bigg(\int  \frac {\HH_w \,dw}{\mm \sigma \BB/\AA}  \bigg).  
\end{eqnarray*}
Now we have worked out the Laurent expansion of the left side, and the order 
of the pole in the integrand.  Differentiating this equation, we arrive at a comparison
between the first non-zero term of the Laurent expansion, and the order of the pole.
The   first  term  on the left is $w^J$ (or $w^{J-1}$ after differentiating),
 where $J$ is at most so much, 
and the order of the pole on the right is at least so much,
and it turns out that the two estimates are incompatible.  The details are below.
\smallskip
  
\begin{theorem}[Finiteness] \label{theorem:main1}
[This claim is withdrawn. The proof is not correct.]
Let $\Gamma$ be a real-analytic Jordan curve in $R^3$.  Then $\Gamma$ cannot bound
  infinitely many immersed disk-type minimal surfaces whose least eigenvalue
satisfies $\lambda_{\min} \ge 2$.
\end{theorem}

\noindent{\em Proof}.
By Theorem \ref{theorem:one-parameter},
 $\Gamma$   bounds a one-parameter family of minimal surfaces $u^t$ with $u^t$ immersed for $t > 0$, and
each $u^t$ having least eigenvalue 2, so the assumptions that have been in force for all 
our calculations are satisfied.  By Lemma~\ref{lemma:exactB/A}, we have
\begin{eqnarray*}
\frac \BB \AA &=& c\sigma^{\K/\mm} \exp \bigg(\int  \frac {\HH_w \,dw}{\mm \sigma \BB/\AA}  \bigg).  
\end{eqnarray*}
Dividing both sides by $\sigma^{\K/\mm}$ we have
\begin{eqnarray*}
\sigma^{-\K/\mm} \frac \BB \AA &=& c  \exp \bigg(\int  \frac {\HH_w \,dw}{\mm \sigma \BB/\AA}  \bigg).  
\end{eqnarray*}
 We  consider the behavior
for large $w$.  By Lemma~\ref{lemma:Laurent2}, there is an integer $J$ with
$1 \le J \le R-2Q$ such that the Laurent expansion of the left side has
a term in $w^{-J}$.  The leading term is $\beta$, the leading coefficient of $\BB$.
Then we have for some $\mu \neq 0$ 
\begin{eqnarray}
\beta + \mu w^{-J} + O(w^{-J-1}) &=& c \, \exp \bigg(\int  \frac {\HH_w \,dw}{\mm \sigma \BB/\AA}  \bigg).  
\label{eq:2728}
\end{eqnarray}
Since $\Im\, \HH$ has one sign in the upper half plane, the Laurent expansion 
of $\HH$ has to begin with $w$ or $w^{-1}$.  Hence the Laurent expansion of 
$H_w$ begins with a constant term or with $w^{-2}$.
Therefore the first term of $H_w$ can be written $\nu w^{-1\pm 1}$ for
some constant $\nu \neq 0$. The integrand can then be written like this:
\begin{eqnarray*}
\frac {\HH_w}{\mm\sigma\BB/\AA} &=& \frac{\AA \HH_w} {\mm \sigma \BB} \\
&=& \frac {\alpha \nu w^{Q-1\pm1}} {\beta w^{2Q+1}  w^{R}} (1 + O(w^{-1}))\\
&=& \frac {\alpha \nu } \beta w^{-1\pm 1 - Q -R-1} (1 + O(w^{-1}))\\
&=& \frac {\alpha \nu} \beta w^{-(2+Q+R \pm 1)}(1 + O(w^{-1}))
\end{eqnarray*}
Summarizing that calculation,
\begin{eqnarray}
\frac {\HH_w}{\mm\sigma\BB/\AA}  &=&\frac {\alpha\nu} \beta w^{-(2+Q+R\pm 1)} (1 + O(w^{-1}))\label{eq:2751}
\end{eqnarray}
Since $Q$ and $R$ are the degrees of $\AA$ and $\BB$, respectively, we have $Q+R \ge 1$.
Therefore the exponent is negative and has magnitude at least $Q+R+1$, which is at 
least 2.   It follows that the integral has a finite limit as $w$ goes to infinity.
Let that finite limit be $L$.  Then the integral is $L + O(w^{-1})$.

Differentiating (\ref{eq:2728}), which is legal since the $O(w^{-J-1})$ term stands
for a Laurent series,  we have
\begin{eqnarray*}
-J\mu w^{-J-1} + O(w^{-J-2}) &=& 
c \, \exp \bigg(\int  \frac {\HH_w \,dw}{\mm \sigma \BB/\AA}  \bigg)
\frac {\HH_w }{\mm \sigma \BB/\AA} \\
  &=& c \, \exp \bigg(L + O(w^{-1})  \bigg)
\frac {\HH_w }{\mm \sigma \BB/\AA} \\
&=& ce^L (1 + O(w^{-1}))\frac {\HH_w }{\mm \sigma \BB/\AA} 
\end{eqnarray*}
Putting in the result from (\ref{eq:2751}) we have
\begin{eqnarray*}
-J\mu w^{-J-1} + O(w^{-J-2}) 
&=& ce^L (1 + O(w^{-1}))\frac {\alpha\nu} \beta w^{-(2+Q+R\pm 1)} (1 + O(w^{-1}))
\end{eqnarray*}
On the left, we have a negative exponent $J+1 \le R-2Q+1$ (by Lemma~\ref{lemma:Laurent2}).
On the right, we have a negative exponent of magnitude at least $R+Q+1$.
If these two exponents were equal, we would have $R-2Q+1 \ge R+Q+1$. Subtracting
the right side from both sides we would have $-Q \ge 0$.  Since $Q$ is the degree of $\AA$,
which is non-negative, we would have $Q=0$.  But that contradicts Lemma~\ref{lemma:someA},
so we have reached a contradiction.
 
That contradiction depends on no other assumptions than the existence of a one-parameter
family of minimal surfaces $u^t$ bounded by $\Gamma$, such that $\lambda_{\min} = 2$ for $t > 0$.
That completes the proof by contradiction that no such family of minimal surfaces exists.

\section{Relative minima in which topology?}

We are going to prove a theorem limiting the number of relative minima of area with a given boundary.
To state such a theorem precisely, we must fix a topology to use in defining ``relative minimum.''
Once a topology is fixed, a relative minimum is a surface $u$ such that there exists a neighborhood
of $u$ (in the chosen topology) such that no surface in that neighborhood has smaller area than $u$.
As usual in Plateau's problem, we work with harmonic surfaces.  

We wish to use the $C^n$ topology. First consider whether it is a stronger theorem for smaller $n$
or larger $n$.   {\em A priori} it might be possible to decrease area by a $C^0$ variation but not by
 a $C^n$ variation, but if area can be decreased by a $C^n$ variation, that same variation counts as 
a $C^0$ variation.  Therefore a $C^0$ relative minimum is automatically a $C^n$ relative minimum, but
not vice-versa.  Thus if there were infinitely many $C^0$ relative minima, there would automatically
be infinitely many $C^n$ relative minima.  Taking the contrapositive, if there cannot be infinitely
many $C^n$ relative minima, there cannot be infinitely many $C^0$ relative minima.  Thus the theorem will be stronger,
when we use the $C^n$ topology for as large an $n$ as possible.  

Unfortunately, our result depends on the regularity result that a relative minimum cannot have an 
interior or boundary branch point,  which at present is known only for the $C^0$ topology.%
 
When and if this regularity result is proved for larger $n$, our proof will immediately apply to the $C^n$
topology; but for now we can prove it for the $C^0$ case.

 \begin{theorem} [Main theorem] \label{theorem:main2}[This claim is withdrawn.  The proof is not correct.]
  Let $\Gamma$ be a real-analytic Jordan curve in $R^3$. 
  Suppose that $n$ is a positive integer and every relative minimum of area 
in the $C^n$ topology bounded by $\Gamma$ has no (interior or boundary) branch points.
Then $\Gamma$ cannot bound infinitely many disk-type minimal surfaces  that are relative 
minima of area in the $C^n$ topology.   
\end{theorem}
\noindent
{\em Remark\/}:   As usual in Plateau's problem, we are counting two surfaces
that differ only by a conformal reparametrization as the same surface.  

\noindent
{\em Proof\/}:  Suppose, for proof by contradiction, that $\Gamma$ bounds
infinitely many  disk-type minimal surfaces furnishing $C^n$ relative minima of 
area.  Then those surfaces are immersed, by hypothesis,  and have $\lambda_{\min} \ge 2$. 
 Then by Theorem \ref{theorem:one-parameter},
 $\Gamma$   bounds a one-parameter family of minimal surfaces $u^t$ with $u^t$ immersed for $t > 0$, and
each $u^t$ having least eigenvalue 2.  But that contradicts Theorem~\ref{theorem:main1}.
That completes the proof.
\medskip

\begin{corollary}   Let $\Gamma$ be a real-analytic Jordan curve in $R^3$. 
Then $\Gamma$ cannot bound infinitely many disk-type minimal surfaces  that are relative 
minima of area in the $C^0$ topology. 
\end{corollary} 

\noindent{\em Proof}.  By Theorem~\ref{theorem:regularity} and the known regularity 
results, which are fully discussed in   
 section~\ref{section:regularity} above. 
 
 \section{Open Problems and Conjectures}

The theorem proved here could possibly be generalized by weakening any of these hypotheses:  (i) that the 
boundary is real-analytic; (ii) that the surfaces considered are of the topological type of the disk; or  (iii)  that the surfaces considered
are relative minima of area.   All questions formed by weakening one or more of these hypotheses are open. 
These questions are all ``equal conclusion'' strengthenings of the theorem, formed by weakening the hypotheses.
We will comment on these questions one by one below.

There is also an interesting ``stronger conclusion'' version of the theorem, in which we ask for an explicit bound on the number of relative minima as a function of the geometry of the boundary curve 
$\Gamma$.  Our ignorance is scandalous:  Not counting curves for which uniqueness is known, there is only one Jordan curve $\Gamma$  known for which anyone can give an explicit bound on the number of solutions of Plateau's problem  (of disk type) bounded by $\Gamma$. That curve is  ``Enneper's wire'' $\Gamma_r$ (the image of the circle of radius $r$ in the usual parametrization of Enneper's surface), for values of $r$ slightly larger
than $1$, and the bound is three.  This bound follows from the theorem of Ruchert \cite{ruchert} that uniqueness holds up to $r=1$ and 
the theorem of \cite{beeson-tromba} characterizing the space of minimal surfaces bounded by $\Gamma_r$ in the vicinity of Enneper's surface for $r=1$.
The following
conjecture, for example, seems presently out of reach:  if the real-analytic Jordan curve  $\Gamma$ has total curvature less than or equal to $6\pi$, then it bounds at most two relative minima of area (of disk type).  We do know that there can 
be no such estimate for the number of minimal surfaces (minimizing or not),  since B\"ohme proved in \cite{bohme3} that 
for any $\epsilon > 0$ and any $n$, there is a real-analytic Jordan curve with total curvature less than $4\pi + \epsilon$
bounding more than $n$ minimal surfaces.   But these are branched surfaces, not relative minima,
so they do not show the impossibility of a bound on the number of immersed minimal surfaces
bounded by $\Gamma$.  

Since Enneper's wire bounds two relative minima when the total curvature exceeds $4\pi$, we could expect to string together 
sequences of Enneper's wires with small bridges connecting them,  obtaining for each $n$ a Jordan curve with 
total curvature $4n\pi + \epsilon$ bounding (at least) $2^n$ relative minima (as well as some unstable minimal surfaces).  
Since we seem to have to increase the total boundary curvature to get many solutions of Plateau's problem, we might hope that for some explicitly given 
function $F$,  we might have,  for all real-analytic Jordan curves 
$\Gamma$ of total curvature less than $2 \pi n$,  there are at most $F(n)$ relative minima of area.   As far as we 
know,  this might be true with $F(n) = 2^{n - 1 }$.  Thus:  (when $n=2$) curvature less than $4\pi$ implies uniqueness (which is known); when $n=3$, curvature less than $6\pi$ implies at most two relative minima (the conjecture mentioned above);
when $n=4$, curvature less than $8\pi$ implies at most 4 relative minima  (also, of course, not known). 

Next we take up the possible ``equal conclusion'' versions.   First,  what about weakening the smoothness of the 
boundary to $C^n$?  Then we have some troubles about the boundary branch point representation, but perhaps that can 
be solved adequately.  A more pressing difficulty is that we then cannot control the dependence of the zeroes, poles, and 
ramification points of the Gauss map $N$ on the parameter $t$.  There is also a third difficulty:
even the proof of Tomi's theorem that there can't be infinitely many absolute minima does not go through for $C^\infty$ 
boundaries, because the regularity result that a minimum of area cannot have a boundary branch point is still an open problem for $C^\infty$ boundaries.    

Second, what about generalizing the theorem to other topological types of surfaces or boundaries? Frank Morgan has 
given an example \cite{morgan81} of a boundary consisting of four circles, such that it bounds a continuum of unstable
minimal surfaces of arbitrarily high genus.  That still leaves open the finiteness problem for relative minima.  I 
conjecture that a system of  real analytic Jordan curves cannot bound infinitely many relative minima of any fixed orientable topological type (that is, any fixed number of boundary components and Euler characteristic).  
For non-orientable surfaces, there is a problem in that the least eigenfunction 
no longer has one sign.  Instead, it has just one nodal line.  Thus even Tomi's argument that there cannot be a loop of 
minimal surfaces all of which are absolute minima of area fails to go through for surfaces of the type of the M\"obius 
strip.  That is the simplest open question for non-orientable surfaces:  can there be a one-parameter family of absolute
minima of surfaces of the type of the M\"obius strip?   

Incidentally, another paper of Morgan \cite{morgan76} gives an example of a curve in $R^4$ that bounds a continuum of 
absolute minima of disk type.  Therefore we cannot generalize the theorem to higher dimensional spaces.

Third,  what about dropping the requirement that the surfaces be relative minima of area?   Then our methods say nothing 
at all, since they all depend on analysis of the least eigenvalue and its eigenfunction.  
The only known results (other than uniqueness theorems) along these lines are two 
``$6 \pi$ theorems'', both of which assume that $\Gamma$ is a real-analytic Jordan curve whose
total curvature is less than or equal to $6 \pi$.  In \cite{beeson-sixpi}, it is proved that such a curve
cannot bound infinitely many minimal surfaces (stable or unstable).  Earlier, Nitsche proved that same conclusion
(see \textsection A29 of \cite{nitsche}, p. 447) using the additional hypothesis that $\Gamma$ 
bounds no minimal surface with branch points.
At the time of Nitsche's
proof,  the forced Jacobi fields had not yet been discovered.  In \cite{beeson-sixpi},  the theory of the forced Jacobi fields
is used.

\section*{Appendix: Dictionary of Notation}
To assist the reader we provide in one place a summary of the notations used in more than one section of this paper.  
This summary is provided as a memory aid, not as a list of complete definitions.  The symbols are more or less in 
order of introduction; there is no natural ``lexicographical'' order for such a variety of symbols.

\begin{tabbing}
principal roots \= \kill \\
$X,Y,Z$ \>   coordinates in $R^3$ \\
$\Gamma$   \> a real-analytic Jordan curve in $R^3$, \\
\>  tangent to the $X$-axis at the origin \\
$u = u^t$ \>  a one-parameter family of minimal surfaces, \\
\> analytic in $t$ in some interval $[0,t_0]$ \\
\> $u^t$ is analytic as a function of $z$ in a closed upper half-disk \\
$t$    \>    the parameter defining that one-parameter family \\
${}^2u$ \>   the second component of the vector $u = ({}^1u, {}^2u,{}^3u)$. \\
$z = x + iy$    \>    the parameter used to define $u^t$ \\
$N$    \>    the unit normal to $u^t$ \\
$g$    \>    the stereographic projection of the Gauss map $N$ \\
$u_t$   \>   subscript denotes differentiation:  $\frac {\partial u} {\partial t}$  \\
$\lambda_{\min}$   \>  least eigenvalue  \\ 
$\phi$  \>   $u_t \cdot N$.   If nonzero, this is an eigenfunction for eigenvalue 2. \\
$n$ and $\chi$ \>   $\phi = t^n \chi$. \\
$f$    \>    ${}^1u_z - {}^2u_z$,  occurs in the Weierstrass representation of $u^t$ \\
$2m$   \>    order of the boundary branch point of $u$ (when $t=0$). \\
$s_i$  \>    branch points in the lower half plane for $t>0$ that \\
\>converge to 0 as $t \to 0$.  There are $S$ of these.\\
$N$   \>    There are $2N$ of the $s_i$.  Hopefully no confusion will arise \\
\>with $N$ the unit normal.\\
$a_i$ \>  zeroes of $f$.  The $a_i$ depend analytically on $t$.  \\
\>There are $Q$ of these.  \\
$b_i$ \>  zeroes of $fg^2$.  The $b_i$ depend analytically on $t$. \\
\>There are $R$ of these. \\
$\gamma_n$   \>  The $n$-th ring of roots consists of $a_i$, $b_i$, and $s_i$ that go to zero as $t^{\gamma_n}$. \\
$\gamma$   \>  One of the $\gamma_n$.  \\
$\alpha_i$  \>  coefficient such that $a_i = \alpha_i t^\gamma + O(t^{\gamma+1})$ \\
$\beta_i$  \>  coefficient such that $b_i = \beta_i t^\gamma + O(t^{\gamma+1})$ \\
$\zeta_i$  \>  coefficient such that $s_i = \zeta_i t^\gamma + O(t^{\gamma+1})$ \\
\>   $\alpha_i$, $\beta_i$, and $\zeta_i$ are used only when non-zero, i.e. \\
\>for the roots that go to zero as $t^\gamma$. \\
$\alpha$\> product of all the $-\alpha_i$ over roots going to zero slower than $t^\gamma$  \\ 
$\beta$\> product of all the $-\beta_i$ over roots going to zero slower than $t^\gamma$ times $B_0(0,0)$ \\ 
$\zeta$\> product of all the $-\zeta_i$ over roots going to zero slower than $t^\gamma$  \\ 
$A $ \> $A = A_0 \prod_{i=1}^{m-N} (z-a_i)$    so $f = A^2$, and $A_0 = 1$ when $t=0$  \\
$B$ \> $B = \beta \prod_{i=1}^{m+k-N} (z-b_i)$ so $fg^2 = B^2$. \\
$C$  \>   complex constant, value of $B_0$ when $t=0$ \\
$S$ \> $S = \prod_{i=1}^{2N} (z-s_i(t))$\\
$w = z/t^{\gamma}$ \>  a ``blow-up'' of the parameter domain \\
$O(t)$   \>   refers to uniform convergence on compact subsets of the  \\
\> $w$-plane away from the $\alpha_i$ \\
$\AA $ \> $\AA = \alpha \prod_{i=1}^{Q}(w-\alpha_i)$  so $A = t^{(m-N)\gamma} \AA (1+\O(t))$ \\
\>   $\AA$ does not depend on $t$ \\
$\BB $ \> $\BB =  \beta  \prod_{i=1}^{R} (w-\beta_i)$   so $B = t^{(m+k-N)\gamma} \BB (1+\O(t))$\\
\> $\BB$ does not depend on $t$, and $B_0$ turns out to be real. \\
$\tAA $ \> $\alpha\prod (\w - \a_i)$.  The product is over $i$ such that $a_i = O(t^\gamma)$\\
$\tBB $ \> $ i\lambda B_0\beta\prod (\w - \b_i)$ \\
$\tSS$  \> $\zeta\prod (\w - \s_i)$ \\
$ \a_i$ \> $\a_i =  a_i/t^\gamma$ \\
$ \b_i$ \> $\b_i = b_i/t^\gamma$ \\
$ \s_i$ \> $\s_i = s_i/t^\gamma$ \\
$\K$ \>  number such that $t^\K \BB/\AA = B/A + O(t^{\K+1})$ \\
$\mm$ \> number such that $t^\mm \int \AA^2 \SS \, dw = \int A^2 S \, dz + O(t^{\mm+1})$ \\ 
$\tAA_0(t,w)$ \> $\tAA_0(t,w) = A_0(t,t^\gamma w)$ \\
$\tBB_0(t,w)$ \> $\tBB_0(t,w) = B_0(t,t^\gamma w)$ \\
$\chi$ \>  $\chi(w)  = \phi(t^{\gamma}w)\ /\ \gamma t^{(2m+k+1)\gamma -1} $\\
\>   $\chi$ shows what happens to the eigenfunction in the $w$-plane.   \\
$\chi^0$ \>  limit of $\chi$ when $t$ goes to $0$ \\
$\HH$ \>  complex analytic function such that $\chi^0 = -\Im \,\HH$ \\
$\sigma$ \> abbreviation for $\int_0^w \tAA^2 \, dw$ in the last two sections \\
$\tau$ \>  the parameter used to parametrize the boundary curve, equal\\
\> to the $X$-coordinate \\
\end{tabbing}

\end{document}